\documentclass[12pt,a4paper]{amsart}
\newif\ifpdf
\ifx\pdfoutput\undefined
	\pdffalse
\else
	\pdfoutput=1
	\pdftrue
\fi

\usepackage{enumerate}
\usepackage{amssymb}
\ifpdf
	\usepackage[pdftex]{graphicx}
	\DeclareGraphicsExtensions{.pdf}
	\usepackage{hyperref}
\else
	\usepackage{graphicx}
	\DeclareGraphicsExtensions{.eps}
	\usepackage{hyperref}
\fi

 \newcommand\bg{{{\mathbf G}_{m,k}}}
\newcommand\bgd{{{\mathbf G}_{m,k}^d}}
\newcommand\bge{{{\mathbf G}_{m,k}^{d+1}}}
\newcommand\bc{{\mathbb C}} \newcommand\ba{{\mathbb A}}
\newcommand\bp{{\mathbb P}} \newcommand\bn{{\mathbb N}}
 \newcommand\bz{{\mathbb Z}}
\newcommand\bq{{\mathbb Q}} \newcommand\ff{{\mathbb F}}
\newcommand\br{{\mathbb R}} 
\newcommand\mm{{\mathcal M}}  
 \newcommand\cl{{\mathcal L}}
 
\newcommand\ww{{\mathcal W}}
\newcommand\tww{{\widetilde{\mathcal W}}}

 \newcommand\bl{{\mathbb L}}
\newcommand\bdn{{\mathbf n}} \newcommand\bdk{{\mathbf k}}
\newcommand\bda{{\mathbf a}} \newcommand\bdx{{\mathbf x}}
\newcommand\bdg{{\mathbf g}} \newcommand\bdz{{\mathbf z}}
\newcommand\bdw{{\mathbf w}}  \newcommand\bdr{{\mathbf r}}
\newcommand\bds{{\mathbf s}}   \newcommand\bdy{{\mathbf y}}
    \newcommand\bdh{{\mathbf h}}
     \newcommand\bdB{{\mathbf B}}
 
\newcommand\barphi{\boldsymbol{\varphi}}\newcommand\bdm{{\mathbf m}}
\newcommand\bpsi{\boldsymbol{\psi}}
\newcommand\blambda{\boldsymbol{\lambda}}
\newcommand\bmu{\boldsymbol{\mu}}
\newcommand\balpha{\boldsymbol{\alpha}}
\newcommand\bbeta{\boldsymbol{\beta}}

\newcommand\zlo{{Z_{\text{\rm top},0}}}

\newcommand\ord{{\operatorname{ord}}}
\newcommand\supp{{\operatorname{supp}}}

\DeclareMathOperator\vol{Vol}
\DeclareMathOperator\dpt{depth}
\DeclareMathOperator\lgt{\mathbf{ev}}
\numberwithin{equation}{section}

\DeclareMathOperator\mult{mult}
\DeclareMathOperator\spec{Spec}
 \DeclareMathOperator\mor{Mor}

\newtheorem{thm}{Theorem}[section]
\newtheorem{lema}[thm]{Lemma}

\newtheorem{cor}[thm]{Corollary}
\newtheorem{spc}[thm]{Support Condition}
\newtheorem{prop}[thm]{Proposition}
\newtheorem*{mmc}{Motivic Monodromy Conjecture}
\newtheorem*{imc}{Igusa Monodromy Conjecture}
\newtheorem*{tmc}{Topological Monodromy Conjecture}

\theoremstyle{remark}
\newtheorem{obs}[thm]{Remark}

\theoremstyle{definition}
\newtheorem{defini}[thm]{Definition}

\newtheorem{ejem}[thm]{Example}
\newtheorem{paso0}{Step}
\newtheorem{paso}{Step}

\newtheorem*{algthm}{Algorithm}
\newtheorem*{ctrb}{Contribution of the $B$-part}
\newtheorem*{newtonser}{Newton diagrams for series}
\newtheorem*{newtonpol}{Newton diagrams for polynomials}
\newtheorem*{stev1}{Strong candidate poles for $\lgt(h)=1$}
\newtheorem*{stev2}{Strong candidate poles for $\lgt(h)=2$}
\newtheorem*{stevd}{Strong candidate poles for $\lgt(h)=d>2$}
\newtheorem*{anlt}{Complex analytic set up}

\title{Quasi-ordinary power series and their zeta functions}

\author[E. Artal]{E. Artal Bartolo}
\address{Departamento de Matem\'aticas\\
Universidad de Zaragoza\\
Campus Pza. San Francisco s/n\\
E-50009 Zaragoza SPAIN}
\email{artal@unizar.es}
\author{Pi. Cassou-Nogu\`es}
\address{ Institut de Math\'ematiques de Bordeaux\\
Universit\'e Bordeaux I\\
350, Cours de la Lib\'era\-tion
33405 Talence Cedex, France.}
\email{cassou@math.u-bordeaux.fr}
\author{I. Luengo}
\author[A. Melle]{A. Melle Hern\'andez}
\address{Departamento de \'Algebra\\
Universidad Complutense\\
Plaza de Ciencias 3\\
E-28040 Madrid SPAIN}
\email{iluengo@mat.ucm.es, amelle@mat.ucm.es}

\date{}

\thanks{
First author is partially supported by
BFM2001-1488-C02-02; the last two authors are partially
supported by BFM2001-1488-C02-01}

\keywords{Motivic, topological and Igusa zeta functions,
monodromy, quasi-ordinary singularities}

\subjclass[2000]{14B05,14E15,32S50}

\begin{document}

\begin{abstract} The main objective of this paper is to
prove the monodromy conjecture for the local Igusa zeta function 
of a quasi-ordinary polynomial of arbitrary dimension 
defined over a number field. In order to do it, 
we compute the local Denef-Loeser motivic zeta function
$Z_{\text{DL}}(h,T)$
of a quasi-ordinary power series $h$ of arbitrary dimension
over an algebraically closed field of characteristic zero from its
characteristic exponents without using embedded
resolution of singularities.
This allows us to effectively represent
$Z_{\text{DL}}(h,T)=P(T)/Q(T)$ such that
almost all the candidate poles given by $Q(T)$
are poles. Anyway, these candidate poles give eigenvalues of the
monodromy action of the complex of
nearby cycles on $h^{-1}(0).$
In particular we prove in this case
the monodromy conjecture made by Denef-Loeser
for the local motivic zeta function and the
local topological zeta function.
As a consequence, if $h$ is a quasi-ordinary polynomial
defined over a number field
we prove
the Igusa monodromy conjecture for its
local Igusa zeta function.
\end{abstract}

\maketitle

\section{Introduction}

Let $h$ be a polynomial in $\bz[x_1,\ldots,x_d]$ and fix $p\in\bz$
a prime number. In order to study the number $N_k$
of solutions
of the congruence $h\equiv 0 \text{ mod } p^k,$
classically one associates with
$h$ the Poincar\'e series
$$
P(T)=\sum_{k=0}^\infty N_k T^k.
$$
J.~Igusa proved in \cite{ig:74}
that $P(T)$ is a rational function on $T$
by relating it with the following $p$-adic integral
\begin{equation}
I(h,s) := \int_{\bz_p^d} \vert h (x) \vert^s \vert d x \vert,
\end{equation}
for $s \in \bc$, ${\rm Re} (s) > 0$,
where $\vert d x \vert$ denotes the Haar measure on
$\bq_p^d$ normalized in such of way that $\bz_p^d$ is of volume 1.

Igusa proved the identity
$P(p^{-d-s})=\frac{1-p^{-s}I(h,s)}{1-p^{-s}}$
and he used an embedded resolution of
$h^{-1}(0)$ to show that
$I(h,s)$ is a rational function on $p^{-s},$
(see \cite{de:84} for a proof without resolution of singularities).
In fact each exceptional divisor of an embedded resolution
gives a candidate pole of $I(h,s)$ but many of them are not
actually poles.

Igusa conjectured that poles
of $I(h,s)$ are  related with eigenvalues
of the complex monodromy at
some point of  $h^{-1}(0),$ see \cite{de:91}.
More precisely he raised the \emph{conjecture}:

\begin{imc}
If $h\in F [x_1, \ldots, x_d]\setminus F$,
for some number field $F\subset \bc,$
then for almost all $p$-adic completion $K$ of $F,$ if $s_0$
is a pole of $I(h,K,s)$, then $\exp(2i\pi\Re(s_0))$
is an eigenvalue of the local monodromy of $h$ at some complex point of
$h^{-1}(0).$
\end{imc}

Since then, some partial results have been obtained, see
the Bourbaki Seminar talk by J. Denef \cite{de:91} for
a survey of these results  until 1991, (\cite{lo:88,lo:90,rv:01,ac:03}).
We recall two cases
where Igusa monodromy conjecture has been proved.
In fact in both cases, the following   stronger version of
the Igusa monodromy conjecture was proved:
for almost all $p$-adic completion $K$ of $F,$ if $s_0$
is a pole of $I(h,K,s)$, then $\Re(s_0)$ is a
root of the Bernstein polynomial
$b_h(s)$ of $h.$
Firstly, F. Loeser in \cite{lo:88}
gave a proof for reduced polynomials
in two variables.
One other interesting case
is the case of  polynomials (in arbitrary dimension)
non-degenerated with respect to their Newton polyhedron 
and verifying some ``resonance''
conditions. This result was also proved by
F. Loeser in \cite{lo:90}.
For plane curves,
the strong candidate poles come from
the rupture components in an embedded resolution
of $h^{-1}(0)$, see \cite{vys:01}.
For non-degenerated polynomials J. Denef gave
a set of strong
candidate poles which comes from faces of codimension $1$
in the Newton polyhedron of $h$, see \cite{de:95}.

In this paper we will prove the Igusa monodromy conjecture
for the local Igusa zeta function $I_0(h,K,s),$
which is the local version
of $I(h,K,s)$, for quasi-ordinary polynomials
in arbitrary dimension.
As we will see quasi-ordinary polynomials
behave in many aspects as plane curves.
In our proof we use in a essential way
motivic integration
on the space of arcs on an algebraic variety.

Motivic measure take values
in a completion of the Grothendieck ring
of algebraic varieties.
Let $K_0(\text{Var}_k)$ be the Grothendieck ring of algebraic varieties
over $k$. Let $\bl=[\ba_k^1]$ denote the class in
$K_0(\text{Var}_k)$ of the affine line. The naive motivic ring $\mm_k$
of algebraic varieties over $k$ is the polynomial ring
$\mm_k= K_0(\text{Var}_k)[\bl^{-1}].$

Let $X$ be a non singular irreducible complex algebraic variety
of pure dimension $d.$
For any $n\in \bn$, let
$\cl_n(X)$ denote
the space of arcs modulo $t^{n+1}$  on $X$; it has a structure of
complex variety.
The arc space $\cl(X)$ of $X$ is the projective limit of
the algebraic varieties $\cl_n(X).$
For any $n\in \bn$, let
$\pi_n:\cl(X)\to \cl_n(X)$ be the natural projection. For any
arc $\barphi\in \cl(X),$ the origin of the arc is $\pi_0(\barphi).$
For each closed point $x\in X,$ let $\cl_x(X)$ (resp. $\cl_{n,x}(X)$) be
the set of  arcs on $X$ (resp. truncated arcs) with origin at $x.$

Consider $X=\ba^d$ the $d$-dimensional complex affine space
and $x=\mathbf 0$ its origin. Let
$h\in \bc[x_1,\ldots,x_d]$
be a complex polynomial, with $h(\mathbf 0)=0.$
Set
$
V_n:=\{\barphi\in \cl_x(X): \,\ord(h\circ \barphi)=n\}.
$
Denef and Loeser in \cite{dl:01} defined
the \emph{naive motivic zeta function} of $h$
by
$$
Z_{\text{naive}}(h,T):=\sum_{n\geq 1}[\pi_n(V_n)]\bl^{-nd}T^n\in{\mm}_k[[T]].
$$
In fact we will consider
the local Denef-Loeser motivic zeta function
which is nothing but $Z_{DL}(h,T)=\bl^{-d}Z_{\text{naive}}(h,T)$.
They showed in \cite{dl:98}, using embedded resolution of
singularities, that
$Z_{DL}(h,T)$ is a rational function.
It belongs to the subring
$\mathcal N$
of the ring $\widehat{\mm}_k[[T]]$ which is generated by the image
in $\widehat{\mm}_k[[T]]$ of $\mm_k[T]$ and $(1-\bl^{-a}T^b)^{-1},\,
a,b \in\bn,b>0 .$
The monodromy conjecture
in this case states, see \cite[section 2.4]{dl:98}:

\begin{mmc} There is a set $S= \{(a,b):\, a,b\in\bn,b>0 \}$
such that $Z_{DL}(h,T) \in \mm_k[T][(1-\bl^{-a}T^b)^{-1}]_{(a,b) \in S}$
and if $q=a/b,(a,b) \in S $, 
 then
$\exp(2 i\pi q)$ is an eigenvalue of the
local complex algebraic monodromy
around zero at some $P\in h^{-1}(0)$.
\end{mmc}

It turns out that $Z_{DL}(h,T)$
is the right function to study several monodromy conjectures
because it specializes to the local Igusa zeta function
and to the local topological zeta function:

$\bullet$ If $h$ is a non-constant polynomial defined
over a number field
$F$ then it follows
from \cite{dl:98}, see also \cite{dl:01},
that for almost all finite places
of $F,$ the real parts $N s+\nu$ of poles
of $I_0(h,K,s)$ come from factors $(1-\bl^{-\nu}T^N)$
in the denominator of $Z_{\text{naive}}(h,T).$

$\bullet$ The local topological zeta function $\zlo(h,s)$ of a complex
polynomial
$h$ was introduced by Denef and Loeser in
\cite{dl:92} as a kind of limit of the
local Igusa zeta function.
Later on, they obtained $\zlo(h,s)$
from $Z_{DL}(h,T)$ by the following procedure, cf.
\cite[\S~2.3]{dl:98}.
First
substitute $T$ by $\bl^{-s}$ in $Z_{DL}(h,T),$
then expand $\bl^{-s}$ and $(\bl-1)(1-\bl^{-\nu+N s})^{-1}$
into series in $\bl-1.$ Finally take the usual Euler characteristic
$\chi_{\text{top}}$
(this works because $\chi_{\text{top}}(\bl)=1).$
Then the local topological zeta function is the following rational function,
see \cite{dl:01,dl:98}:
\begin{equation*}
\zlo(h,s):=\chi_{\text{top}}(Z_{DL}(\bl^{-s})).
\end{equation*}
It is clear that poles of $\zlo(h,s)$ induce
poles of $Z_{DL}(h,T).$ The other way around is not true,
see Example \ref{nash}.
In this case  Denef and Loeser conjectured in \cite{dl:92}:

\begin{tmc}
If $s_0=\nu+s N$ is a pole of $\zlo(h,s)$
then $\exp(2i\pi(-\nu/N))$ is an eigenvalue
of the monodromy action at some point of $h^{-1}(0)$.
\end{tmc}

One result of this paper is  that if
$h$ defines a quasi-ordinary singularity of hypersurface
of arbitrary dimension
then the monodromy conjecture
for $Z_{DL}(h,T)$ is true, see Corollary \ref{maincoro}.
This fact implies that the Igusa monodromy conjecture
for $I_0(h,K,s)$ and the Denef-Loeser monodromy conjecture for
$\zlo(h,s)$. In fact we prove a stronger result.

A germ of a complex analytic variety  $(V,0)$ is called
a \emph{quasi-ordinary  singularity} if there is a finite morphism
(proper with finite fiber map) of analytic germs
$\pi:(V,0)\to (\bc^d,0)$  whose discriminant locus
is contained in $x_1x_2\dots x_d=0$,
for some local coordinates
$(x_1,\ldots, x_d).$
We are only interested in quasi-ordinary hypersurface singularities
$(V,0)\subset(\bc^{d+1},0).$

A convergent power series
$h\in\bc\{\bdx,z\},$ $h(0)=0$,
defines a \emph{quasi-ordinary
singularity} at zero if the germ of its
zero locus $(V,0)\subset(\bc^{d+1},0)$ is
a \emph{quasi-ordinary  singularity}.
These singularities behave in many aspects
as
singularities of plane curves mainly because,
after Jung-Abhyankar theorem, they admit fractional power
series parameterizations and a finite set of
\emph{characteristic exponents}, see \cite{ab:55},  \cite{lu:83}.
J.~Lipman and Y.~Gau proved that
in the irreducible case
these exponents determine the embedded topology
of $(V,0)$, \cite{ga:88,li:88}.
Recently, several embedded resolutions
constructed from the characteristic exponents
have been obtained, see the works of O. Villamayor
\cite{vi:02} or the Ph.D. of P.D.~Gonz\'alez P\'erez
\cite{go:01}.

In \cite{lu:83}, I.~Luengo showed that the usual Newton-Puiseux
method for curves can be used in arbitrary dimension
to find the roots with fractional exponents of a quasi-ordinary
polynomial by means of \emph{Newton maps} (see paragraph
\ref{newtonmaps}). The key point to find the roots is that after a Newton map
we get a new quasi-ordinary polynomial with less
characteristic exponents
and one proceeds by finite induction.
In this paper we do not use embedded resolution
of quasi-ordinary singularities to compute
$Z_{DL}(h,T)$ but
the above procedure given by  Newton maps.
Such a method allow us to decompose $Z_{DL}(h,T)$ as
a sum of two rational motivic zeta functions:
$$
Z_{DL}(h,T)=Z_{DL}^A(h,T)+Z_{DL}^B(h,T).
$$

The $A$-part corresponds to arcs $\barphi \in \cl_0(X)$
such that the $t$-order $\ord_t(h\circ \barphi)$ can be computed
from the (degenerated or not) Newton polytope of $h.$
The computation of $Z_{DL}^A(h,T)$ from the Newton polytope
follows ideas of J.~Denef and K.~Hoornaert in the $p$-adic
case, see \cite{dh:01}. In fact we will show
an interesting description of $Z_{DL}^A(h,T)$
in terms of generating functions of  some rational 
polyhedra obtained from the Newton polytope of $h,$
see sections \ref{secmotint} and \ref{secnondeg}.
As an application
a formula for $Z_{DL}(h,T)$ for any germ of
complex analytic function $h$ with non-degenerated Newton polytope
is given, see theorem~\ref{nondeg}.
G.~Guibert has recently obtained a similar formula, see~\cite{gu:03}.

In order to compute the measure of the arcs $\barphi \in \cl_0(X)$ in the $B$-part we pull back
these arcs under Newton maps. In particular $Z_{DL}^B(h,T)$
is the sum of some motivic zeta functions depending on the pull-back
of $h$ under all its Newton maps. There is here one major technical problem. 
In dimension higher that $1,$
usually there exist some arcs which cannot be lifted
under usual Newton maps. To solve this problem
we need to consider Newton maps with coefficients in $\bc\{t\}.$
Thus our results really deal with
quasi-ordinary power series with coefficients in  $\bc\{t\}$
instead of  $\bc.$ Now Newton maps are $\bc\{t\}$-morphisms
in the terminology introduced by Denef and Loeser in
\cite{dl:02} and we can apply the change of variables formula. 
The differential form plays a role here
but throughout this introduction we omit it. 
In fact, we perform all these computations 
for an algebraically closed field $k$ of characteristic zero.

Essentially in $Z_{DL}^B(h,T)$ we get quasi-ordinary
singularities with less characteristic exponents
and we can apply recursively this formula.
In this way, for $\zlo(h,s)$
a very effective and closed recursive formula is given
only depending on the tree of characteristic exponents,
cf. Theorem \ref{mainformula}.
For $Z_{DL}(h,T)$
the formula, also recursive, 
is enough to give a short list of candidate
poles.
Each characteristic exponent is a rational $d$-tuple.
Each non-zero coordinate of each characteristic exponent
will give a candidate pole of $Z_{DL}(h,T).$
In section $6$ some of them are excluded
to get
a smallest set $SCP(h)$ of strong candidate poles.
Our main result is
\begin{equation}\label{zcoeff}
Z_{DL}(h,T)\in\bz[\bl,\bl^{-1},(1-\bl^{-\nu}T^{N})^{-1}][T]_{(N,\nu)\in SCP(h)}.
\end{equation}

The proof of the above result gives some extra information.
Namely all elements in $SCP(h)$, but one concrete case $(*)$
when $d=2$ (cf. Proposition \ref{long2poles}), appear also as
strong candidate poles of a transversal section
of $h$ at some point of the singular
locus of $h^{-1}(0).$
Transversal sections at generic points
are also quasi-ordinary singularities,
this fact allows us to prove the motivic monodromy conjecture
by induction on the dimension. For curves, this method
gives a closed formula for $Z_{DL}(h,T)$, cf. (\ref{fcurvascase}), and from this formula  a new proof of the monodromy conjecture follows directly.
In the $(*)$-case, there is only one strong candidate pole not appearing
in the transversal sections. We use the formula,
proved by
P.D.~Gonz\'alez P\'erez, L.J.~McEwan and A.~N\'emethi in \cite{g:03},
for the zeta function of
the monodromy at the origin of quasi-ordinary singularities
to show that this strong candidate pole
also gives an eigenvalue of the monodromy
 of $h$ at the origin.

\bigskip

The computation of $\zlo(h,s)$  for quasi-ordinary
singularities gives a very effective way to compute the poles
of $\zlo(h,s)$ for a general surface singularity. To short this
work we will provide the details in a forthcoming paper.
The basic idea is to use the well-known \emph{Jung-method}
as follows.

Let $p:(S,0)\subset (\bc^3,0)\to (\bc^2,0)$ be a finite morphism of a (hyper)surface
singularity $(S,0),$ defined by $h,$
with discriminant locus $(\Delta,0)\subset (\bc^2,0).$
Take an embedded resolution
$\pi:(\widetilde{\bc^2},\mathcal D)\to(\bc^2,0)$ of the germ of
curve $(\Delta,0)$ where the exceptional locus
$\mathcal D=\pi^{-1}(0)$ is a normal crossing divisor.
The pull-back
$\tilde p: (\tilde S,\mathcal E)\to (\widetilde{\bc^2},\mathcal D)$ is a
finite map whose discriminant is contained in the pull-back of
$(\Delta,0).$ Thus the singularities
of $(\tilde S,\mathcal E)$ are all quasi-ordinary.
We have proved that the poles of $\zlo(h,s)$
are contained in the set of poles of either the transversal
sections (now curves) at the rupture components
of the resolution or components of the singular locus
of $(S,0).$ This gives a short list of strong candidate poles
which are in general poles of $\zlo(h,s).$ Next step
will be to prove that they induce eigenvalues of the monodromy.


\tableofcontents

{\bf Conventions.} Throughout this paper
we denote by $\bn$ the set of the nonnegative integers,
$\bp$ the set of positive integers,
$\br_+=\{x\in\br:\,x\geq0\,\}$
and $\br_{>0}=\{x\in\br:\,x>0\,\}.$
To shorten the notation we will use
bold symbols for $d$-tuples,
for instance $\bdx=(x_1,\ldots,x_d).$

In this paper we work over a field $k$ of characteristic zero.
A variety over $k$ will mean a reduced separated scheme of finite type
over the field $k$, $\bgd:=\spec k[x_1,\ldots,x_d,
x_1^{-1},\ldots,x_d^{-1}]$ denote the
$d$-dimensional torus over $k$ and
$\ba_k^d:=\spec k[x_1,\ldots,x_d]$ the
$d$-dimensional affine space over $k.$

\section{Motivic integration}\label{secmotint}

In this section we recall several results form
\cite{dl:98,dl:99,dl:01,dl:02}.
We refer to these papers for the proofs of such results.
In the first two section
we work over a field $k$ of characteristic zero.

\subsection{Grothendieck ring of varieties}
\mbox{}

The Grothendieck ring of algebraic varieties
over $k$, denoted by $K_0(\text{Var}_k)$, is the free Abelian group on isomorphism classes $[X]$
of algebraic varieties $X$ over $k$ subject to the relations $[X]=[X-Y]+[Y]$
where $Y\subset X$ is a closed subvariety of $X$.
The Cartesian product of varieties
gives the ring structure.
The following properties that we will freely use throughout the paper,
hold in $K_0(\text{Var}_k)$.

\begin{enumerate}[(1)]
\item If $f:Y\to Z$ is a fibre bundle with fibre $F$ which is locally trivial
in the Zariski topology, then $[Y]=[F][Z].$
\item If a variety $X$ is partitioned by locally closed subvarieties
$X_1,\ldots,X_n,$ then $[X]=[X_1]+\ldots+[X_n].$
\item If $f:Y\to Z$ is a bijective morphism, then $[Y]=[Z].$ The proof
of this property is deduced from the proof
of the same property at the level of virtual Hodge polynomials which
can be found for instance in \cite{ch:96}.
\end{enumerate}

Let $\bl=[\ba_k^1]$ denote the class in
$K_0(\text{Var}_k)$ of the affine line. The naive motivic ring $\mm_k$
of algebraic varieties over $k$ is the polynomial ring
$\mm_k= K_0(\text{Var}_k)[\bl^{-1}].$
Let $F^m \mm_k$ denote the subgroup of $\mm_k$ generated
by $[X]\bl^{-i}$ with $\dim X-i\leq -m$ and $\widehat{\mm}_k$
denote the completion of $\mm_k$ with respect to the filtration
$F^\cdot.$ This completion was first introduced by M.~Kontsevich.

\subsection{The arc space of a variety}
\mbox{}

Let $X$ be a nonsingular irreducible algebraic
variety over $k$ of pure dimension $d.$
For any $n\in \bn$, let
$\cl_n(X)$ denote
the space of arcs modulo $t^{n+1}$  on $X$ which has a structure of
$k$-variety, whose $K$-rational points,
for any field $K$ containing $k,$  are the $K[t]/(t^{n+1})K[t]$-rational
points of $X.$
The arc space $\cl(X)$ of $X$ is the projective limit of
the algebraic varieties $\cl_n(X).$
For any $n\in \bn$, let
$\pi_n:\cl(X)\to \cl_n(X)$ be the natural projection. For any
arc $\barphi\in \cl(X),$ the origin of the arc is $\pi_0(\barphi).$
For any closed point $x\in X,$ let $\cl_x(X)$ (resp. $\cl_{n,x}(X)$) be
the arcs (resp. truncated arcs) with origin at $x.$
The above definitions extend to the case where $X$ is
a reduced and separated scheme of finite type over $k[t].$ For any
$n\in\bn$, $\cl_n(X)$ is the $k$-scheme which represents the
functor
$$
R\mapsto \mor_{k[t]-schemes} (\spec R[t]/t^{n+1}R[t], X)$$
defined in the category of $k$-algebras and again $\cl (X)$
is its projective limit. The truncation map will be also denoted by
$\pi_n:\cl(X)\to \cl_n(X).$

Let $A$ be a \emph{semialgebraic}, resp.
$k[t]$-\emph{semialgebraic}, subset of $\cl(X)$;
it is called
\emph{stable at level} $n\in \bn$ if $A=\pi_n^{-1} \pi_n (A).$
We remark that if $A$ is stable at level $n$
then it is stable
at level $n'\geq n.$
The set $A$ is called \emph{stable} if it is stable at some level $n.$
A subset $A\in \cl(X)$ is \emph{cylindrical at level} $n$ if $A=\pi^{-1}_n(C)$
with $C$ a constructible set, and $A$ is \emph{cylindrical} if it is cylindrical
at some level. Denote by $\bdB^t$ the set of all
$k[t]$-semialgebraic subsets of $\cl(X).$

The \emph{motivic measure} on $\cl(X)$
is the unique map $\mu_X:\bdB^t\to \widehat{\mm}_k$
such that:
\begin{enumerate}[(a)]
\item If $A\in \bdB^t$ is stable at level $n,$ then
$\mu_X(A)=[\pi_n(A)]\bl^{-(n+1)d}.$
\item If $A\in \bdB^t$ is contained in $\cl(S)$ with $S$ a reduced closed
subscheme of $X\otimes_k k[t]$ with $\dim_{k[t]}S<\dim X,$ then $\mu_X(A)=0.$
\item Let $A_i\in \bdB^t$ for all $i\in \bn.$ Assume that the $A_i$'s are mutually
disjoint and that $A:=\cup_{i\in \bn} A_i$ is $k[t]$-semialgebraic. Then
$\sum_{i\in \bn} \mu_X(A_i)$ converges in $\widehat{\mm}_k$ to $\mu_X(A).$
\end{enumerate}

Because there might exist cylindrical subsets of $\cl(X)$
which are not semialgebraic the motivic measure has been extended
to a measure, also denoted $\mu_X,$ defined over the Boolean
algebra of the \emph{measurable} subsets of $\cl(X)$,
see \cite[Appendix]{dl:02}. The above properties hold for
measurable subsets of $\cl(X)$ too.

For a measurable subset $A$ in $\cl(X)$ and a function
$\alpha:A\to \bz\cup\{\infty\}$,
we say that $\bl^{-\alpha}$ is \emph{integrable} on $A$ if the fibres
of $\alpha$ are measurable,
$\alpha^{-1}(\infty)$ has measure zero and the
\emph{motivic integral}
$$
\int_{A}\bl^{-\alpha} d\mu_X:=\sum_{n\in \bz}
\mu_X(A\cap \alpha^{-1}(n))\bl^{-n}\in \widehat{\mm}_k
$$
converges in $\widehat{\mm}_k.$

\begin{defini}\label{t-morphism} Let $X$ and $Y$ be $k$-varieties. A function
$\pi:\cl(Y)\to \cl(X)$ will be call a
$k[t]$-\emph{morphism} if it is induced
by a morphism of $k[t]$-schemes $Y\otimes_k k[t]\to X\otimes_k k[t].$
\end{defini}

\begin{thm}[Change variables formula] \label{changevariables} Let $X$ and $Y$ be smooth
$k$-varieties of pure dimension $d.$ Let
$\pi:\cl(Y)\to \cl(X)$ be a $k[t]$-morphism.
Let $A$ and $B$ be $k[t]$-semialgebraic subsets of $\cl(X)$ and $\cl(Y),$ respectively.
Assume that $\pi$ induces a bijection between $B$ and $A$. Then, for any
function
$\alpha:A\to \bz\cup\{\infty\}$ such that
$\bl^{-\alpha}$ is integrable on $A$, we have
$$
\int_{A}\bl^{-\alpha} d\mu_X=\int_B\bl^{-\alpha\circ \pi-\ord_t J_\pi(y)}
d\mu_Y,
$$
where $\ord_t J_\pi(y)$, for any $y\in \cl(Y),$ denotes the $t$-order of
the Jacobian of $\pi$ at $y.$
\end{thm}

\subsection{Local Denef-Loeser motivic zeta function}
\mbox{}

Let
$h\in k[[x_1,\ldots,x_d]]$
be a formal power series in the maximal ideal of the formal
power series ring.
Let $X:=\ba_k^d$ be the $d$-dimensional affine space
and $x=\mathbf 0$ its origin.
Set
$
V_n:=\{\barphi\in \cl_x(X): \,\ord(h\circ \barphi)=n\}.
$
The \emph{local Denef-Loeser motivic zeta function of} $h$ is the power series
\begin{equation}
Z_{DL}(h,T):=\sum_{n\geq 1} \mu_X(V_n)T^n\in{\mm}_k[[T]].
\end{equation}
Since $V_n$ is a stable semialgebraic set at level $n$
of $\cl_{n,x}(X)$, then we have
$\mu_X(V_n)=[\pi_n(V_n)]\bl^{-(n+1)d}\in\mm_k$.
In fact $Z_{DL}(h,T)$ belongs to the subring
$\mathcal N$
of the ring $\widehat{\mm}_k[[T]]$ which is generated by the image
in $\widehat{\mm}_k[[T]]$ of $\mm_k[T]$ and $(1-\bl^{-a}T^b)^{-1},\,
a\in\bn,b\in\bp.$
Denef and Loeser in \cite{dl:01} introduced
the \emph{naive motivic zeta function of} $h$
as
$$
Z_{\text{naive}}(h,T):=\sum_{n\geq 1}[\pi_n(V_n)]\bl^{-nd}T^n,
$$
then $\bl^d Z_{DL}(h,T)=Z_{\text{naive}}(h,T)$.

We will work in a slightly more general
set up. In order to be able to use the change variables formula
we consider $Z_{DL}(h,\omega,T)$ where $\omega$ is a regular differential form on $X$ and
the pair $(h,\omega)$ will verify the following condition.

\begin{spc}\label{spc}  The pair $(h,w)$
satisfies the support condition if and only if
\begin{enumerate}[\rm(1)]
\item $h(\bdx)=\prod_{j=1}^d x_j^{N_j} f(\bdx),\,N_j\in\bn$, with
$x_j$ does not divide $f$
for any $j=1,\ldots,d$,
\item $\omega$ is a regular differential form of type
$\omega=\prod_{j=1}^d x_j^{\nu_j-1}d x_1\wedge\ldots\wedge
d x_d, \,\nu_j\geq 1$,
\item $N_j=0$ implies  $\nu_j=1$, for any $j=1,\ldots,d$.
\end{enumerate}
\end{spc}

Set $
V_{n,m}:=\{\barphi\in V_n\ |\ \ord(\omega\circ\barphi)=m\}.$
For a given $n$,
there are finitely many $m$ such that $ V_{n,m}\neq \emptyset$, because of the support
condition, cf. \cite{vy:01}. The \emph{local Denef-Loeser motivic zeta function
of a pair
$(h,\omega)$} is the rational
function 
\begin{equation}
\label{zdl}
Z_{DL}(h,w,T):=\sum_{n\in\bn} \left(\sum_{m\in\bn}
\bl^{-m}\mu_X(V_{n,m})\right)T^n\in \mathcal{N}.
\end{equation}

The local topological zeta function $\zlo(h,\omega,s)$ is obtained
from $Z_{DL}(h,w,T)$ by the following procedure, see
\cite[\S~2.3]{dl:98}.
First
substitute $T$ by $\bl^{-s}$ in $Z_{DL}(h,\omega,T),$
then expand $\bl^{-s}$ and $(\bl-1)(1-\bl^{-\nu+N s})^{-1}$
into series in $\bl-1.$ Finally take the usual Euler characteristic
$\chi_{\text{top}}$
in \'etale ${{\mathbf Q}}_{\ell}$-cohomology,
this works because $\chi_{\text{top}}(\bl)=1.$
Then $\zlo(h,\omega,s)$ is the  rational function
\begin{equation} \label{defintopol}
\zlo(h,\omega,s):=
\chi_{\text{top}}(Z_{DL}(h,\omega,\bl^{-s})).
\end{equation}
We will use the symbol
$\chi_{\text{top}}(\bullet (\bl^{-s}))$
to denote the composition of the above three
operations whenever it has sense.

\begin{obs}\label{kt} See \cite[Remark 1.19]{dl:02} to generalize the results presented here
to schemes over $k[[t]]$ instead over $k[t].$
\end{obs}

\section{Generating functions
and Newton polyhedron}\label{secnondeg}

\subsection{Generating functions for integer points in
rational polyhedra}
\label{gener-cones}
\mbox{}

In this section some well known facts about generating functions of
rational polyhedra are reviewed.
We use as a reference \cite[Section 4.6]{st:86}  and \cite{bp:99}.

Let $\br^d$ be the Euclidean $d$-space with the
standard scalar product $\bdx \cdot \bdy =\sum_{l=1}^d x_l y_l$.
A \emph{rational polyhedron} $P\subset \br^d$ is the set of
solutions of a finite system of linear inequalities with integer coefficients:
$$
P:=\{\bdx\in \br^d: \bbeta_i \cdot \bdx \leq c_i \text{ for }
i=1,\ldots,m\},
\text{ where } \bbeta_i\in\bz^d \text{ and } c_i\in\bz.$$
A bounded rational polyhedron is called a \emph{polytope}.
A nonempty polyhedron is called a \emph{cone} if $\lambda \bdx\in P$
whenever $x\in P$ and $\lambda\geq 0.$ A pointed polyhedral cone
is a cone which does not contain a line.

The \emph{algebra of polyhedra} ${\mathcal P}(\br^d)$
is the $\bq$-vector space
spanned by the indicator functions $[P]$
of all polyhedra
$P\subset \br^d,$ where
the \emph{indicator function} $[P]:\br^d\to \br$ of $P$ is
defined by
$$
[P](x)=
\begin{cases}
1&\text{ if }x\in P,\\
0&\text{ if }x\not\in P.
\end{cases}
$$

We will use the same notation $[\bullet]$ for indicator functions
and for elements in the Grothendieck ring of algebraic varieties,
nevertheless we hope no confusion will arise.

Let $P\subset\bq^d$ be a rational polyhedron, with the set
of integral points  in $P$ we associate
the generating function
$$
\Phi_P(\bdx):=
\sum _{\balpha \in P \cap \bz^d} x_1 ^{\alpha _1}\dots
x_d ^{\alpha _d}.
$$
These series define a map $\Phi:{\mathcal P}(\bq^d)\to \bq(\bdx)$
with the following properties:

\begin{enumerate}
\item if $P_1,\ldots,P_r\subset \br^d$  are rational polyhedra whose
indicator functions satisfy a linear identity
$
\alpha_1[P_1]+\ldots+\alpha_r[P_r]=0,$ with $\alpha_i\in\bq,$
then
$$
\alpha_1\Phi_{P_1}(\bdx)+\ldots+\alpha_r\Phi_{P_r}(\bdx)=0.$$
\item If $\bdg+P$ is a translation of $P$ by an integer vector
$\bdg\in\bz^d$ then $\Phi_{\bdg+P}(\bdx)=\bdx^\bdg\, \Phi_P(\bdx)$.

\item $\Phi_{\{\mathbf{0}\}}(\bdx)=1$.
\end{enumerate}

Let $\mathcal C$ be a pointed polyhedral cone.
The one dimensional faces of
$\mathcal C$ are called \emph{extreme rays}.
A pointed polyhedral cone has only
finitely many extreme rays and it is the convex
hull of its extreme rays.
A \emph{simplicial} cone $\sigma$ is an $e$-dimensional
pointed convex polyhedral cone with $e$ extreme rays;
it may be also defined as a cone
generated by $e$ linearly independent
integer vectors $\bbeta_1,\ldots,\bbeta_e$, thus
$\sigma=\{\lambda_1\bbeta_1+\ldots+\lambda_e \bbeta_e,
\lambda_i\in \br_+\}$.

A triangulation of $\mathcal C$ consists of a finite collection
$\Gamma =\{ \sigma _1, ..., \sigma _t \}$ of simplicial cones such that
\begin{enumerate}[i)]

\item $ \cup \sigma _i =\mathcal C$;

\item if $\sigma \in \Gamma$, then every face of $\sigma$ is in $\Gamma$;

\item $\sigma _i \cap \sigma _j$ is a common face of $\sigma _i $ and
$ \sigma _j$.
\end{enumerate}
It is proved that a pointed convex polyhedral cone $\mathcal C$
possesses a triangulation $\Gamma$ whose $1$-dim
elements are the extreme rays of $\mathcal C$.
We will always consider these triangulations in the following.

We are mainly interested in \emph{positive} points in a
pointed convex polyhedral cone $\mathcal C$. Define
$E:=\mathcal C\cap \bn^d$, resp.  $\overline {E}:=\mathcal C\cap \bp^d.$
Their generating functions are computed using triangulations.
The boundary of $\mathcal C $, denoted $\partial \mathcal C ,$
is the union of all the facets of $\mathcal C $.
If $\Gamma $ is a triangulation of $\mathcal C $,
let
$\partial \Gamma$ denote
the set $\{ \sigma \in \Gamma, \sigma \in \partial \mathcal C \}$, and
$\overline {\Gamma }= \Gamma \setminus \partial \Gamma.$
Let $\sigma \in \Gamma$ be a simplicial cone, we set
$E_\sigma:=\sigma\cap \bn^d$ and
$\overline {E}_{\sigma }:=\{ v \in E_\sigma: \,v
\notin E_\tau, \forall \tau \subset \sigma \}$. Then $\overline {E}$
is the disjoint union
$\cup_{\sigma\in \overline \Gamma} \overline {E}_{\sigma}$
and
$$ \Phi_{\overline {E}}(\bdx)=\sum _{\sigma \in \overline {\Gamma }}
\Phi_{\overline {E}_{\sigma}}(\bdx).$$

Therefore to compute generating functions of cones,
we  compute generating
functions of simplicial cones.
Let $\bda_1,...,\bda_t$ be a set of linearly independent integer vectors
which generate the cone
$${F}:=\{ \blambda \in \bn^d:\, n\blambda= \lambda _1\bda_1
+\ldots+\lambda_t\bda_t,\, n \in \bn,\,  \lambda_i \in \bn \},$$
$\bda _1,...,\bda_t$ will be called a
\emph{set of quasi generators} of $F$.
Define the interior ${\overline F}$ of $F$ as
$${\overline F}:=\{ \blambda \in \bp^d:\, n\blambda= \lambda _1\bda_1
+\ldots+\lambda_t\bda_t,\, n \in \bp,\,  \lambda_i \in \bp \}.$$
Consider the finite set
${\overline D}_F:=\{\blambda \in {F}:\, \blambda=\lambda_1\bda_1
+\ldots+\lambda_t\bda_t,\, 0< \lambda_i \leq 1 \}.$
For any $\blambda\in {\overline F},$ there exist unique
$\bbeta\in{\overline D}_F$ and $\lambda_1,\ldots,\lambda_t \in \bn$
such that
$\blambda=\bbeta+\lambda _1\bda_1
+\ldots+\lambda_t\bda_t.$
Then \cite[Prop. 4.6.8]{st:86} yields:
\begin{equation} \label{simplicial}
\Phi_{{\overline F}}(\bdx)=\frac{
\left(\sum _{\bbeta \in {\overline D}_F} \bdx ^{\bbeta}\right)}
{\prod _{i=1}^{t}\left(1-\bdx^{\bda_i}\right)}.
\end{equation}

Given $F$, there is a unique set $CF({\overline F}):=
\{\bbeta_1,...,\bbeta_t\}$
of primitive quasi generators.
We say that ${\overline F}$ is strictly generated by
$\bbeta_1,...,\bbeta_t$ and call
$$ {G}_{\overline F}:=\{\blambda \in \bp^d:\, \blambda=\lambda_1\bbeta_1
+\ldots+\lambda_t\bbeta_t,\, 0< \lambda_i \leq 1 \}$$
the fundamental set of ${\overline F}$.

If $\mathcal C $ is a pointed polyhedral cone and $\Gamma $
a triangulation of
$\mathcal C $, from \cite[Prop. 4.6.10]{st:86},
we know that $CF({\overline E}):=\cup _{\sigma \in \overline {\Gamma }}
CF({{\overline F}}_{\sigma})$
is the set of $\bbeta \in \mathcal C \cap \bp^d$,
which lie on extreme rays of $\mathcal C$
such that  $\bbeta \neq n\bbeta '$ for some $
n >1$ and $\bbeta ' \in \mathcal C \cap \bp^d$.
Furthermore in \cite[Theorem 4.6.11]{st:86} it is proved
that the rational function  $\Phi_{{\overline E}}(\bdx) \in \bq [\bdx][D(\bdx)^{-1}]$
where
\begin{equation} \label{denom}
D(\bdx)=\prod_{\bbeta\in CF({\overline E})}\left(1-\bdx^{\bbeta}\right).
\end{equation}

The last result we need is the following.
Let $\mathcal C $ be a pointed polyhedral cone,
let $b_1,...,b_d \in \bz$
such that for each $r \in \bn$,
the number $g(r)$ of points in
${\overline E}=\mathcal C \cap \bp^d $
such that $b_1\bda_1+...+\bda_d b_d=r$ is finite.
Let $G(\lambda)= \sum _{r \in \bn} g(r) \lambda ^r$.
Then
\begin{equation}\label{stanley}
G(\lambda)=\Phi_{\overline E}(\lambda^{b_1},...,\lambda ^{b_d}).
\end{equation}

\subsection{Motivic zeta function and Newton polyhedra}
\label{sec-mot-newton}
\mbox{}

\begin{newtonser}
Let $h=\sum_{\bdn\in\bn^d}a_\bdn \bdx^\bdn\in k[[\bdx]]$
be a formal power series with $h(0)=0.$
The \emph{support}  of $h$ is the set
$\supp (h)=\{\bdn\in\bn^d:\, a_\bdn\ne 0\}.$
The \emph{Newton polyhedron} $\Gamma(h)$ of $h$ is the convex hull
in $\br_{+}^d$ of the set $\bigcup_{\bdn\in \supp (h)} (\bdn+(\br_+)^d).$
The \emph{Newton polytope} or \emph{Newton diagram}
$ND(h)$ of $h$ is the union of all compact
faces of $\Gamma(h)$; the set of all compact faces
is denoted by $CF(h)$.
The \emph{principal part} of $h$ is the polynomial
$h|_{ND(h)}:=\sum_{\bdn\in ND(h)} a_\bdn x^\bdn$. For any 
$\tau\in CF(h)$ we denote by $h_\tau$ the polynomial
$\sum_{\bdn\in \tau} a_\bdn x^\bdn.$
The principal part of $h$ is called \emph{non-degenerated}
if for each
closed proper face
$\tau\in CF(h)$,
the subscheme of $\bgd$  defined by
$$\frac{\partial h_\tau}{\partial x_1}=\ldots=
\frac{\partial h_\tau}{\partial x_d}=0$$
is empty.
\end{newtonser}

\begin{newtonpol}
Let $h:\ba_k^d\to\ba_k^1$ be a regular morphism,
$h(\bdx)=\sum_{\bdn\in\bn^d}a_\bdn \bdx^\bdn$.
The \emph{support}
of $h$ is the set
$\supp (h)=\{\bdn\in\bn^d:\, a_\bdn\ne 0\}$.
The \emph{global Newton polytope} $\Gamma_{\text{gl}}(h)$
of $h$ is the convex hull
in $\br_+^d$ of the set $\supp (h).$
The polynomial $h$ is called \emph{$0$-non-degenerated} if for each
closed face $\tau\subset \Gamma_{\text{gl}}(h)$, including
$\tau=\Gamma_{\text{gl}}(h)$, the subscheme of $\bgd$ defined by
$h_\tau=0$ is smooth over $k.$

The \emph{Newton polyhedron $\Gamma_\infty(h)$ of $h$
at infinity} is the convex hull of $\supp (h)\cup\{0\}.$
The polynomial $h$
is \emph{nondegenerated with respect to} $\Gamma_\infty(h)$
if for every face $\tau$ of $\Gamma_{\infty}(h)$
(of any dimension), which does not contain the origin,
the subscheme of $\bgd$ defined by
$$\frac{\partial h_\tau}{\partial x_1}=\ldots=
\frac{\partial h_\tau}{\partial x_d}=0$$
is empty. If $k$ is algebraically closed then the subscheme of $\bgd$ defined by
$h_\tau=0$ is smooth over $k$ if and only if  $h_\tau$,
$x_1\frac{\partial h_\tau}{\partial x_1},\ldots,
x_d\frac{\partial h_\tau}{\partial x_d}$
have no common zero on
the torus $\bgd$.
\end{newtonpol}

Let $h\in k[[\bdx]]$
be a formal power series and let
$\omega$ be a regular differential form such that
$(h,\omega)$ verifies the support condition
\ref{spc}. Assume that
$h(\bdx)=\prod_{j=1}^d x_j^{N_j}f(\bdx)$, $N_j\in\bn$,
where $x_j$ does not divide $f$
for any $j=1,\ldots,d$,
and the form $\omega$ equals $\left(\prod_{j=1}^d x_j^{\nu_j-1}\right)d x_1\wedge\ldots\wedge
d x_d$, $\nu_j\geq 1$.

We recall more known definitions and properties.
For $\bdk=(k_1,\ldots,k_d)\in \br_+^d,$ we define
$m_h(\bdk):=\inf_{\bdx\in \Gamma(h)} \{\bdk\cdot\bdx\}$ and
$\sigma_\omega(\bdk):=\nu_1 k_1+\ldots+\nu_d k_d.$
Since  $h$ is obtained from  $f$ multiplying
by a monomial then $\Gamma(h)$ is a translation
of $\Gamma(f)$. In particular
$m_h(\bdk)=m_f(\bdk)+N_1 k_1+\ldots+N_d k_d$.

The first meet locus of $\bdk\in \br_+^d$ is
$
F(\bdk):=\{\bdx\in  \Gamma(h):\,\bdk\cdot \bdx=m_h(\bdk)\,\}.$
If $\tau$ is a face of $ \Gamma(h)$ (or $ \Gamma(f)$)  the cone
associated with $\tau$ is the convex polyhedral cone, in the dual space,
defined by
$\Delta_\tau:=\{ \bdk\in \br_+^d:\, F(\bdk)=\tau\,\}$.

It is well-known that the cones associated with
the elements of $CF(h)$ give a partition
of (the dual space) $\br_{>0}^d$ in a disjoint union
$\bigcup_{\tau\in CF(h)} \Delta_\tau$.
It turns out that for
each $\bdk=(k_1,\ldots,k_d)\in\bp^d$
there exists a unique
compact face $\tau$ such that $\bdk\in \Delta_\tau.$
Given $\tau\in CF(h)$
we define, see \cite{dh:01} for a similar definition
in the $p$-adic case,
\begin{equation} \label{stnt}
S_{\Delta_\tau}(h,\omega,T):=\sum_{\bdk\in \bp^d\cap\Delta_\tau}
\bl^{-\sigma_\omega(\bdk)}T^{m_h(\bdk)}.
\end{equation}
Using the recalled results on generating functions,
in particular (\ref{stanley}), one has
\begin{equation} \label{motgen}
S_{\Delta_\tau}(h,\omega,T)=
\Phi_{\bp^d\cap\Delta_\tau}(\bl^{-\nu_1}T^{p_1},
\ldots,\bl^{-\nu_d}T^{p_d}),
\end{equation}
where $(p_1,\ldots,p_d)\in\tau$,
for instance one of its vertices.
In what follows we write
$\Phi_{\Delta_\tau}(\bdx):=
\Phi_{\bp^d\cap\Delta_\tau}(\bdx)$.
The term $S_{\Delta_\tau}(h,\omega,T)$ can be computed as follows.
Take a partition of the cone $\Delta_\tau$ into
rational simplicial cones $\Delta_i$, $i=1,\dots,s$, then
$
S_{\Delta_\tau}(h,\omega,T)=\sum_{i=1}^s S_{\Delta_i}(h,\omega,T).
$
If $\Delta_i$ is the cone strictly generated by linearly
independent vectors $\bda_1,\ldots,\bda_r\in \bn^d$ then
(\ref{simplicial}) implies
\begin{equation}\label{sdeltai}
S_{\Delta_i}(h,\omega,T)=\left(\sum_{\bdg\in G_i} \bl^{-\sigma_\omega(\bdg)}T^{m_h(\bdg)}\right)
\prod_{j=1}^r\frac{1}{1-\bl^{-\sigma_\omega(\bda_j)}T^{m_h(\bda_j)}},
\end{equation}
where $G_i$ is the fundamental set of $\Delta_i\cap \bp^d$:
$$
G_i:=\bn^d\cap \left\{\sum_{j=1}^r \mu_j\bda_j\ \Big|\ \, 0<\mu_j\leq 1
\text{ for }
j=1,\ldots,r\right\}.
$$

The multiplicity mult$(\Delta_i)$ of $\Delta_i$ is the
cardinality of
$G_i$. It is also equal to the volume of
the parallelepiped spanned by $\bda_1,\ldots,\bda_r$ with respect to
the volume form $\tilde \omega$ on the vector space $V$ generated by
$\{\bda_1,\ldots,\bda_r\}$ normalized such that the parallelepiped spanned by
a lattice basis of $\bz^d\cap V$ has volume $1$.
We define
\begin{equation}\label{Jsimpl}
J_{\Delta_i}(h,\omega,s):=
\frac{\mult (\Delta_i)}{\prod (\sigma_\omega(\bda_j)+m_h(\bda_j) s)}.
\end{equation}
If $\Delta_\tau$ is a $r$-dimensional rational convex cone
and $\Delta_\tau=\cup \Delta_i$ is a decomposition in
rational simplicial cones $\Delta_i$ of dimension
$r$ such that
$\dim (\Delta_i\cap\Delta_j)<r$ for $i\ne j$,
then we define
$J_{\Delta_\tau}(h,\omega,s):=\sum J_{\Delta_i}(h,\omega,s)$.
In fact $J_{\Delta_\tau}(h,\omega,s)$
is nothing but
$J_{\Delta_\tau}(h,\omega,s)=
\chi_{\text{top}}((\bl-1)^r S_{\Delta}(h,\omega,\bl^{-s}))$
and it does not depend on the decomposition.

\medskip
The following lemma follows from (\ref{denom}),
(\ref{sdeltai}) and the main theorem in \cite{de:95}.

\begin{lema} \label{motana}
$S_{\Delta_\tau}(h,\omega,T)\in \bz[\bl,\bl^{-1},
(1-\bl^{-\sigma_\omega(\bda)}T^{m_h(\bda)})^{-1}][T],$
where $\bda$ belongs to the set of
vectors such that $\bda \cdot \bdx=M$
is a reduced integral equation of an affine hyperplane
containing $\tau\in CF(h)$.
\end{lema}

We denote by $N_\tau$
the subvariety of $\bgd$
defined by $\{h_{\tau}=0\}$.
The symbol $[N_\tau]$ means its class
in the ring $K_0(\text{Var}_k)$.
Define $
L_{\tau}^A(h):=\bl^{-d
}\left((\bl-1)^d-[N_{\tau}]\right)\in {\mathcal M}_k.$
Since 
$\Gamma(h)$ is a translation of $\Gamma(f),$ there exists
a natural bijection  between $CF(h)$ and
$CF(f)$. Thus we can also write
$[N_\tau]=[\bgd\cap \{f_{\tau}=0\}].$
Let $\tau\in CF(h)$ with $\dim(\tau)=d-r'$,
$1\leq r'\leq d.$

If $r'=d,$ it means $\tau$ is a vertex of $\Gamma(h)$,
then $[N_\tau]=0,$ $\Delta_\tau$ is a $d$-dimensional rational convex
polyhedron
and $L_\tau^A(h)=\bl^{-d}(\bl-1)^d.$
In such a case
\begin{equation} \label{nondeg-top0}
\chi_{\text{\rm top}}\left(
L_\tau^A(h) S_{\Delta_\tau}(h,\omega,\bl^{-s})\right)=
\chi_{\text{\rm top}}\left((\bl-1)^d S_{\Delta_\tau}
(h,\omega,\bl^{-s})\right)=
J_{\Delta_\tau}(h,\omega,s).
\end{equation}

If $r'\in\{1,\ldots,d-1\},$ then $h_\tau(\bdx)$
is a weighted homogeneous polynomial with more than one monomial.
The quasi-projective variety $W:=\ba_k^d\setminus \{h_\tau(\bdx)=0\}$
is $k$-isomorphic to the affine algebraic variety
$Y:=\{(\bdx,z)\in \ba_k^d\times \ba_k^1\,:\,zh_\tau(\bdx)=1\,\}.$
Under such an isomorphism, $W\cap \bgd$ is isomorphic to $Y\cap \bge$.
It implies that
$[N_\tau]=[\bgd]-[W\cap \bgd]=(\bl-1)^d-[Y\cap \bge ]$
in $K_0(\text{Var}_k)$.
Thus
$L_\tau^A(h)=\bl^{-d}[Y\cap \bge]$.

\begin{lema}[see e.g. \cite{dh:01}] \label{dh01}
 Let $g\in k[x_1,\ldots,x_d].$ If
$d-r':=\dim \Gamma_{\infty}(g)<d$
then there exists a coordinate change $T$ on the torus $\bgd$
such that $g(\bdx)=(g\circ T^{-1})(\bdy)=
{\tilde g}(y_1,\ldots,y_{d-r'})$, where
${\tilde g}\in k[y_1,\ldots,y_{d-r'},
(y_1\dots y_{d-r'})^{-1}].$

Moreover if $\tau$ is a face of $\Gamma_\infty(g)$ then
$T^{-1}(\tau)$ is a face of $\tilde g$ and
$g$ is non-degenerated with respect to $\Gamma_{\infty}(g)$
if and only if ${\tilde g}$ is non-degenerated with
respect to $\Gamma_{\infty}({\tilde g})$
\end{lema}

In our case  $g_\tau:=zh_\tau(\bdx)$
is a weighted homogeneous polynomial
whose Newton polyhedron
has dimension $d-r'+1<d+1$.
Applying Lemma \ref{dh01} we find  a
homogeneous polynomial ${\tilde g}_\tau(\bdy)$ in $d-r'+1$ variables
such that the variety $Y\cap \bge$
is $k$-isomorphic to
${\mathbf G}_{m,k}^{r'}\times G^*$ where $G^*_\tau=\{\bdy\in
{\mathbf G}_{m,k}^{d-r'+1}: {\tilde g}_\tau(\bdy)=1\}.$
It turns out that $\bl^{d}L_\tau^A(h)=[Y\cap \bge ]=(\bl-1)^{r'}[G^*_\tau]$
in $K_0(\text{Var}_k).$ The convex rational rational cone
$\Delta_\tau$ has dimension $r'$, thus
\begin{equation} \label{nondeg-topr}
\chi_{\text{\rm top}}\left(
L_\tau^A(h) S_{\Delta_\tau}(h,\omega,\bl^{-s})\right)=\chi_{\text{\rm top}}([G^*_\tau])
J_{\Delta_\tau}(h,\omega,s).
\end{equation}

Later on we will consider
germs degenerated with respect to their
Newton polyhedra. In the following definition we collect the
terms corresponding to the $A$-part in the decomposition
of the arc space according to these polyhedra.

\begin{defini}\label{zdla} Let $h\in k[[\bdx]]$ be a power series
and $\omega$ a differential form such that
$(h,\omega)$ satisfies
condition
(\ref{spc}).
The $A$-part or the part corresponding
to the Newton polyhedron of $h$ is defined by
the rational function
$$
Z_{DL}^A(h,\omega,T):=\sum_{\tau \in CF(h)}
L_{\tau}^A(h) \,S_{\Delta_\tau}(h,\omega)\in
{\mathcal M}_k[(1-\bl^{-\sigma_\omega(\bda)}T^{m_h(\bda)})^{-1}][T],$$
where $\bda$ belongs to the set of
vectors such that $\bda \cdot \bdx=M$
is a reduced integer equation of an affine hyperplane
containing some $\tau\in CF(h)$.
In the same way we define the $A$-part of $(h,\omega)$
for the local topological zeta function as
the rational function
$$
Z_{\text{\rm top},0}^A(h,\omega,s)=
\sum_{\tau\in CF(h)} \chi_{\text{\rm top}}([G^*_\tau])
J_{\Delta_\tau}(h,\omega,s),
$$
where $\chi_{\text{\rm top}}([G^*_\tau]):=1$ if $\dim(\tau)=0.$
\end{defini}

If $h$ has non-degenerated Newton principal part,
 $Z_{DL}(h,\omega,T)$ is written in terms
of some invariants
of $CF(h)$ and $\omega.$
This formula has been proved by
J. Denef and K. Hoornaert in the $p$-adic setting, \cite{dh:01}.
Our proof is based on their results,
so we follow their notation too. In fact
we provide the proof because it will give
some light in the arc decomposition that we have to do
in the quasi-ordinary case. G. Guibert
has communicated to the authors that he has recently obtained
a similar formula, \cite{gu:03}.

\begin{thm}
\label{nondeg} Let $h(\bdx)=\prod_{j=1}^d x_j^{N_j}f(\bdx),\,N_j\in\bn,$
be a regular function with $h$
non-degenerated with respect to
(all the compact faces of)
its Newton polyhedron and let
$\omega=\prod_{j=1}^d x_j^{\nu_j-1}d x_1\wedge\ldots\wedge
d x_d, \,\nu_j\geq 1$ be a differential form
verifying the support condition \ref{spc}.
Then
$$
Z_{DL}(h,\omega,T)=\sum_{\tau \in CF(h)}
L_{\tau}(h)S_{\Delta_\tau}(h,\omega,T)
$$
where $L_{\tau}(h):=L_{\tau}^A(h)+L_{\tau}^B(h)$ and
$
L_{\tau}^B(h):=\bl^{-d}(\bl-1)[N_{\tau}]
\frac{\bl^{-1}T}{1-\bl^{-1}T}.
$
\end{thm}

We break the proof in several steps.

\begin{paso0} Classifying arcs.
\end{paso0}
Let $S$ be the affine hypersurface $x_1\cdots x_d=0.$ Since
$S$ has dimension less than $d$ then by property $(b)$ of the motivic measure
the set $\cl_0(S)$ has measure zero. Then
we only consider arcs $\barphi\in\cl_0(\ba_k^d)\setminus\cl_0(S)$, i.e.,
$\barphi=(\varphi_1,\ldots,\varphi_d)$
where $\varphi_i(t)=a_{k_i}t^{k_i}+\text{higher degree terms}$.
As usual we write
$\bdk(\barphi):=(k_1,\ldots,k_d)\in\bp^d$ and
$\bda(\barphi):=(a_{k_1},\dots,a_{k_d})\in \bgd$.

Let $\tau$ be
the unique compact face of $ND(h)$ such that $\bdk=\bdk(\barphi)\in \Delta_\tau\cap\bp^d$.
If
$\bda=\bda(\barphi)\in
\bgd\setminus N_\tau$ then
$\ord (h\circ\barphi)=
\ord (h_\tau\circ\barphi)=m_h(\bdk)$.
Otherwise, if $\bda\in N_\tau$ then
$\ord(h\circ\barphi)>m_h(\bdk)$.
\smallbreak
The set
${\overline V}_{n,m}=\{\barphi\in \cl_0(\ba_k^d)\setminus\cl_0(S):
\,\ord(h\circ \barphi)=n,\,\ord(\omega\circ\barphi)=m\,\}$
can be decomposed as
${\overline V}_{n,m}=\bigcup_{\tau\in CF(f)}
\left(V_{n,A,m}^\tau\cup V_{n,B,m}^\tau\right)$,
where
$$
V_{n,A,m}^\tau:=\{\barphi\in {\overline V}_{n,m}: \,\bdk(\barphi)\in\Delta_\tau,\,
\bda(\barphi)\in
\bgd\setminus N_\tau\,\},
$$
$$
V_{n,B,m}^\tau:=\{\barphi\in {\overline V}_{n,m}: \,\bdk(\barphi)\in\Delta_\tau,\,
\bda(\barphi)\in
N_\tau\,\}.
$$
Rewrite the formula (\ref{zdl}) as:
\begin{equation}
\label{zdlab}
\begin{split}
Z_{DL}(h,w,T)=&\sum_{\tau \in CF(h)}\left(\sum_{n\geq 1}
\left(\sum_{m\in\bp} \bl^{-m}\mu_X(V_{n,A,m}^\tau)\right)T^n+\right.\\
+&\left.\sum_{n\geq 1}\left(\sum_{m\in\bp} \bl^{-m} \mu_X(V_{n,B,m}^\tau)\right)
T^n\right).
\end{split}
\end{equation}

\begin{paso0} Computation of the $A$-part of the series.
\end{paso0}

Let $\tau\in CF(h)$ and $\bdk\in\Delta_\tau$.
Define
$$
V_A^\bdk:=\{\barphi\in\cl_0(\ba_k^d)\setminus\cl_0(S): \bdk(\barphi)=\bdk,\,
f_\tau(\bda(\barphi))\neq 0\}.
$$
Given $\barphi\in V_A^\bdk$ then  $\ord(h\circ\barphi)$ and
$\ord(\omega\circ\barphi)$ depend only on $\bdk$, therefore
$V_A^\bdk\subset V_{n_\bdk,A,m_\bdk}^\tau$ where
$n_\bdk:=m_h(\bdk)=
m_f(\bdk)+N_1 k_1+\ldots N_d k_d$ and
$m_\bdk:=\sigma_\omega(\bdk)-(k_1+\dots+k_d).$
The stable semialgebraic set $V_A^\bdk$ has measure
$\mu_X(V_A^\bdk)=[\pi_{n_\bdk}(V_A^\bdk)]\bl^{-(n_\bdk+1)d}$. Since
$$
[\pi_{n_\bdk}(V_A^\bdk)]=[\bgd\setminus
\{f_\tau=0\}]\bl^{d n_\bdk-k_1-\ldots-k_d}
$$
then
$
\bl^{-m_\bdk}\mu_X(V_A^\bdk)=\bl^{-d}
\left((\bl-1)^d-[N_{\tau}]\right)\bl^{-\sigma_\omega(\bdk)}.
$

It turns out that
\begin{equation*}
\begin{split}
\sum_{n\geq 1}
\left(\sum_{m\in\bn} \bl^{-m} \mu_X(V_{n,A,m}^\tau)\right) T^n=
\sum_{\bdk\in \bp^d\cap\Delta_\tau}
\bl^{-m_\bdk}\mu_X(V_A^\bdk) T^{n_\bdk}=
\\
=\sum_{\bdk\in \bp^d\cap\Delta_\tau}\bl^{-d}
\left((\bl-1)^d-[N_{\tau}]\right)\bl^{-\sigma_\omega(\bdk)}
T^{m(\bdk)+N_1k_1+\ldots
N_d k_d}=\\
=\bl^{-d
}\left((\bl-1)^d-[N_{\tau}]\right)\sum_{\bdk\in \bp^d\cap\Delta_\tau}
\bl^{-\sigma_\omega(\bdk)}T^{m_h(\bdk)}=L_{\tau}^A(h) S_{\Delta_\tau}(h,\omega,T).
\end{split}
\end{equation*}

\begin{paso0}\label{parteB} Computation of the $B$-part of the series.
\end{paso0}
Let $\tau\in CF(h)$.
As before for any $\bdk\in \bp^d\cap
\Delta_\tau$
we define the sets $V_{n,B}^\bdk$
of arcs $\barphi$ in $V_{n,B}^{\tau}$
such that $\bdk(\barphi)=\bdk$.
Note that $V_{n,B}^{\bdk}\subset V^\tau_{n,B,m_\bdk}$.

\begin{lema}\label{lemab} Fix $\bdk\in \bp^d\cap\Delta_\tau$ as before.
Then
$$
\sum_{n\geq 1}
\bl^{-m_\bdk} \mu_X(V_{n,B}^\bdk) T^{n}=\bl^{-d}(\bl-1)[N_{\tau}]
\frac{\bl^{-1}T}{1-\bl^{-1}T}\bl^{-\sigma_\omega(\bdk)}T^{m_h(\bdk)}.$$
\end{lema}

The formula in Theorem \ref{nondeg} is deduced from the  following equality which is a consequence of the above lemma:
$$
\sum_{\bdk\in \bp^d\cap\Delta_\tau}
\sum_{n\geq 1}
\bl^{-m_\bdk}\mu_X(V_{n,B}^\bdk) T^{n}
=\bl^{-d}(\bl-1)[N_{\tau}]
\frac{\bl^{-1}T}{1-\bl^{-1}T}S_{\Delta_\tau}(h,\omega,T).
$$


\begin{proof}[Proof of Lemma \ref{lemab}]
Given $\bdk\in\bp^d\cap\Delta_\tau$, consider the $k[t]$-morphism
$\pi_\bdk:Y\to \ba_{k[t]}^d$
defined by
$
\pi_\bdk(y_1,\ldots,y_d)=(t^{k_1}y_1,\ldots,t^{k_d}y_d),
$
where $Y:=\ba_k^d.$ Let $\barphi\in V_{n,B}^\bdk$ be an arc, then the equalities
$\varphi_i(t)=t^{k_i}\psi_i(t)$, $i=1,\dots,d,$
define a unique arc
$\bpsi(t):=(\psi_1(t),\ldots,\psi_d(t))\in\cl(Y)$ centered at
$\bpsi(0)=\bda(\barphi)\in \bgd\cap\{f_\tau=0\}$.
It verifies $\ord(f\circ\pi_\bdk\circ \bpsi)=\ord(f\circ\barphi)=n$
and because of the quasi-homogeneity of $f_\tau$:
$$
f\circ\pi_\bdk(y_1,\ldots,y_d)=t^{m(\bdk)}(f_\tau(y_1,\ldots,y_d)+tg(y_1,\ldots,y_d)),
$$
where the function $g$ has coefficients in $k[t]$.
In particular $n=m_f(\bdk)+n'$ with $n':=\ord(f_1\circ \bpsi)\geq 1$
where  $f_1(y_1,\ldots,y_d):=f_\tau(y_1,\ldots,y_d)+tg(y_1,\ldots,y_d)$.
The algebraic set
$f_\tau(y_1,\ldots,y_d)=0$
is nonsingular at points of
$\bgd\cap\{f_\tau=0\}$
because of the non-degeneracy condition.
In particular at any of such points there are
coordinates
such that $f_\tau={\bar y}_1$ and in the same way
at such a point $f_1=0$ will be non-singular.
Applying the change variables formula
the proof of the lemma is finished.
Remark that the pull-back by $\pi_\bdk$ of the regular differential
form $d x_1\wedge\ldots\wedge d x_d$
is $t^{\sigma_\omega(\bdk)}dy_1\wedge\ldots\wedge dy_d.$
\end{proof}

\mbox{}

A formula for $Z_{\text{\rm top},0}(h,\omega,s)$ for
a non-degenerated function $h$
was given in \cite{dl:92}. We apply Theorem \ref{nondeg}
to reprove that formula.
We need several standard definitions.

\begin{defini}\label{volgamma}Let $\gamma$ be the convex hull in $\br^d$ of some subset of $\bz^d$. The \emph{volume form $\omega_\gamma$ on the affine space} $\text{Aff}(\gamma)$
generated by $\gamma$ is defined such that the parallelepiped generated by
a base of the lattice $\text{Aff}(\gamma)\cap\bz^d$ has volume $1$.
For any $\tau\in CF(h)$, $V(\tau)$
is the volume of $\tau$ relative to the induced lattice;
i.e. to the volume form $\omega_\tau$.
\end{defini}

\begin{thm}
\label{nondeg-top} \cite{dl:92}
If $h$ has
non-degenerated principal part and the pair $(h,\omega)$
verifies the support condition (\ref{spc})
then
\begin{multline*}
Z_{\text{\rm top},0}(h,\omega,s)=
 \sum_{\tau\in CF(h),
 \dim\tau=0} J_{\Delta_\tau}(h,\omega,s)\\
 +\left(\frac{s}{s+1}\right) \sum_{\tau\in CF(h),
 \dim\tau \geq 1}
 (-1)^{\dim(\tau)}\dim(\tau)! V(\tau) J_{\Delta_\tau}(h,\omega,s).
\end{multline*}
\end{thm}

The proof follows the same ideas as in equations (\ref{nondeg-top0})
and (\ref{nondeg-topr}). Let $\tau\in CF(h)$.
If $\tau$ is zero dimensional we have the identity (\ref{nondeg-top0}).
Otherwise, following the same notation as in equation (\ref{nondeg-topr}),
if $\dim(\tau)=d-r'$ then
$$
L_\tau(h)=\bl^{-d}\left((\bl-1)^{d+1}\frac{\bl^{-1}T}{1-\bl^{-1}T}+
(\bl-1)[Y\cap \bge]  \frac{1-T}{1-\bl^{-1}T}\right).
$$
Since $d-r'<d,$ the denominator of
$S_{\Delta_{\tau}}(h,\omega,T)$ has at most $d-1$
factors and
$$
\chi_{\text{\rm top}}
\left(\bl^{-d}(\bl-1)^{d+1}\frac{\bl^{-1}T}{1-\bl^{-1}T}
S_{\Delta_\tau}(h,\omega,\bl^{-s})\right)=0.
$$
Therefore we have again the term $[Y\cap \bge].$
Applying Lemma \ref{dh01} to
$g_\tau:=zh_\tau(\bdx)$ we have $[Y\cap \bge ]=(\bl-1)^{r'}[G_\tau^*]$
in $K_0(\text{Var}_k)$ where $G_\tau^*=\{\bdy\in
{\mathbf G}_{m,k}^{d-r'+1}: {\tilde g}_\tau(\bdy)=1\}$
for some non-degenerated weighted
homogeneous polynomial ${\tilde g}_\tau(\bdy)$.
Finally we also need the following result.

\begin{lema} \label{0nodeg} \cite{bkk:76,dl:91}
Let $g(y_1,\ldots,y_n)$ be
a $0$-non-degenerated polynomial and let $Z^*=\{y\in
{\mathbf G}_{m,k}^n: g(\bdy)=0\}.$
Then the $\ell$-adic Euler characteristic
$\chi_{\text{\rm top}}(Z^*)$ is given by
$$ \chi_{\text{\rm top}}(Z^*)=(-1)^{n-1}n!\vol_n(\Delta(g)).
$$
Here $\vol_n(\Delta(g))$ is the $n$-dimensional Euclidean volume.
\end{lema}

In our case the polynomial $g={\tilde g}_\tau-1$
is $0$-non-degenerated and verifies the hypothesis
of the Lemma \ref{0nodeg}, then
$ \chi_{\text{\rm top}}(G^*)=(-1)^{r}(r+1)!\vol_{r+1}(\Delta(g)).
$
Since
$\Delta(g)$ is a cone over the origin,
$(r+1)!\vol_{r+1}(\Delta(g))=r! V(\tau)$ (this last volume as in definition \ref{volgamma}). Hence we get the formula in Theorem \ref{nondeg-top}.

\begin{obs}
All the results in this section are also valid in the complex analytic set up.
\end{obs}

\begin{ejem}\label{nash} It is clear that the poles of $\zlo(h,\omega,s)$
come from the poles of $Z_{DL}(h,\omega,T)$
but
in general the other way around is not true. 
Consider for instance
$h=x_1^3+x_2^3+x^3_3+x^3_4+x_5^6$
and the differential form $\omega=d\bdx.$ This example has been recently 
studied by Ishii-Kollar in \cite{ik:02} to disprove a J.~Nash conjecture.

Then $h$ has non-degenerated
principal part and its Newton polyhedron has only one compact facet
$\tau$
which, in the dual, corresponds with the extreme ray $v_\tau:=(2,2,2,2,1).$
Take the vertex $(3,0,0,0,0)$ of the face ${\tau},$ then
$$
S_{\Delta_{\tau}}(h,\omega,T)=
\Phi_{\Delta_{\tau}}(\bl^{-1}T^3,\bl^{-1}T^0,\bl^{-1}T^0,
\bl^{-1}T^0,\bl^{-1}T^0)=\frac{\bl^{-9}T^6}{1-\bl^{-9}T^6}.$$
After Lemma \ref{motana}, $(1-\bl^{-9}T^6)$ is a possible
factor of the denominator
of $Z_{DL}(h,\omega,T)$.
In fact it is case.
To prove it, one considers the Igusa $p$-adic zeta function of $h$
and apply \cite[Theorem 5.17]{dh:01} which guarantees
that $t_0=-3/2$ is one of its  poles   for $p\gg0.$
Nevertheless the local topological zeta function is 
$\zlo(h,\omega,s)=1/(1+s)$.
This shows an example where poles of $p$-adic and
local topological zeta function are not the same.
\end{ejem}

\section{Quasi-ordinary power series}

In this section we recall some known properties
of the quasi-ordinary power series. We give the necessary details
to describe the method we will use to compute their local
zeta functions. Let $k$ be an algebraically
closed field of characteristic zero.

\mbox{}

Let $h\in k[[\bdx]][z]$
be a $z$-polynomial of degree $s$ with coefficients
in the formal power series $k[[\bdx]],$
$\bdx=(x_1,\ldots,x_d)$. Assume that
$h(\bdx,z)=x_1^{N_1}\dots x_d^{N_d}g(\bdx,z)$
where $N_l\geq 0$ and $x_l$ doesn't divide $g(\bdx,z),$ $\forall l, l=1,\ldots,d,$.
Let's write the decomposition of  $g({\mathbf 0},z),$
$g({\mathbf 0},z)=u \prod (z-z_b)^{s_b}$,
with
$s=\sum s_b$,
$u\in {\mathbf G}_{m,k}$ and $z_b\in \ba_k^1.$

The power series $h$
is called \emph{quasi-ordinary}
(or $k$-\emph{quasi-ordinary} if we want to emphasize the base field)
if its $z$-discriminant $D_z(g)$ (or equivalently $D_z(h)$) is
$$
D_z(g)=x_1^{\alpha_1}x_2^{\alpha_2}\dots x_d^{\alpha_d}\varepsilon(\bdx),
$$
where $\varepsilon({\mathbf 0})\not=0$. 
The condition on the $z$-discriminant implies
that $g(\bdx,z)$ is squarefree
in the ring $k[[\bdx,z]]$.

For each root $z_b$ of $g({\mathbf 0},z)$,
the Jung-Abhyankar Theorem, \cite{ab:55,lu:83},
states that
there are exactly $s_b$ distinct roots of $g(\bdx,z)$ centred at
$z_b$ in
$k[[x_1^{1/m_b},\ldots,x_d^{1/m_b}]],$ for some $m_b \in \mathbb P.$
Thus
for each $j\in\{1,\ldots,s_b\},$ there exists a fractional power series
$\zeta_{b_j}\in k[[x_1^{1/m_b},\ldots,x_d^{1/m_b}]]$ with
$\zeta_{b_j}({\mathbf 0})=z_b$
such that
\begin{equation}\label{descrip0}
g(\bdx,z)=\prod_{b}\prod_{j=1}^{s_b}
\left(z-\zeta_{b_j}(x_1^{1/m_b},\ldots,x_d^{1/m_b})\right).
\end{equation}
In fact, for each root $z_b,$ the product $\prod_{j=1}^{s_b}
\left(z-\zeta_{b_j}(x_1^{1/m_b},\ldots,x_d^{1/m_b})\right)$
is a well defined element in $k[[\bdx]][z].$
Moreover, this power series is quasi-ordinary  too
because of the properties of the discriminant.

In principle we are interested just
in one of the roots, say $z_b=0$, otherwise
we can a translation of type $\bar z=z-z_b$ to study
the corresponding root.
Later on when the transversal sections will be described
we will deal with other roots different from $z_b=0,$
see section \ref{transversal}.

Therefore we may assume
\begin{equation}\label{desc-w}
h(\bdx,z)=x_1^{N_1}\dots x_d^{N_d}f(\bdx,z)u(\bdx,z),
\end{equation}
where $f(\bdx,z)$
is a degree $n$ quasi-ordinary $z$-Weierstrass polynomial,
$n$ being the multiplicity
of the root $z_b=0$ in $g({\mathbf 0},z),$
$u(\bdx,z)$ is a unity in
$k[[\bdx,z]]$, $u(\bdx,z)$ being the product
of all other roots centred at points different from $z_b=0$, and
$N_l\geq 0$ for any $l=1,\ldots,d$.
Therefore,
for each $j\in\{1,\ldots,n\}$, there exist
$\zeta_j\in k[[x_1^{1/m},\ldots,x_d^{1/m}]]$ (for some positive integer
$m$) with
$\zeta_j({\mathbf 0})=0$
such that
\begin{equation}\label{descrip}
f(\bdx,z)=\prod_{j=1}^{n}
\left(z-\zeta_j(x_1^{1/m},\ldots,x_d^{1/m})\right).
\end{equation}

Let $\zeta$ be a root of $f,$  all its conjugates
$\tilde{\zeta}_{\boldsymbol{\varepsilon}}:=
\zeta(\varepsilon_1 x_1^{1/m},\ldots,\varepsilon_d x_d^{1/m}),$
where $\varepsilon^m_l=1,$ for all $l=1,\ldots,d,$ are also roots of $f.$
If $f$ is irreducible in  $k[[\bdx]][z],$
then all its roots are conjugate to
one of them, say $\zeta,$ and the set
$\{\tilde{\zeta}_{\boldsymbol{\varepsilon}}\}$
has exactly $n$ distinct elements $\{\zeta=\zeta_1,
\ldots,\zeta_n\}$ which verify
$$
f(\bdx,z)=\prod_{p=1}^n (z-\zeta_p).
$$

The irreducible factors of $f$ in
$k[[\bdx]][z]$ are also $k$-quasi-ordinary
power series. Let $\{f^{(i)}\}_{i\in I}$ be the set
of irreducible factors of $f.$
Then 
$$
f(\bdx,z)=\prod_{i\in I} \prod_{j=1}^{n^{(i)}} (z-\zeta^{(i)}_j)
$$
The set of indexes $I$ is decomposed in two disjoint subsets:
$I'$ and $I\setminus I'$ such that $\#(I\setminus I')\leq 1$, 
where in $\#(I\setminus I')=1$  if and only if  $z$ is an irreducible component
of $f$ (i.e. if $0$ is root of $f$).

\begin{defini}We will say that a variable $x_i$ is \emph{essential} for $h$ if
$x_i$ divides $D_z(h)$.
The \emph{number of essential variables}
will be denoted by $\lgt(h)$.
\end{defini}

\begin{obs} \label{nonred} All our results  for
the local Denef-Loeser motivic zeta function of quasi-ordinary
power series are proved in the case where $g$ is reduced. Ne\-ver\-the\-less
they
can be proved in the
nonreduced case.
This means that $g$ can have multiple components in
$k[[\bdx,z]]$.
In such a case the power series $h(\bdx,z)$ will be called
quasi-ordinary
if $x_1^{N_1}\dots x_d^{N_d}g_{\text{red}}(\bdx,z)$
is quasi-ordinary in the above sense.
\end{obs}

We are only interested in the local Denef-Loeser motivic zeta function.
Since for any arc $\barphi(t),$ with $\barphi(\mathbf 0)=0,$ the
$t$-order $\ord(h\circ \barphi)$ does not depend
on the unity $u(\bdx,z)$ then we usually suppose that
$u(\bdx,z)=1$.
\smallskip

\medskip

\subsection{Characteristic exponents}
\mbox{}

Given
$h(\bdx,z)=\prod_{l=1}^d x_l^{N_l}f(\bdx,z)=
\prod_{l=1}^d x_l^{N_l}
\prod_{i\in I}
(z-\zeta_{i})$, one has $D_z(f)=\prod (\zeta_k-\zeta_j)$. For two
different roots $\zeta_k$ and $\zeta_j$ of $f$, using
the
unique factorization of the discriminant, we have
$$
\zeta_k-\zeta_j=x_1^{\lambda_{kj,1}}
x_2^{\lambda_{kj,2}}\dots x_d^{\lambda_{kj,d}}
\varepsilon_{kj}(x_1^{1/m},\ldots,x_d^{1/m}),
$$
where $\varepsilon_{kj}(0)\not=0$, and $\lambda_{kj,t}\in\frac{1}{m}
\bdz_{\geq 0}.$

\begin{defini} We set $\Lambda_{CE}(f):=\{ \blambda_{kj}\}=
\{(\lambda_{kj,1},\ldots,
\lambda_{kj,d})\}_{kj}\subset  \frac{1}{m}
\bdz_{\geq 0}^d$. The elements of the finite set $\Lambda_{CE}(f)$ are called
\emph{characteristic exponents} of $f$. We usually identify
$\Lambda_{CE}(h)$ with $\Lambda_{CE}(f)$ because $\Lambda_{CE}(f)$
is the set of characteristic exponents
of $h$ at the root
zero. We usually omit ``at the root zero" whenever it will
be clear.

For each characteristic exponent
$\blambda=(\lambda_{1},\ldots,
\lambda_{d}),$ the corresponding monomial $\bdx^{\blambda}$
is called \emph{characteristic monomial} of $f$.
In the same way we identify characteristic monomials of $h$ at zero
with characteristic monomials of $f$.
\end{defini}

Since the irreducible factors $f^{(i)}$ of $f$ are also quasi-ordinary,
 the set of  characteristic exponents $\Lambda_{CE}(f^{(i)})$ verifies
$\Lambda_{CE}(f^{(i)})\subset\Lambda_{CE}(f) $
for every $i\in I.$ The other characteristic exponents of $f$ measure the order of coincidence
between distinct irreducible factors of $f$ (as in the case of plane curves).  

The following partial ordering will be used in the paper. For any
$\blambda,\bmu\in\bq^d,$ then $\blambda\leq \bmu$ if
$\lambda_l\leq \mu_l$ for all $l=1,\ldots,d.$ We say
that $\blambda<\bmu$ if
$\blambda\leq \bmu$ and there is an $l$ such that $\lambda_l< \mu_l.$
If $f(\bdx,z)$ is irreducible in $k[[\bdx,z]],$ then the elements
of $\Lambda_{CE}(f)$
are totally ordered, see \cite{li:88} and the references therein.
But this is no longer true in the non-irreducible case.

Up to know to have the same characteristic exponents is
just that two sets coincide. Later on, we will give
a new definition of having the same characteristic exponents
(which will extend the usual definition for curves).

\subsection{Newton polyhedron and good coordinates}
\mbox{}\label{Newton-good}

In this section we show that for every quasi-ordinary power series $h(\bdx,z)$
there exist  coordinates, called \emph{good coordinates},
such that the compact faces of $ND(h)$ have only dimension $1$ and $0$.
This means that the $1$-dimensional compact faces of $ND(h)$ are totally ordered
by their slopes and they form a \emph{monotone polygonal path}. Our definition 
is slightly more general than the definition of good coordinates given by P.~Gonz\'alez in \cite{go:01}, \cite[Lemma 3.16]{g:03};
we will call them $P$-good coordinates. He proves
the existence of \emph{$P$-good coordinates} for quasi-ordinary power series; in fact his $P$-good coordinates are a particular case
of our good  coordinates, see Remark \ref{pedro}. One of the main properties of   $P$-good coordinates is that they are uniquely defined, which is not the case for good coordinates. The advantage of these last ones is that they are preserved when passing to transversal sections, see subsection \ref{transversal}; this is not the case for $P$-good coordinates.

\begin{defini}
A quasi-ordinary power series $h\in k[[\bdx]][z]$ is in \emph{good coordinates}
if $ND(h)$ is a monotone polygonal path such that 
either $z$ is a component of $h$ or $h_{\gamma}$
is not the product of a monomial in $\bdx$ and a power of a linear form, where $\gamma$ is the compact $1$-dimensional face of $ND(h)$
which meets the plane $z=0,$ i.e. $h_{\gamma}\ne x_1^{a_1}\ldots x_d^{a_d}
(z-\alpha\, x_1^{b_1}\ldots x_d^{b_d})^m$ for all $\alpha\in {\mathbf G}_{m,k}$. 
\end{defini}

\begin{prop}\label{nuestras-buenas}
There exists a system of coordinates which are good for $h$.
\end{prop}

\begin{proof}
Since $h$ and $f$
essentially differ by  a monomial, the Newton polyhedron $\Gamma(h)$ is obtained from $\Gamma(f)$ by the 
translation induced by the corresponding monomial. Thus it is enough
to prove the result for $f$ being a Weierstrass quasi-ordinary polynomial of degree
$n.$ We proceed by induction on $n.$ The case $n=1$ is clear: we have 
$f=z+a(\bdx),$ with $a(\bdx)\in k[[\bdx]]$, 
and the change $z_1=z+a(\bdx)$ is enough.

In order to proceed by  induction we need the following notion introduced by Hironaka
and studied in \cite{lu:83}.

\begin{defini} Let $f(\bdx,z)=z^n+a_{n-1}(\bdx)z^{n-1}+\ldots+a_0(\bdx)$ be a quasi-ordinary
Weierstrass polynomial, set $\tau_0=(0,\ldots,0,n)$ the corresponding
vertex in $ND(f).$ We say that $f$ is \emph{$\nu$-quasi-ordinary with respect to $z$}
if there is a vertex $\tau_1$ in $ND(f),$ $\tau_1\ne \tau_0,$
such that if $\tau$ is the projection of $\tau_1$ over $\bn^d\times \{0\}$
from $\tau_0$ and $\gamma$ is the segment joining $\tau_0$ and $\tau$
then  $ND(f)$ is contained in $\cup_{{\mathbf n}\in \gamma} ({\mathbf n}+\br_+^{d+1})$ and $f_\gamma$ is a polynomial which is not a power of a linear form.
\end{defini}


In \cite{lu:83} it is proved that if we make the change 
$z=z_1-\frac{a_{n-1}(\bdx)}{n}$ the polynomial $f$ becomes $\nu$-quasi-ordinary.
If $ND(f)=\gamma $ (with the above notations) then $\gamma$ is the only $1$-dimensional
compact face, therefore the last one, and $f_\gamma$ is not a power of a linear form
so we are done. Otherwise, in \cite{lu:83} it is shown that $f$ can be decomposed as
a product $f=f_0(\bdx,z)\,f_1(\bdx,z), \, f_i(\bdx,z)\in k[[\bdx]][z],$
 such that $f_\gamma$ corresponds with the factor 
$f_0,$ recall that the Newton polytope of $f$ is the Minkowski sum of the Newton polytopes
of each factor. The factor $f_1$ is a quasi-ordinary Weierstrass polynomial of degree
less than $n$ and we apply induction to conclude. In fact the change of variables
we need for the induction step are of type $z_1=z_2+b(\bdx)$ and they do not
affect the first segment of Newton polytope of $f$ because of the condition $ND(f)$ is contained in $\cup_{{\mathbf n}\in \gamma} ({\mathbf n}+\br_+^{d+1})$ in the
definition of $\nu$-quasi-ordinary.   
\end{proof}

\begin{obs}
\label{bueno} Observe that from the proof of the proposition
\ref{nuestras-buenas} one gets
that it is enough to make change of variables of type
$z=z_1+a(\bdx),$ $a(\bdx)\in k[[\bdx]]$ to get good coordinates.
The procedure to get good coordinates is far from being unique. 
For instance $f(x,z)=(z^2-x^5)(z-x^2)(z^3-x^2)$ is already in good coordinates
but the described method gives us different good coordinates. 
\end{obs}

\medbreak

Let $h$ be a quasi-ordinary power series in good coordinates. Since
$ND(h)$ is a monotone polygonal path we order the set $\{\gamma\}$ of 
compact $1$-dimensional faces according with the slopes (which are rational $d$-tuples)
of the edges
using the partial order defined
before. In fact they are totally ordered
and the smallest one will be the $z$-highest.
Assume that there are exactly $r$ of them  which are
ordered as $\gamma_1,\ldots,\gamma_r.$
For
$q\in\{1,\ldots,r\},$ since $\gamma_q$ is 
$1$-dim, the weighted-homogeneous polynomial $f_{\gamma_q}=z^k\displaystyle
\prod x_l^{a_l}\,(z^{n_1^q m_q}+\ldots+\alpha (x_1^{b_1^q}\dots x_d^{b_d^q})^{m_q}),$
for some non negative integers $k,a_l$ and
$m_q:=\gcd(n_1^qm_q,b_1^qm_q,\ldots,b_d^qm_q).$ 
Since $\gamma_q$ is an edge of the Newton polytope
then $n_1^q$ is a positive integer, the $d$-uple $(b_1^q,\ldots, b_d^q)\in\bn^d$
is non-zero  and $\alpha\in{\mathbf G}_{m,k}$. The polynomial
$f({\mathbf 1},z^{1/n_1^q})$ of degree $m_q$ in $k[z]$
can be factorized.
There are some
positive integers ${m_{q,j}},j=1,\ldots,v(q),$
such that $f_{\gamma_q}=z^k\displaystyle
\prod x_l^{a_l}f_{\blambda_q}$ where
$$
f_{\blambda_q}:=\prod_{j=1}^{v(q)}
(f_{q,j})^{m_{q,j}},
\text{ where }
f_{q,j}:=
z^{n_1^{q}}-\beta_j^{q}
\bdx^{n_1^{q}\blambda_q}=z^{n_1^q}-\beta_j^{q}x_1^{b_1^q}\dots x_d^{b_d^q},
$$
$\lambda_{q,l}=\frac{b_l^q}{n_1^q}$
for each $l=1,\ldots,d$. We define the rational $d$-uple $\blambda_q:=(\lambda_{q,1},\ldots,\lambda_{q,d})\in \bq^d.$
We will say that $f_{\blambda_q}$ 
is the weighted-homogeneous polynomial associated with the
compact $1$-dimensional face $\gamma_q$ of $\Gamma(f)$.

\begin{obs} \label{pedro} Let us see that we can find good coordinates which are not $P$-good coordinates
in the example $f(x,z)=(z^2-x^5)(z-x^2)(z^3-x^2)$. The problem arises
in the characteristic exponents which are integers. The definition-construction 
of $P$-good coordinates shows that there exists
a change of coordinates of type $z_1=z+a(\bdx)$, $a(\bdx)\in k[[\bdx]],$ such that in the set $\Lambda_{ND}(h):=\{\blambda_1,\ldots, \blambda_{r}\}$,
the greatest characteristic exponent $\blambda_r$ is the unique  which can be in $\bz^d$; moreover, its Newton initial form $f_{\blambda_r}$ cannot be a power of a linear form.
In the above example
doing the change $z_1=z+2x^2$ we get a new Newton polytope
which is a polygonal path  with only two $1$-dimensional faces. 
We observe two facts. The first one is that for each quasi-ordinary power series in good coordinates,
the method described in  \cite[Lemma 3.16]{g:03} provides $P$-good coordinates. The second one is that
in any case, the characteristic exponents appearing in both Newton polyhedra which are smaller than the first integer characteristic exponent coincide in both polytopes. 
Our constructions work in $P$-good coordinates
but  proofs must be slightly modify.
\end{obs}

\begin{prop} If $h$ is in good coordinates then 
each $\blambda_q$ is a characteristic exponent
of $f.$  
\end{prop}

\begin{proof} Since $h$ is in good coordinates we can
use the process given in \cite{lu:83}
to compute the roots of $f$ from the
Newton polytope. In fact each root, different from zero,
is found in only one of the edges of $ND(h).$ Assume $\gamma_1,\ldots,\gamma_r$
are the edges ordered as before. Then any pair $\zeta\in\gamma_q$
and $\zeta'\in \gamma_{q+1}$ gives 
$\blambda_q$ as characteristic exponent,
that is
$\zeta-\zeta'=\bdx^{\blambda_q} \epsilon(\bdx),$ $\epsilon(\mathbf 0)\ne 0.$   
If $\gamma_q=\gamma_r$, as $h$ is in good coordinates,
either $z$ is a factor of $h$ and $ 
\zeta-0=\bdx^{\blambda_r}$, $\zeta \in \gamma_r$
or $f_{\gamma_r}$ is not a power of a linear form
and taking two distinct roots in $\gamma_r$ gives
${\blambda_r}$.
\end{proof}

Thus each root 
 $\zeta_i$ of $f$ is written in a unique
way as 
\begin{equation}\label{raiz}
\zeta_{i}(\bdx)=\alpha_1^{(i)}\bdx^{\blambda_{(i)}}
+\sum_{\blambda_{(i)}<\blambda}\alpha_{\blambda}^{(i)}
\bdx^{\blambda},\,\,(\alpha_1^{(i)}\ne 0).
\end{equation}
We define
$f_q$ to be  the product of all roots $\zeta$ of $f$ whose initial term
gives
the Newton polytope of $f_{\blambda_q}.$ This means that 
$$
f_q(\bdx,z):=\prod_{\text{roots}:\blambda_{(i)}=\blambda_q}
(z-\zeta_i).
$$
\begin{obs}\label{fq} Since all the conjugates under the Galois group
of a root $\zeta_i$ have the same
$\blambda_{(i)},$  the irreducible component in
$k[[\bdx]][z]$ of $f(\bdx,z)$ determined by  $\zeta_i$ divides
$f_q(\bdx,z).$ Thus $f_q\in k[[\bdx]][z]$
because it is a product of some
irreducible components of $f(\bdx,z).$
\end{obs}

Thus we define the set 
$\Lambda_{ND}(h):=\{\blambda_1,\ldots, \blambda_{r}\}$ which is a subset
of the characteristic exponents of $h$ and which
are ordered 
$\blambda_{r}>\ldots> \blambda_{1}.$

Let $f_j^q$ be the element in $k[[\bdx]][z]$ defined by
$$
f_j^q(\bdx,z):=
\prod_{\text{roots}:\blambda_{(i)}=\blambda_q\text{ and }
\beta^{(i)}=\beta_j^{q}}
(z-\zeta_i).
$$
The same ideas as in Remark \ref{fq} show
that $f_j^q\in k[[\bdx]][z]$.
Hence
$$
f_q=\prod_{j=1}^{v(q)}f_{j}^q
\quad\text{ and }\quad
f=\begin{cases}\displaystyle
\prod_{q=1}^{r}\prod_{j=1}^{v(q)}f_{j}^q,&
\text{if $I=I'$ and}\\
z\displaystyle\prod_{q=1}^{r}\prod_{j=1}^{v(q)}f_{j}^q,&
\text{otherwise.}
\end{cases}
$$

\begin{defini} \label{NC}
If $h$ is in good coordinates then we call $\{f_{j}^q\}_{q,j}$
the \emph{Newton components of} $f$; each $f_{j}^q$
defines a quasi-ordinary
power series. 
\end{defini}

\subsection{Dual decomposition}\label{dualsection}
\mbox{}

Assume a pair $(h,\omega)$ verifying the support condition ~\ref{spc}
is given,
where $h$ is a quasi-ordinary power series
in good coordinates as in (\ref{desc-w})
and
$\omega=\prod_{j=1}^d x_j^{\nu_j-1}
d x_1\wedge\ldots\wedge d x_d\wedge dz, \,\nu_j\geq 1.$
The compact faces of
$\Gamma({f})$ (or $\Gamma({h})$)
are $1$-dim edges
$\gamma_1,\dots,\gamma_r$
with their corresponding vertices
$\tau_0,\tau_1,\dots,\tau_r$.

\begin{obs}\label{order}
The main point in the Newton polytope is that its
set of compact faces (and then its vertices) is totally ordered.
In principle they are ordered by the order imposed to the characteristic exponents.
From now on, if we do not emphasize the contrary, we assume they are ordered
with the reverse order that is
such that
$\tau_r$ has the highest $z$-coordinate.
\end{obs}

Let $e_1,\dots,e_{d+1}$
denote the canonical basis of the dual space $(\br^{d+1})^*$
where we choose coordinates $(v_1,\ldots,v_d,v_{d+1})$. The
fan $\Sigma(\Gamma(h))\subset (\br^{d+1})^*$ is obtained subdividing
the cone $(\br_+^{d+1})^*$
with linear hyperplanes $l^q:\eta_q=0$ where
\begin{equation}\label{eq-hyper}
\eta_q:=\sum_{l=1}^d b_l^q v_l -n_1^q v_{d+1}=
n_1^q\left(\sum_{l=1}^d \lambda_{\kappa_{q,l}} v_l-v_{d+1}\right),\,\,q=1,\ldots,r.
\end{equation}
Let $\Delta_{\gamma_q}=\{\bdk\in \br_+^{d+1}\
|\  \bdk\cdot (b_1^q,\ldots,b_d^q,-n_1^q) = 0\}$ be its dual cone.

\begin{lema} \label{1connue}
$\Delta_{\gamma_q}$ is the strictly positive simplicial cone
$
\Delta_{\gamma_q}=
\{\lambda_1 \bdw_1^q+\ldots+\lambda_d \bdw_d^q,\,\lambda_i\in \br_+\,\},
$
where
$
\bdw_l^q:=\frac{1}{c_l^q}\left(n_1^q e_l+b_l^q e_{d+1}\right)\in\bn^{d+1},
\,c_l^q:=\gcd(n_1^q,b_l^q),\,{\bar b}_l^q:=\frac{b_l^q}{c_l^q},\,
p_l^q:=\frac{n_1^q}{c_l^q}
$ with $l=1\ldots,d.$
\end{lema}

The linearly independent vectors
$
\{\bdw_l^q\}_{
l=1,..d}$ are primitive, so the proof is clear.
Let $G_q$ be the fundamental set of 
$\Delta_{\gamma_q},$ then (\ref{simplicial}) implies
\begin{equation}
\Phi_{\Delta_{\gamma_q}}(\bdy):=\Phi_{\Delta_{\gamma_q}\cap \bp^{d+1}}(\bdy)=
\frac{\sum _{\bbeta \in G_q} \bdy ^{\bbeta}}
{\prod _{l=1}^{d}\left(1-\bdy^{\bdw_l^q}\right)}.
\end{equation}
If $p=(p_1,\ldots,p_d,p_{d+1})$ is an element
in the closure of ${\gamma}_q$ then
(\ref{motgen}) gives
\begin{equation}\label{motgen1}
S_{\Delta_{\gamma_q}}(h,\omega,T)=
\Phi_{\Delta_{\gamma_q}}
(\bl^{-\nu_1}T^{p_1},
\ldots,\bl^{-\nu_{d+1}}T^{p_{d+1}}).
\end{equation}
The edge $\gamma_q$ is defined
by the affine equations  $
n_1^q x_l+ b_l^q z=M_l^q,\quad l=1,\ldots,d,
$ for some positive integers $M_l^q$.

\begin{lema} \label{1connuevo}
$S_{\Delta_{\gamma_q}}(h,\omega,T)
\in\bz[\bl,\bl^{-1},
(1-\bl^{-(\nu_lp^q_l+\bar{b}_l^q)}T^{\frac{M_l^q}{c_l^q}})^{-1}]
[T],\, l\in \{1,\ldots,d\}$
and
$$
J_{\gamma_q}(h,\omega,s)=\frac{(n_1^q)^{d-1}}{\tilde{T_1}^q\dots \tilde{T_d}^q},
$$
where $\tilde{T_l}^q:=(M_l^q s+\nu_ln_1^q+b_l^q),\,l=1,\ldots,d$
and $q=1,\ldots,r.$
\end{lema}

\begin{proof} 
Since $m_h(\bdw_l^q)=\bdw_l^q \cdot \bda,$
where $\bda\in \gamma_q,$
then $m_h(\bdw_l^q)=\frac{M_l^q}{c_l^q}.$
On the other hand,
$\sigma_w(\bdw_l^q)=\bdw_l^q\cdot(\nu_1,\dots,\nu_d,1)
=\frac{1}{c_l^q}(\nu_l n_1^q+ b_l^q).$
Let us compute $J_{\gamma_q}(h,\omega,s).$
Because of $\gcd(b_1^q,\dots,b_d^q,-n_1^q)=1$  we choose
$\beta:=(a_1,\dots,a_l,u)$
such that $(b_1^q,\dots,b_d^q,-n_1^q)\cdot \beta=1$.
Computing $\text{mult}(\Delta_{\gamma_q})=
|\det(\bdw_1^q,\dots,\bdw_d^q,\beta)|$ by induction one can prove that
$$
\det(\bdw_1^q,\dots,\bdw_d^q,\beta)=
-(n_1)^{d-1}\frac{(b_1^q,\dots,b_d^q,-n_1^q)\cdot \beta}{\prod_{l=1}^d
c_l^q}=-\frac{(n_1^q)^{d-1} }{\prod_{l=1}^d c_l^q}.
$$
The result follows from (\ref{Jsimpl}).
\end{proof}

\medskip

Consider next the first vertex $\tau_0$ with coordinates $(\bda,d_0)$
in $\Gamma(h).$
The rational simplicial cone
$\Delta_{\tau_0}$ is the $(d+1)$-dimensional
cone:
$\br_+ \bdw_1^1+\dots+\br_+ \bdw_d^1+\br_+ e_{d+1}.$
In a similar manner as above if
$G_{\tau_0}$ is the fundamental set
of $\Delta_{\tau_0}$ then
\begin{equation*}
\Phi_{\Delta_{\tau_0}}(\bdy)=\frac{
\left(\sum _{\bbeta \in G_{\tau_0}} \bdy ^{\bbeta}\right)}
{(1-y_{d+1})\prod _{l=1}^{d}\left(1-\bdy^{\bdw_l^1}\right)},
\end{equation*}
and
\begin{equation}\label{motgen0}
S_{\Delta_{\tau_0}}(h,\omega,T)=\Phi_{\Delta_{\tau_0}}
(\bl^{-\nu_1}T^{a_1},
\ldots,\bl^{-\nu_d}T^{a_d},\bl^{-1} T^{d_0}).
\end{equation}
Recall that $I=I'$ if and only if 
$d_0=0$; otherwise $d_0=1$.

\begin{lema} \label{01con} $$S_{\Delta_{\tau_0}}(h,\omega,T)\in\bz[\bl,\bl^{-1},(1-\bl^{-1}T)^{-\varepsilon},
(1-\bl^{-(\nu_lp^1_1+\bar{b}_l^1)}
T^{\frac{M_l^r}{c_l^1}})^{-1}][T]
$$ 
with $l\in\{1,\ldots,d\}$
and
$$
J_{\tau_0}(h,\omega,s)=\frac{(n_1^1)^{d}}{(s+1)^\varepsilon\tilde{T_1^1}\dots \tilde{T_d^1}},\quad
\text{see Lemma \ref{1connuevo}},$$
where $\varepsilon=0$ or $1$ accordingly with $I=I'$ or not.
\end{lema}

\begin{proof}
If $I=I',$
$m_h(e_{d+1})=0$ and $\sigma_\omega(e_{d+1})=1$.
Furthermore, $\text{mult}(\Delta_{\tau_0})=
|\det(\bdw_1^1,\dots,\bdw_d^1,e_{d+1})|=\dfrac{(n_1^1)^d}
{\prod_{l=1}^d c_l^1}$, and we are done.
The other case is similar.
\end{proof}

Finally consider any other vertex $\tau_q$ of $\Gamma(h)$,
different
from $\tau_0$, which is the intersection of the edges
$\gamma_q$ and $\gamma_{q+1},$ $q=1,\ldots,r.$
If $q=r$ then $\gamma_{r+1}$ is a  (non-compact) 
 line parallel to $z$-axis.

The dual cone $\Delta_{\tau_q}$ is
$\sum_{i=1}^d (\br_+ \bdw_i^q +\br_+ \bdw^{q+1}_i),$
where if $q=r$ then $\bdw_l^{r+1}:=e_l,$ for $l=1,\ldots,d.$
For any $q=1,\ldots,r,$ consider the strictly positive simplicial cone
$$
\Delta_{q}^c:=\sum_{l=1}^d \br_+ \bdw_l^q+ \br_+ e_{d+1}.
$$
For instance $\Delta_{r+1}^c$ is 
the whole positive cone $\bp^{d+1}$.
Moreover its multiplicity verifies
\begin{equation}
\label{multi}
\text{mult}(\Delta_{q}^c)=
|\det(\bdw_1^q,\dots,\bdw_d^q,e_{d+1})|=\dfrac{(n_1^q)^d}{\prod_{l=1}^d
c_l^q}.
\end{equation}
The following relation holds for the indicator functions
of these cones:
\begin{equation}
\label{desmulti}
[\Delta_{\tau_q}]-[\Delta_{q+1}^c]+[\Delta_{\gamma_q}]+[\Delta_{q}^c]=0,
\qquad q=1\ldots,r,
\end{equation}
which implies the same relation among the corresponding
generating functions. Since $\Delta_{q}^c$ and $\Delta_{\gamma_q}$
are simplicial cones, 
\begin{equation} \label{genvert}
\Phi_{\Delta_{\tau_q}}(\bdy)=\Phi_{\Delta_{q+1}^c}(\bdy)-
\left(\Phi_{\Delta_{\gamma_q}}(\bdy)+\Phi_{\Delta_{q}^c}(
\bdy)\right).
\end{equation}
If $G_q^c$ is the fundamental set of $\Delta_{q}^c$ then
\begin{equation*}
\begin{split}
\Phi_{\Delta_{\tau_q}}(\bdy)&=
\frac{\left(\sum _{\bbeta \in G_{q+1}^c} \bdy ^{\bbeta}\right)}
{(1-y_{d+1})\prod _{l=1}^{d}\left(1-\bdy^{\bdw_l^{q+1}}\right)}-\\
&-\left(\frac{\left(\sum _{\bbeta \in G_{q}} \bdy ^{\bbeta}\right)}
{\prod_{l=1}^{d}\left(1-\bdy^{\bdw_l^q}\right)}+
\frac{\left(\sum _{\bbeta \in G_{q}^c} \bdy ^{\bbeta}\right)}
{(1-y_{d+1})\prod _{l=1}^{d}\left(1-\bdy^{\bdw_l^{q}}\right)}\right).
\end{split}
\end{equation*}
If $\tau_q\in\Gamma(h)$ has $z$-height $d_q$ 
 then
\begin{equation}\label{genver0}
S_{\Delta_{\tau_q}}(h,\omega,T)=\Phi_{\Delta_{\tau_q}}
(\bl^{-\nu_1}T^{a_1},
\ldots,\bl^{-\nu_d}T^{a_d},\bl^{-1} T^{d_q}).
\end{equation}
Hence since
$m_h(e_j)=N_j$ and $\sigma_\omega(e_j)=\nu_j$ then
\begin{equation}
\label{01con2}
S_{\Delta_{\tau_r}}(h,\omega,T)\in\bz\left[\bl^{\pm 1},
\frac{1}{1-\bl^{-\nu_j}T^{N_j}},
\frac{1}{1-\bl^{-(\nu_lp^r_l+\bar{b}_l^r)}
T^{\frac{M_l^r}{c_l^r}}}\right][T],
\end{equation}
with $j,l\in\{1,\ldots,d\}.$ In the same way,
for each $q\in\{1,\ldots,r-1\},$
\begin{equation}
\label{01con3}
S_{\Delta_{\tau_q}}(h,\omega,T)\in\bz\left[\bl^{\pm 1}, \frac{1}{1-\bl^{-(\nu_lp^q_l+\bar{b}_l^q)}T^{\frac{M_l^q}{c_l^q}}},
\frac{1}{1-\bl^{-(\nu_lp^{q+1}_l+\bar{b}_l^{q+1})}
T^{\frac{M_l^{q+1}}{c_l^{q+1}}}}\right][T],
\end{equation}
with $l=1,\ldots,d.$

To compute more explicitly
the local topological zeta function we argue as follows.
Recall that $J_{\tau_q}(h,\omega,s)$ can be computed using an adequate
simplicial decomposition of $\Delta_{\tau_q}.$ We use the decomposition
from (\ref{genvert}). Then
$$J_{\tau_q}(h,\omega,s)=\chi_{\text{top}}\left((\bl-1)^{d+1}
\left(\Phi_{\Delta_{q+1}^c}(\bda)-
\left(\Phi_{\Delta_{\gamma_q}}(\bda)+
\Phi_{\Delta_{q}^c}(\bda)\right)\right)\right),
$$
where $\bda:=(\bl^{-(\nu_1+a_1s)},
\ldots,\bl^{-(\nu_d+a_ds)},\bl^{-(1+d_qs)})$. Since
$\Phi_{\Delta_{\gamma_q}}(\bda)$ has at most $d$ poles then
$$
J_{\tau_q}(h,\omega,s)=\chi_{\text{top}}
\left((\bl-1)^{d+1}
\left(\Phi_{\Delta_{q+1}^c}(\bda)-
\Phi_{\Delta_{q}^c}(
\bda)\right)\right).
$$
Moreover we know from
equation (\ref{multi})
that $
\text{mult}(\Delta_{q}^c)=\dfrac{(n_1^q)^d}{\prod_{l=1}^d
c_l^q}.$

\begin{cor} \label{0con} For each $q\in \{1,\ldots,r-1\},$
then
$$J_{\tau_q}(h,\omega,s)
=\frac{1}{(1+d_q s)}\left(\frac{(n_1^{q+1})^{d}}{\tilde{T_1}^{q+1}
\dots\tilde{T_d}^{q+1}}-\frac{(n_1^q)^{d}}
{\tilde{T_1}^{q}\dots\tilde{T_d}^{q}}\right),$$
and if $q=r$ then
$$
J_{\tau_r}(h,\omega,s)=\frac{1}{(1+d_rs)}
\left(
\frac{1}{T_1\dots T_d}-
\frac{(n_1^r)^d}{\tilde{T_1}^{r}\dots\tilde{T_d}^{r}}
\right),$$
where $T_i:=N_i s+\nu_i$, $\tilde{T_l^r}$ have been defined in
Lemma \ref{1connuevo}, and $i,l\in\{1,\ldots,d\}.$
\end{cor}

\begin{obs}
It follows from the above discussion, definition \ref{zdla}
and Remark \ref{classqo} 
that
\begin{multline}
\label{eqA}
Z_{DL}^A(h,\omega,T) =
\sum_{j=0}^r L_{\tau_j}^A(h) S_{\Delta_{\tau_j}}(h,\omega,T)
+\sum_{q=1}^r L_{\gamma_q}^A(h) S_{\Delta_{\gamma_q}}(h,\omega,T)= \\
=\!\bl^{-(d+1)}(\bl-1)^{d+1}\left( \sum_{j=0}^r
S_{\Delta_{\tau_j}}(h,\omega,T)\!+\!
\sum_{q=1}^r
\left(1-\frac{v(q)}{\bl-1}\right)
S_{\Delta_{\gamma_q}}(h,\omega,T)\right),
\end{multline}
\begin{equation}
\label{eqAtop}
Z_{\text{top,0}}^A(h,\omega,s)=\sum_{j=0}^r J_{\tau_j}(h,\omega,s)
-\sum_{q=1}^r
{v(q)}
J_{\gamma_q}(h,\omega,s).
\end{equation}
After Lemma \ref{motana},
each $(1-\bl^{-1} T^{d_q})$ or
$(1+d_qs)$ appearing in
(\ref{genver0}) and Corollary \ref{0con} do not contribute
to the denominator of the $A$-part of the corresponding zeta functions of $(h,\omega)$.
\end{obs}


\subsection{Newton map associated with a Newton component}
\label{newtonmaps}\mbox{}

\bigbreak

We have already discussed the dual decomposition
associated with the Newton polyhedra of quasi-ordinary power series in good coordinates. In general
these power series are degenerated with respect to their Newton polyhedron.
In this paragraph we describe how to \emph{improve} a  quasi-ordinary
power series in terms of the complexity of the series.

Each compact $1$-dimensional face $\gamma_q$ of $\Gamma(f)$
corresponds to a polynomial
$$
f_{\blambda_q}=\prod_{j=1}^{v(q)}
(z^{n_1^{q}}-\beta_j^{q}
\bdx^{n_1^{q}\blambda_q})^{m_{q,j}}=
\prod_{j=1}^{v(q)}
(z^{n_1^q}-\beta_j^{q}x_1^{b_1^q}\dots x_d^{b_d^q})^{m_{q,j}}.
$$

Fix $j\in\{1,\dots,v(q)\}$, the corresponding polynomial
$f_{q,j}=z^{n_1^q}-\beta_j^{q}\bdx^{n_1^q\blambda_q}$ and the
series $f_j^q$, see definition
\ref{NC}.
We denote by $r\leq d$, the number of essential variables, $\lgt(f_{q,j})$, of the irreducible factor $f_{q,j}$. For the sake of simplicity we assume
that $x_1,\dots,x_r$ are the essential variables of $f_{q,j}.$
Fix $\pi_j^q:\ba_k^r\to V_j^q\subset \ba_k^{r+1}$
a bijective morphism (a parametrization)
of the irreducible variety 
defined by $f_{q,j}$:
$$
(s_1,\ldots,s_r)\mapsto (s_1^{p_1^q},\ldots,
s_r^{p_r^q},\alpha_j^q s_1^{{\bar b}_1^q}\dots s_r^{{\bar b}_d^r}),
$$
where $(\alpha_j^q)^{n_1^q}=\beta_j^{q}$.

\begin{obs} \label{classqo}
In $K_0(\text{Var}_k),$ one has
$[V_j^q]=\bl^r.$ Later on we will
consider the subvariety $(V_j^q\times\ba_k^{d-r})\cap\bge$ on the torus; its class in $K_0(\text{Var}_k)$ is
$(\bl-1)^d$.
Since $f_{\blambda_q}$ has $v(q)$ irreducible components,
$[N_{\gamma_q}]=v(q)(\bl-1)^d.$
\end{obs}

\begin{defini} \label{newtonmap}
Let $X=\spec k[[\bdx,z]], Y=\spec k[[\bdy,z_1]]$. The \emph{Newton map}
associated with the Newton component $f_j^q$
(or with the polynomial $f_{q,j}$) of $f$
is the morphism
$\pi_{q,j}:Y\to X$,
$\pi_{q,j}(\bdy,z_1):=(\bdx,z)$ defined
by the equations $x_l=y_l^{p_l^q}$ for every $l=1,\ldots,d,$ and $z=(z_1+\alpha_j^q)\prod_{l=1}^d y_l^{{\bar b}_l^q}$. 
\end{defini}

The  integers ${\bar b}_l^q$, $p_l^q$
were defined in Lemma \ref{1connue} and
if $l=r+1,\dots,d$, then ${\bar b}_l^q:=0$ and
$p_l^q:=1$.

\begin{lema}  \label{simpler}
Let $h$ be a quasi-ordinary power series in good coordinates.
The pull-back $\bar{h}_{q,j}$
of $h$ under the Newton map $\pi_{q,j}$
associated with a Newton component $f_j^q$ of $f$
defines a quasi-ordinary power series at the root $0.$
Moreover under any
Newton map $\lgt(h)$ does not increase. More precisely the
set of essential variables of $\bar{h}_{q,j}$ is contained in
the one of $h$ (identifying $x_l$ and $y_l$).
\end{lema}

\begin{proof} From the well-known properties of the discriminant
the following
identities will give the proof of the lemma.
$$
D_{z_1}(h\circ\pi_{q,j})=y_1^{m_1}y_2^{m_2}\dots y_d^{m_d}(D_z(h)\circ \pi_{q,j})=
y_1^{r_1}y_2^{r_2}\dots y_d^{r_d}\varepsilon(\bdy),
$$
where $\varepsilon(0)\ne 0, m_i\geq 0, r_i\geq 0$. The conditions on $\lgt(h)$ and the essential variables are also clear.
\end{proof}

Let $f_k^i$ be a Newton component of $f$,
hence it defines a quasi-ordinary power series.
The previous lemma shows that
the pull-back
$f^{i}_k\circ\pi_{q,j}$ can be decomposed as $y_1^{a_1}y_2^{a_2}\dots y_d^{a_d}\bar{f^{i}_k}$,
for some power series $\bar{f^{i}_k}\in k[[\bdy,z_1]].$

\begin{lema} \mbox{}

\begin{enumerate}[\rm(1)]
\item If $\blambda_{(i)}=\blambda_q$
and $\beta_k^i=\beta_j^q$
then
$\bar{f^{i}_k}({\mathbf 0})=0,$ therefore $\bar{f^{i}_k}$ defines
a quasi-ordinary power series.
\item Otherwise,
$\bar{f^{i}_k}({\mathbf 0})\ne 0,$ and $\bar{f^{i}_k}$
is a unit in $k[[\bdy,z_1]].$
\end{enumerate}
\end{lema}

\begin{proof}
The proof of the lemma will follow from the following description.
Take one of the roots $\zeta$ of $f$ different from $z=0$.
Since we are in a good system of coordinates, $\zeta$ is written as, see
identity (\ref{raiz}),
$$
\zeta=\alpha_1^{(i)}\bdx^{\blambda_{(i)}}
+\sum_{\blambda_{(i)}<\blambda}\alpha_{\blambda}^{(i)}
\bdx^{\blambda}.
$$
The pull-back of $\zeta$
under the Newton map $\pi_{q,j}$
is as follows.
\begin{enumerate}
\item If $\blambda_{(i)}=\blambda_q$
and
$(\alpha_1^{(i)})^{n_1^{(i)}}=\beta_j^q$ then
\begin{equation*}
(z-\zeta)\circ \pi_{q,j}(\bdy,z_1)
=\bdy^{\blambda_q}
(z_1-
\sum_{\blambda_q<\blambda}\alpha_{\blambda}^{(i)}
\bdy^{\blambda-\blambda_q} ).
\end{equation*}
\item If $\blambda_{(i)}=\blambda_q$
but
$(\alpha_1^{(i)})^{n_1^{(i)}}\not=(\alpha_j^{q})^{n_1^q}=
\beta_j^q$ then
\begin{equation*}
\left(z-\zeta\right)\circ \pi_{q,j}(\bdy,z_1)
=\bdy^{\blambda_q}
(z_1+\alpha_j^{q}-\alpha_1^{(i)}-
\sum_{\blambda_q<\blambda}\alpha_{\blambda}^{(i)}
\bdy^{\blambda-\blambda_q} ).
\end{equation*}
\item If $\blambda_{(i)}<\blambda_q$
then
\begin{equation*}
\left(z-\zeta\right)\circ \pi_{q,j}(\bdy,z_1)
=\bdy^{\blambda_{(i)}}
((z_1+\alpha_j^{q})
\bdy^{\blambda_q-\blambda_{(i)}}
-\alpha_1^{(i)}-
\sum_{\blambda_{(i)}<\blambda}\alpha_{\blambda}^{(i)}
\bdy^{\blambda-\blambda_{(i)}}).
\end{equation*}
\item If $\blambda_q<\blambda_{(i)}$
then
\begin{equation*}
\left(z-\zeta\right)\circ \pi_{q,j}(\bdy,z_1)
=\bdy^{\blambda_q }
(z_1+\alpha_j^{q}
-\alpha_1^{(i)}\bdy^{\blambda_{(i)}-\blambda_q}-
\sum_{\blambda_{(i)}<\blambda}\alpha_{\blambda}^{(i)}
\bdy^{\blambda-\blambda_q}  ).
\end{equation*}
\end{enumerate}
\end{proof}

The pull-back $\bar{h}_{q,j}$ of $h$ is $\prod_{l=1}^d
y_l^{N_l} \bar{f}_{q,j}(\bdy,z_1)w(\bdy,z_1)$
where
$w({\mathbf 0})\not=0.$
Moreover
the characteristic exponents of the
quasi-ordinary power series $\bar{f}_{q,j}(\bdy,z_1)$
are easily deduced from the characteristic exponents of $f$.
In general the new quasi-ordinary power series
$\bar{h}_{q,j}$ is not given
in a good system of coordinates. Nevertheless
the change of coordinates we need to put it in good coordinates is described
in Remark \ref{bueno}.

\begin{obs} \label{raices}
If  $\beta_j^q$
is a simple root of the polynomial
$f_{\blambda_{\kappa_q}}({\mathbf 1},z)$ then
$\bar{f_j^q}({\mathbf 0},z_1)=z_1+\cdots.$
Using the Implicit
Function Theorem the root of $\bar{f_j^q}$
is a series in $k[[\bdy]]$.
As in the Puiseux Theorem ($d=1$),
all  roots of $f$ can be found
by using a sequence of Newton maps,
(see \cite[Theorem 2]{lu:83}  for details).
\end{obs}

If $f$ (or $h$)
has non-degenerated Newton principal part then
all  roots of  all $f_{\blambda_q}({\mathbf 1},z)$
are simple roots.
Otherwise $f$ has degenerated principal part and
following the proof of  \cite[Theorem 2]{lu:83}
one can prove that
by successive Newton maps and maps of
type $z_1=z-m(\bdx),$
to get again a good system of coordinates,
we reach a point
in which all  pull-back of $f$
have non-degenerated principal part.

\begin{defini} \label{depth}
We define the \emph{depth} of a quasi-ordinary power series $h$,
and denote it by $\dpt(h)$,  as follows. First we put $h$ in $P$-good coordinates;
assume $h$ is represented using its Newton components as follows:
$$h=x_1^{N_1}x_2^{N_2}\dots x_d^{N_d}
z^\varepsilon \prod_{q=1}^{r}\prod_{j=1}^{u(q)}f_{j}^q (\bdx,z) u(\bdx,z),\qquad\text{
with } u(\mathbf 0,0)\ne 0.$$
If $h=x_1^{N_1}x_2^{N_2}\dots x_d^{N_d}z^\varepsilon u(\bdx,z)$,
$\varepsilon=0,1$, then we will say that
$\dpt(h):=0$.
Otherwise
$$
\dpt(h):=\max\{\dpt(\bar{h}_{q,j})\}+1,
$$
where the maximum is taken over all  pull-backs $\bar{h}_{q,j}$ under the
Newton maps associated with  Newton components $f_{j}^q$ of $f$.
If $h$ appears in good coordinates and a choice $\eta$ of good coordinates has been given after each Newton mapping, then we define in the same way $\dpt_\eta(h)$.
\end{defini}

\begin{obs}
From the above discussion $\dpt(h)\in\bn$
and decreases under the Newton maps.
Moreover $h$
has non-degenerated Newton principal part and
$h\not=x_1^{N_1}x_2^{N_2}\dots x_d^{N_d}z^\varepsilon u(\bdx,z)$
($u$ unity) if and only if
$\dpt(h)=1$.
The finiteness of the invariant $\dpt(h)$ will be used in the
proof of the main result. It is easily seen that for any recursive choice $\eta$ of good coordinates, $\dpt_\eta(h)\leq\dpt(h)$. Note that the components which correspond to edges appearing after an edge with integer multi-slope need less Newton mappings to be \emph{improved}.
\end{obs}

We can keep all the informations on the Newton process in a tree, the same
way, as Eisenbud and Neumann diagrams do for curves. 
Assume  we have chosen a system of coordinates such that the Newton polygon of $h$ is a polygonal path. In our tree, each compact face will be represented by a vertex, and two vertices are connected by an edge (a vertical edge) if and only if the two faces
intersect. The non compact faces are represented by arrows connected to the face they intersect. The arrow $\mathcal F _{0,..,0}$ representing $x_i=N_i, i=1,...,d$ is decorated by $(N_1,...,N_d)$, and the arrow  $\mathcal F _{0}$, representing the hyperplane $z=0$ or $z=1$, is decorated by $0$ or $1$. The edges are decorated as follows: let $\gamma _q$ be a face (represented by the vertex $v_q$) 
and $\blambda _q =(...,b_l ^q/n_1 ^q,...)$ its multislope, then the extremity of the edge attached to $v_q$ in the direction of  $\mathcal F _{0}$, is decorated
by $n_1 ^q$, and the edge attached to $v_q$ in the direction of  $\mathcal F _{0,...,0}$ is decorated by $(b_1^q,...b_d^q)$. Notice that  the edge multideterminants are all strictly positive. Now if $f$ is non degenerated with respect to $\gamma _q $, we attached to $v_q$ as non vertical edges ending with  arrows as the number of integer points on $\gamma _q$ minus one (recall that this is nothing but $v(q)$. If not, we consider a Newton map attached
to $\gamma _q $ and one of its factors. We consider the pull back of $h$ under this map, and put it in coordinates such that the Newton polygon is a polygonal path. We consider the part of the tree corresponding to the new 
Newton polygon and we delete its arrow $\mathcal F _{0,..,0}$ and attached the edge to the vertex $v_q$, so that the edge is not vertical. Because after a finite number of steps we have something non degenerated,
at the end we have a tree whose edges are attached to vertices or arrows.
We can make this diagram minimal by first erasing the edges at the bottom of the tree ending with an arrow decorated with a $0$ (we erase also this arrow) and decorated with $1$, and then the vertex on the other end of the edge if it is connected only to two other vertices. In good coordinates the diagrams are always minimal. If we are in $P$-good coordinates, there are no $1$ on vertical edges. If we erase the decoration $1$ on edges, and forget about vertical or horizontal edges, the diagram in good coordinates or $P$- good coordinates are the same. 

\begin{defini} \label{candidate} Let $h$ be a quasi-ordinary power series with a recursive choice $\eta$ of good coordinates.
As usual suppose
$h=x_1^{N_1}x_2^{N_2}\dots x_d^{N_d}f(\bdx,z)u(\bdx,z)$
and assume $f$ is represented using its Newton components as:
$f=
z^\varepsilon\prod_{q=1}^{r}\prod_{j=1}^{v(q)}f_{j}^q$, $\varepsilon=0$ or $ 1.$
Let $\omega$ be a differential form such that
the pair $(h,\omega)$ verifies the support condition
(\ref{spc}).
According to (\ref{1connuevo}), (\ref{01con}),
(\ref{01con2}) and (\ref{01con3}) we consider
a set $CP_\eta(h,\omega)$ of \emph{candidate poles}
for 
$(h,\omega)$ recursively
as follows:

$$
CP_\eta(h,\omega):=\left\{(1,1)^\varepsilon,(N_i,\nu_i)\right\}_{i=1}^d\cup \bigcup_{q=1}^r
\left\{\left(\frac{M_l^q}{c_l^q}, \nu_lp^q_l+\bar{b}_l^q\right)
\right\}_{l=1}^d\cup\, \bigcup\, CP_\eta({\bar{h}_{q,j}},{\bar{\omega}_{q,j}}),
$$
where the last union runs over all  pull-back
$(\bar{h}_{q,j},{\bar{\omega}_{q,j}})$ of $(h,\omega)$ under
the Newton maps associated with Newton components $f_{j}^q$ of $f,$
$q=1,\ldots,r$ and $j=1,\ldots,v(q)$. We will drop $\eta$ if $P$-good coordinates have been chosen.
\end{defini}

\begin{obs} Let $h\in k[[\bdx]][z]$ be a quasi-ordinary function
with $h(\mathbf{0},0)=0$
such that
$
D_z(h)=x_1^{\alpha_1}x_2^{\alpha_2}\dots x_d^{\alpha_d}\varepsilon(\bdx),
$
where $\varepsilon\in k[[\bdx]],\, \varepsilon(0)\not=0.$
After relabeling the variables $x_i,$ one of the following conditions
occurs.
Either $\lgt(h)=0,$ therefore $\alpha_i=0$ for all $i=1,\ldots,d$
(which is is equivalent to $h=z^\varepsilon$, $\varepsilon=0,1$, in good coordinates) or
%
$\alpha_k>0$ for all $k=1,\ldots,m\leq d$
and $\alpha_k=0$ for all $k\in\{m+1,\ldots,d\}.$
In such a case the last $d-m$ coordinates
of any characteristic exponents of $f$ are always zero.
In particular compact faces of $\Gamma(h)$
are contained in the $(m+1)$-dimensional coordinate
plane $x_{m+1}=\ldots=x_d=0$.
The Newton maps never involve the last $d-m$ coordinates
and therefore all  compact faces of the new
Newton polyhedra just only depend on the pull-back under the Newton maps
of the coordinates $x_1,\ldots,x_m,z$.
\end{obs}

\begin{lema} The set of candidate poles does not depend on $\eta$,
i.e., it equals always the result obtained for $P$-good coordinates.
\end{lema}

This is a consequence of the fact that the candidate poles can be read
from the tree we have defined; note that the multiplicities of the differential form can be deduced from the decorated tree. This tree is the same for good coordinates and $P$-good coordinates.




\begin{defini}
We will say that two quasi-ordinary power series with recursive choices $\eta,\eta'$ of good coordinates
\emph{have the same $(\eta,\eta')$-characteristic exponents} if their graphs are equal. We will drop the term $(\eta,\eta')$ if 
$P$-good coordinates is the common choice or if $(\eta,\eta')$
are clear in the context.
\end{defini}

Next lemma states the behaviour of these graphs under
the Newton mappings.

\begin{lema}\label{charexp}
Let $f\in k[[\bdx]][z]$ be a quasi-ordinary power series
with a recursive choice $\eta$ of good coordinates. Let $g\in
k'[[\bdx]][z]$ be 
other quasi-ordinary power series 
(maybe defined over distinct algebraically closed fields of characteristic
zero $k'$) with a recursive choice $\eta'$. Then $f$ and $g$ have the same $(\eta,\eta')$-characteristic exponents if and only if
\begin{enumerate}[\rm(1)]
\item $\Lambda_{ND}(f)=\Lambda_{ND}(g)$,
\item there exists a bijection between the roots of the polynomials
$f_{\blambda_q}$ and $g_{\blambda_q}$ and
\item under this bijection the Newton component $f_j^q$ (or the root
$f_{q,j}$) corresponds to the Newton component $g_j^{q}$ (or the root
$g_{q,j}$) then the quasi-ordinary power series $\bar{f}_{q,j}$ and
$\bar{g}_{q,j}$ have the same  $(\eta,\eta')$-characteristic exponents.
\end{enumerate}
\end{lema}

It is possible to prove that the data $\eta$-characteristic exponents, that is the 
tree defined above, is the same as the characteristic exponents of each irreducible 
component of $f$ and the order of coincidence
between distinct irreducible factors of $f$ (as in the case of plane curves).

\begin{obs} \label{Achange}
Under any linear change of coordinates of type $z=\beta_{d+1} z_1$,
$x_l=\beta_l y_l$, $l=1\ldots,d$, with
$\beta_1,\ldots,\beta_{d+1}\in {\mathbf G}_{m,k},$ the characteristic exponents,
resp. monomials, of $f$ do not change.
\end{obs}

\subsection{Transversal sections of a quasi-ordinary power series}
\mbox{}\label{transversal}

Consider a quasi-ordinary 
$h\in k[[\bdx]][z]$ with
$h(\bdx,z)=x_1^{N_1}x_2^{N_2}\dots x_d^{N_d}g(\bdx,z)$,
where no $x_i$ divides $g(\bdx,z)$
and $N_l\geq 0$ for any $l=1,\ldots,d.$
As usual we write
$h(\bdx,z)=x_1^{N_1}x_2^{N_2}\dots x_d^{N_d}f(\bdx,z)u(\bdx,z),$
$f(\bdx,z)$ being the quasi-ordinary Weierstrass $z$-polynomial of degree
$n$
in $k[[\bdx]][z]$
whose roots are centred at $({\mathbf 0},0)$
 and $u({\mathbf 0},0)\not=0$.
Assume that $h$ is given in good coordinates. 
Suppose that a pair $(h,\omega)$
verifying the support condition ~\ref{spc} is given,
where
$\omega=\prod_{j=1}^d x_j^{\nu_j-1}
d x_1\wedge\ldots\wedge d x_d\wedge dz, \,\nu_j\geq 1.$

Fix $i\in\{1,\ldots,d\},$ and suppose that $d\ge 2$. Let $k_i$ denote
an algebraic closure of the fraction field
of the formal power series ring $k[[x_i]].$
The $i$-th transversal section of $h$ will be
a finitely many  disjoint union
of $k_i$-quasi-ordinary formal power series obtained from
$h$ and their corresponding differential forms
verifying the support condition ~\ref{spc}.
Let ${\mathbf 0}_i$ denote the $d$-tuple which
has all coordinates but the $i$-th one, which is $x_i,$ equal $0.$
Consider the polynomial $f({\mathbf 0}_i,z)$ over $k_i$ of degree $n:$
$$f({\mathbf 0}_i,z)=z^n+\text{ lower degree terms }
=\prod_{m=1}^{v_i}(z-\alpha_m)^{\delta_m}\in k_i[z],$$
$\alpha_m\in k_i$ and $n=\sum \delta_m.$
Since the discriminant specializes,
for each root $\alpha_m$ of $f({\mathbf 0}_i,z)=0$ in $k_i[z],$
Jung-Abhyankar Theorem shows that
there exists a $k_i$-quasi-ordinary formal power series
$h_i^{\alpha_m}\in k_i[[{\widehat \bdx}_i]][z]$
centred at $({\widehat {\mathbf 0}}_i,\alpha_m),$
 ${\widehat \bdx}_i$ being all variables but $x_i.$

\begin{defini}
The $i$-\emph{transversal section} of $h$
consists in the set $\{h_i^{\alpha_m}: m=1,\ldots,v_i\}$
of $k_i$-quasi-ordinary formal power series,
 $h_i^{\alpha_m}$
centred at $({\widehat {\mathbf 0}}_i,\alpha_m).$
It is clear that $\lgt(h_i^{\alpha_m})<d.$
If the root $\alpha_m\neq 0$ then we need a translation
of type $z_1=z-\alpha_m$ to describe the corresponding $h_i^{\alpha_m}.$
On the other hand, if $x_i$ is not an essential variable for $h$
then the only root
is $\alpha=0.$

Since we are also interested
in keeping the differential form, the
$i$-\emph{transversal section} of a pair $(h,\omega)$
consists in the set of  finitely many pairs $(h_i^{\alpha_m},\omega_i)$
where $\omega_i=\prod_{j=1,j\ne i}^d x_j^{\nu_j-1}
d x_1\wedge\ldots\wedge {\widehat{d x_i}} \wedge \ldots\wedge d x_d\wedge dz$.
\end{defini}

\begin{ejem}
Consider $f(x_1,x_2,z)=((z^2-x_1^3)^2+x_1^7)^3+x_1^{25}x_2^3.$
For $i=1,$ the $1$-transversal sections are the $k_1$-quasi-ordinary
power series centred at the roots
of $((z^2-x_1^3)^2+x_1^7)^3=0$ in $k_1.$
In fact these roots can be found
using the Newton-Puiseux algorithm
in $k_1$ for the previous polynomial.
For $i=2,$ the $2$-transversal section is
the $k_2$-quasi-ordinary
power series $f_2^{0}=((z^2-x_1^3)^2+x_1^7)^3+(x_2^3)x_1^{25}\in
k_2[x_1,z]$ centred at the root $\alpha=0$.
In fact it is clear that the Newton polyhedron
of $f_2^{0}$ is the projection
over the plane $x_2=0$ of
$\Gamma(f)$. 
\end{ejem}

A description of the components $h_i^{\alpha}$ of the $i$-transversal
section is as follows.
Assume that
\begin{equation}\label{descrip1}
f(\bdx,z)=\prod_{j=1}^{n}
\left(z-\zeta_j(x_1^{1/m},\ldots,x_d^{1/m})\right),
\end{equation}
with
$\zeta_j({\mathbf 0})=0.$
If $\zeta_j\ne 0$, since $f$ is given in good coordinates,
 then
\begin{equation*}
\zeta_{j}(\bdx)=\alpha_1^{(j)}\bdx^{\blambda_{(j)}}
+\sum_{\blambda_{(j)}<\blambda}\alpha_{\blambda}^{(j)}
\bdx^{\blambda},\,\,(\alpha_1^{(j)}\ne 0).
\end{equation*}

The characteristic exponents of type $({\mathbf 0},\frac{b_i}{n_1},{\mathbf 0})\in \bq^d,$
correspond to $1$-dimensional faces in $CF(f)$ 
contained in the $(x_i,z)$-plane.
We write in a unique way
$\zeta_j(\bdx)$ as sum of two formal
power series $p_j(x_i^{1/m})+t_j(\bdx)$
such that $p_j(x_i^{1/m})\in k[[x_i^{1/m}]]$
is the sum of monomials of $\zeta_j$
depending only on $x_i$ and $t_j=\zeta_j-p_j(x_i) .$ It is clear that
$t_j({\mathbf 0}_i)=0.$

If $\alpha\in k_i$
is a root of $f({\mathbf 0}_i,z)=0$,
$h_i^{\alpha}({\widehat \bdx}_i,z)$ is 
\begin{equation}\label{descrip2}
h_i^{\alpha}=
u(\bdx,z)\prod_{l=1}^d x_l^{N_l}
\prod_{j:p_j(x_i^{1/m})=\alpha}
\left(z-(\alpha+
t_j(x_1^{1/m},\ldots,{\widehat {x_i}},\ldots,x_d^{1/m})\right)\in
k_i[[{\widehat \bdx}_i]][z],
\end{equation}
where $t_j$ is seen as an element in
$k_i[[x_1^{1/m},\ldots,{\widehat {x_i}},\ldots,x_d^{1/m}]].$
The last product gives a decomposition of $h_i^{\alpha}$
in its irreducible roots too, because each factor has
$z$-degree one.
To work with $h_i^{\alpha}$ we first do, if necessary, the translation $z_1=z-\alpha,$
$\alpha\in k_i.$
Such a translation is a composition of some Newton maps
of $h$ associated with their corresponding
Newton components (all of them of type $({\mathbf 0},a_i/w_i,{\mathbf 0})$) and some
change of variables for having good coordinates.
The Newton map $\pi_1$ associated with a
Newton component of $h$ of type $({\mathbf 0},b_i/n_1,{\mathbf 0})$ is given by
$x_i=y_i^{p_i},\, z=(z_1+\alpha)y_i^{{\bar b}_i}$ and $x_j=y_j,$ if $j\ne i.$
After $\pi_1$, we need a translation                                                                                     
of type $z_1=\tilde z-\phi(\bdx)$ to get good coordinates.
This translation does not change
the position of the $i$-th variable. If in order to get $\alpha$ maybe we will
need more Newton maps of type  $({\mathbf 0},a_i/w_i,{\mathbf 0})$ and 
we finish when no more characteristic exponent of the mentioned type
appear to get $\alpha$. 
Therefore the composition of
these Newton maps and their corresponding translations, which we
denote by $\pi,$  is the same as
the Newton-Puiseux algorithm
to find the root $\alpha$ of $f({\mathbf 0}_i,z)$ in $k_i[z].$
The map $\pi$ does nothing on the set of variables different from $i$-th one.
Let $\tilde f$ (resp. ${\tilde h}$) denote the pull-back of $f$ (resp. $h$)
under the map $\pi$.
Up to a factor $x_i^{k/m},$ 
${\tilde f}(\bdx,z)=\displaystyle\prod_{j=1}^{n}
\left(z-t_j(x_1^{1/m},\ldots,x_d^{1/m}\right)$ and
\begin{equation*}
h_i^{\alpha}({\widehat \bdx}_i,z)={\tilde u}\prod_{l=1}^d x_l^{N_l}u(\bdx,z)\prod_{j:p_j(x_i^{1/m})=\alpha}
\left(z-
t_j(x_1^{1/m},\ldots,{\widehat {x_i}},\ldots,x_d^{1/m})\right).
\end{equation*}
Moreover, $\omega_i$ and the pull-back $\pi^* \omega$ are the same,
up to a constant,
as differential forms over $k_i.$

Once this translation is done, the Newton polyhedron
$\Gamma_i(h_i^{\alpha})$
of $h_i^{\alpha}$ as $k_i$-quasi-ordinary power series is the
projection over the hyperplane $x_i=0$ of the Newton polyhedron
of the  pull-back $\tilde f$; this is clear for $\alpha=0$
and for the other case is similar. In particular $\Gamma_i(h_i^{\alpha})$
is a monotone polygonal path. In fact,
the $k_i$-quasi-ordinary power series $h_i^{\alpha}$
is in good coordinates because the condition not to be a power of a linear form
is generic. 
\begin{obs}
Note that if $h$ is given in $P$-good coordinates, we cannot deduce that
$h_i^{\alpha}$ is also in $P$-good coordinates. 
This is the main reason for introducing the more general concept of good coordinates.
\end{obs}

\begin{ejem}
(1) Consider $f(x_1,x_2,z)=(z^2-x_1^2x_2^5)(z^2-x_1^2x_2^7).$
In principle $\Gamma(f)$
has two edges which are projected
over just only one edge if we project over $x_2=0.$
The $2$-transversal
section has at least two different components
with the same characteristic exponent.
(2)  $f(x_1,x_2,z):=(z^2-x_1^2x_2^5)^2+x_1^4x_2^{17}.$
Seeing this as a $2$-transversal
section one has $f_2^{0}=z^4-(2x_2^5)x_1^4z+(x_2^{10}+x_2^{17})x_1^4.$
In this case some more monomials of $x_1^4$
have appeared. These monomials are hidden in $ND(f)$.
\end{ejem}

To sum up:

\begin{prop} \label{transsection} With the above notations,
\begin{enumerate}[{\rm(1)}]
\item The formal power
series $h_i^{\alpha_m}\in k_i[[{\widehat \bdx}_i]][z]$
centred at $({\widehat {\mathbf 0}}_i,\alpha_m),$
defines a quasi-ordinary power series with $\lgt(h_i^{\alpha_m})<d.$

\item The Newton process for $h$ induces Newton process for
$h_i^{\alpha_m}.$  Each step of the Newton process for
$h$ is necessary for at least one $h_i^{\alpha_m}$, $i\in\{1,\ldots,d\},$
and $m\in\{1,\ldots,v_i\}.$

\item Fix $i\in\{1,\ldots,d\}$ and $m\in\{1,\ldots,v_i\}.$ Then

\begin{enumerate}[\rm i)]
\item If $\alpha_m=0$ then the quasi-ordinary power series
$h_i^{0}$ is in good coordinates.

\item At each step of the Newton process for $h,$
$\Gamma(h_i^{0})$ is the projection of $\Gamma(h)$ on the hyperplane $x_i=0$.

\item If $\alpha_m\ne 0,$
then some steps of the Newton process for $h$ firstly provide 
the translation $z=z_1-\alpha$
we need to put $\alpha$
at zero. After that we apply $i)$ and $ii).$

\item Every pair $(h_i^{\alpha_m},\omega_i)$
verifies the support condition (\ref{spc}). 

\end{enumerate}
\item In fact
\begin{equation}
\label{cpp}
CP(h,\omega)\subset \bigcup_{i=1}^d  \bigcup_{m=1}^v
CP(h_i^{\alpha_m},\omega_i),
\end{equation}
Recall that the second union runs over all 
roots of $f({\mathbf 0}_i,z)=0$ in $k_i[z].$
\end{enumerate}
\end{prop}

The only fact in the proposition that
is not yet proved is the equality (\ref{cpp}).
The set $CP(h,\omega)$ of \emph{candidate poles}
for
$(h,\omega)$ can be written as
$CP(h,\omega)=CP_1(h,\omega)\cup\, \bigcup\,
CP({\bar{h}_{q,j}},{\bar{\omega}_{q,j}}),
$ see definition \ref{candidate},
where
the set $CP_1(h,\omega)$
is completely described from $ND(h)$ and $\omega$.
Since this monotone polygonal path
is determined by (and determines) its projections over all
hyperplanes $\{x_i=0\}, i=1,\ldots,d,$ then $CP_1(h,\omega)=
\bigcup_{i=1}^d  CP_1(h_i^{0},w_i).
$
In the previous identity we forget the roots which are not centred at
zero because some Newton maps are needed to find them.

In general, (\ref{cpp}) follows
from the fact that candidate poles are determined
by the characteristic monomials of the quasi-ordinary power series.
The characteristic monomials are determined by
difference of its roots.
In particular for each of the $i$-transversal section
the difference of two  of its roots in
$k_i[[x_1^{1/m},\ldots,{\widehat {x_i}},\ldots,x_d^{1/m}]]$ can be computed
as the specialization
of the difference of the corresponding roots
in $k[[x_1^{1/m},\ldots,{x_i}^{1/m},\ldots,x_d^{1/m}]].$

In fact we can say more.
Fix a characteristic exponent of $f$
and a Newton component with such a characteristic
exponent. Fix $i\in\{1,\ldots,d\}.$
If the characteristic exponent is $(\mathbf 0,b_i/n_1,\mathbf 0)$
then we have already described what is going on with the
Newton map $\pi$ associated with the Newton component.
Otherwise, reading the characteristic exponents and the roots
from equations
(\ref{descrip1}) and (\ref{descrip2}),
the characteristic exponents of the $i$-transversal section
of $f\circ \pi$ are the same as those one for
the pull-back under the corresponding Newton map
for the $i$-transversal section.
In particular   at each step
the Newton polyhedron of the $i$-transversal section
is the projection of the Newton polyhedron of $h$.

The fact that we get the same differential form, up to a constant factor,
in $x_i$ is clear because this is a general fact no related with
quasi-ordinary formal power series.

\section{Denef-Loeser motivic zeta function under the Newton maps}

We keep the field $k$ being algebraically closed and denote an algebraic
closure of the quotient field $k((t))$ of the domain $R:=k[[t]]$ by $K.$
The goal is to compute zeta functions
for a $k$-quasi-ordinary power series using induction on its depth.
In the last section we have computed the $A$-part.
In order to compute the $B$-part we must take into account that  a 
quasi-ordinary power series may be degenerated
with respect to its Newton polyhedron. The main idea is to use 
Newton maps to measure some sets of arcs in terms of the 
quasi-ordinary power series issued after these Newton maps.
The problem is the existence of some arcs which cannot be lifted
under such maps.
We need to combine Newton maps and  $k[t]$-morphisms
to solve this problem.
Therefore we introduce technical objects,
called $\ww$ and $\tww$-quasi-ordinary series
in order to have  families of power series closed under Newton maps
to deal with.
We will compute the Denef-Loeser motivic zeta functions for these series following the ideas of subsections \ref{sec-mot-newton}
and \ref{dualsection}; the computation of the $A$-part is similar to the one in \ref{dualsection} but there are some small differences
and we will compute inductively the whole zeta function.

\begin{defini} Let
$h^t(\bdx,z)\in R[[x_1,\dots,x_d]][z]$
be a formal power series such that
$h^t=t^\theta x_1^{N_1}x_2^{N_2}\dots x_d^{N_d} f^t(\bdx,z) u^t(\bdx,z)$
where $x_j$ does not divide $f^t(\bdx,z)$ 
$\forall j=1,\ldots,d$, $t$ does not divide $f^t(\bdx,z) u^t(\bdx,z)$, $\theta\in\bn$ and $u^t(\mathbf 0,0)\neq 0$. 
We say that $h^t$ is
\emph{$\ww$-quasi-ordinary} if
its $z$-discriminant is
\begin{equation} \label{discriminant}
D_z(f^t)=t^\beta x_1^{\alpha_1}x_2^{\alpha_2}\dots
x_d^{\alpha_d}\varepsilon(t,\bdx),
\end{equation}
where $\varepsilon(0,{\mathbf 0})\not=0$; therefore if we consider
$t$ as a variable,
$f^t\in k[[t,\bdx]][z]$
is a  $k$-quasi-ordinary power series.
Let $J\subset\{1\dots,d\}$ and let $J'$ be its complement; we say that $h^t$ is $J$-bounded if its roots $\zeta^t$ verify
\begin{equation}\label{acotado}
\zeta^t\in k[\bdx_{J}^{1/m}][[\bdx_{J'}^{1/m},t]], \text{ for some m.}
\end{equation}
\end{defini}

\begin{obs}\label{translat} The main property of a $J$-bounded $\ww$-quasi-ordinary series $h^t$ is that $\forall\bdy_J^0\in\bc^d$ such that the entries in
$J'$ are zero the series $h^t(\bdx+\bdy_J^0,z)$ is a well-defined $\ww$-quasi-ordinary series.
\end{obs}

Let $f(\bar{\bdx},\bar{z})\in k[\bar{x}_1,\dots,\bar{x}_d][\bar{z}]$
be a quasi-ordinary polynomial
whose roots only contain
the (finite set of) characteristic monomials
of $f^t$ as $K$-quasi-ordinary power series
in a new set of variables $\bar{x}_1,\dots,\bar{x}_d,\bar{z}$.
Let $h:=\bar{x}_1^{N_1}\dots \bar{x}_d^{N_d} f(\bar{\bdx},\bar{z}).$
In particular the $K$-quasi-ordinary power series $h^t$ and the
$k$-quasi-ordinary power series $h$
have the same characteristic exponents and the same
characteristic monomials.

If $f^t(\bdx,z)\in k[[t]][[\bdx]][z]$
is not in \emph{good} coordinates, by Remark \ref{bueno},
the change of coordinates
to put $f^t(\bdx,z)$ in good
coordinates is an automorphism of
$k[[t]][[\bdx]][z]$
of type $z\mapsto z+m^t(\bdx)$,
$m\in k[[t]][[\bdx]]$. In particular
the $t$-variable does not change and condition (\ref{discriminant})
on the discriminant is preserved.
It is easily
seen that if $h^t$ is $J$-bounded then condition (\ref{acotado}) is also preserved because 
$m^t(\bdx)$ comes from some monomials of the roots of $h^t$.
Then we assume $f^t$ and $f(\bar{\bdx},\bar{z})$
are in good coordinates.

Consider differential forms
$\omega=\prod_{j=1}^d x_j^{\nu_j-1}d\bdx\wedge dz$
and $\bar{\omega}=\prod_{j=1}^d \bar{x}_j^{\nu_j-1}d\bar\bdx\wedge d\bar{z},$ with $\nu_j\geq 1$.
From now on
we assume that
$(h,\bar{\omega})$ verifies the support condition~\ref{spc}.
Let $X:=\ba_k^{d+1}$ and set
$$
V_{n,m}^t:=\{\barphi\in \cl_0(X): \,
\ord(h^t\circ \barphi)=n+\theta,
\ \ord(\omega\circ\barphi)=m\};
$$
define $V_{n,m}$ in the same way for $h$
(forgetting $\theta$).
The sets $V_{n,m}^t$ are measurable because of the
support condition. We define the power series
\begin{equation*}
\begin{split}
Z_{DL}^{\ww}(h^t,\omega,T)&:=T^\theta \sum_{n\in\bn} \left(\sum_{m\in\bn}
\bl^{-m}\mu_X(V_{n,m}^t)\right)T^n\in \widehat{\mm}_k[[T]] ,\\
Z^{\ww}_{\text{top},0}(h^t,\omega,s)&:=
\chi_{\text{top}}(Z^{\ww}_{DL}(h^t,\omega,\bl^{-s})).
\end{split}
\end{equation*}
Recall that
\begin{equation*}
\begin{split}
Z_{DL}(h,\bar{\omega},T)=\sum_{n\in\bn} \left(\sum_{m\in\bn}
\bl^{-m}\mu_X(V_{n,m})\right)T^n,\\
\zlo(h,\bar{\omega},s)=
\chi_{\text{top}}(Z_{DL}(h,\bar{\omega},\bl^{-s})).
\end{split}
\end{equation*}

We know from the results of Denef and Loeser that
$Z_{DL}(h,\bar{\omega},T)$ is a rational function.
The last part of the section is devoted to prove  the
rationality of $Z_{DL}^{\ww}(h^t,\omega,T)$ and to provide
a \emph{small set} of candidate poles for such functions, see
definition \ref{candidate}.
Recent interesting results of J. Sebag \cite{se:02}
prove the rationality of this power series.

\begin{thm}\label{thm-polos}
If $h^t$ is $\ww$-quasi-ordinary then the zeta function
$$Z_{DL}^{\ww}(h^t,\omega,T)\in\bz\left[\bl,\bl^{-1},
\left(1-\bl^{-\nu}T^N\right)^{-1}\right][T]
$$
where $(N,\nu)\in CP(h,\bar{\omega}).$
Moreover, $Z^{\ww}_{\text{\rm top},0}(h^t,\omega,s)$ exists.
\end{thm}

We need to define a subfamily of $\ww$-quasi-ordinary series
in order to obtain more information on the local topological zeta function.

\begin{defini} A symbol $\tww:=(p;g_1,\dots,g_d)$, where
$p\in \bp,g_i\in\bn$, $i=1,\dots,d,$ is called a \emph{weight}.
The semigroup $\Gamma_\tww\subset\bn^d$ associated with $\tww$
is the additive set of $\balpha\in\bn^d$ such that
$$
w_{\tww}(\balpha):=\sum_{i=1}^d \frac{\alpha_i g_i}{p}\in\bn.
$$
\end{defini}

\begin{defini} Let
$h^t(\bdx,z)\in R[[x_1,\dots,x_d]][z]$
be a formal power series such that
$h^t=t^\theta x_1^{N_1}x_2^{N_2}\dots x_d^{N_d} f^t(\bdx,z) u^t(\bdx,z)$
with $x_j$ does not divide $f^t(\bdx,z)$ 
$\forall j=1,\ldots,d$, $u^t(\mathbf 0,0)\neq 0$,
$t$ does not divide $f^t(\bdx,z) u^t(\bdx,z)$ and
 $\theta\in\bn$.
We call it $\tww$-quasi-ordinary if it is
$\ww$-quasi-ordinary and all  monomials
in $h^t$ are of type $t^{w_\tww(\balpha)}\bdx^{\balpha} z^n$, where $\balpha\in\Gamma_\tww$ (and some non-zero coefficient in $k$).
\end{defini}

\begin{obs} If $h^t$ is $\tww$-quasi-ordinary then the power
series $h^1:=h^t_{|t=1}$ is $k$-quasi-ordinary and have the same characteristic exponents as $h^t$.
If $h^t$ is $\tww$-quasi-ordinary,
it is true also for $f^t$ because
the $d$-tuple $(N_1,\dots,N_d)\in\Gamma_\tww$ and
$\theta=w_\tww(N_1,\dots,N_d)$. Moreover, if we put $h^t$ in good coordinates, the new series is also $\tww$-quasi-ordinary.
\end{obs}

\begin{obs} Let $h^t$ be a $\tww$-quasi-ordinary power series, $\tww:=(p;g_1,\dots,g_d)$. Then, if $J:=\{j\mid
1\leq j\leq d,\ g_j\neq 0\}$ then $h^t$ is $J$-bounded.
\end{obs}

\begin{defini} Let $h^t$ be a $\tww$-quasi-ordinary power series.
We say that $h^t$ is \emph{Newton compatible}
if the number of irreducible factors
(different from $t$) of $h^t_{|ND(h^t)}$
in $R$ equals the number of irreducible factors of $h^1_{|ND(h^1)}$ in $k$.
For instance $z^2-tx^2y^2$ is not Newton compatible and
it is $\tww$-quasi-ordinary for some adequate weight.
\end{defini}

\begin{thm}\label{thm-top}
If $h^t$ is $\tww$-quasi-ordinary and Newton compatible then
$$
Z^{\ww}_{\text{\rm top},0}(h^t,\omega,s)=\zlo(h,\bar{\omega},s).
$$
\end{thm}

From now on we suppose that $h^t$ is $\ww$-quasi-ordinary;
we will indicate explicitly when we consider it as
$\tww$-quasi-ordinary.
The rest of this section is devoted to the proof of both theorems. We will focus on the proof of Theorem \ref{thm-polos} and we will point out the special arguments required for the proof of Theorem \ref{thm-top}. Since the results
do not depend on the factor $T^\theta$, we can omit it.
From now on we assume that $h^t$ (or $f^t$) is
given in good coordinates. The proofs are given by induction
on the depth.

\begin{paso}\label{start-ind}
{\bf $\dpt(h)=\dpt(h^t)=0$.}
\end{paso}
Then $ND(h)$ has only one compact face
which is $0$-dimensional,
$$Z_{DL}^{\ww}(h^t,\omega,T)=\bl^{-(d+1)}(\bl-1)^{d+1}
\frac{\bl^{-(1+\sum \nu_j)} T^{1+\sum N_j}}{1-\bl^{-1}T}\prod_{j=1}^d
\frac{1}{1-\bl^{-\nu_j}T^{N_j}},$$
and $Z_{DL}^{\ww}(h^t,\omega,T)=Z_{DL}(h,\bar{\omega},T)$.
Then $Z^{\ww}_{\text{\rm top},0}(h^t,\omega,s)=\zlo(h,\bar{\omega},s)$.

\begin{paso} Assume that $\dpt(h)>0$.
\end{paso}

We keep the notation of \S\ref{secnondeg}.
Let $S$ be defined by
$x_1\cdots x_d z=0,$ (or
$\bar{x}_1\cdots \bar{x}_d \bar{z}=0$);
we consider arcs in $\barphi\in\cl_0(X)\setminus\cl_0(S)$
and define $\bdk(\barphi)$ and $\bda(\barphi)$ as in \S\ref{secnondeg}.

The following easy remark will be a key point
to understand the dual decomposition
associated with $f^t(\bdx,z)\in R[[\bdx]][z]$.
Each $\barphi\in \cl_0(X)$ univocally defines
an arc
$\bpsi:=(t,\barphi)\in \cl_0(\ba_k^{1}\times X)$ and
$\ord_t(f^t\circ \barphi)=\ord_t(g\circ \bpsi),$
where $g(y,\bdx,z)\in k[[y,\bdx]][z]$ is the same function
as $f^t$ but $t$ is substituted by a new variable $y.$

\begin{paso} Newton polyhedron of $h\in k[\bar{\bdx},\bar{z}]$.
\end{paso}
This Newton polyhedron $\Gamma(h)$
and its dual decomposition
have been described in paragraph \ref{dualsection}.
Denote by $\bar{\gamma}_1,\dots,\bar{\gamma}_r$
its $1$-dimensional
faces with corresponding vertices
$\bar{\tau}_0,\bar{\tau}_1,\dots,\bar{\tau}_r$.
From now on we assume that $\bar{\tau}_r$
is the $\bar z$-highest vertex.
Thus we can write
$$f_{ND(f)}=z^\varepsilon \prod_{q=1}^r \prod_{{\bar j}=1}^{v(q)} (f_{q,{\bar j}})^{m_{q,{\bar j}}}
\,\text{ where }  f_{q,{\bar j}}:=
\bar{z}^{n_1^q}-\bar{\beta}_{\bar j}^q \bar{x}_1^{b_1^q}
\ldots \bar{x}_d^{b_d^q},$$
with $\varepsilon=0$ or $1$ depending if $I=I'$ or not.
Recall that in the dual space, the  reduced integer equation of the hyperplane
$l^q$ is
$\eta_q(v_1,\dots,v_d,v_{d+1})=0,$ see equation (\ref{eq-hyper}).
The intersection of these hyperplanes with $\br_{>0}^{d+1}$ are
the $d$-dimensional cones $\Delta_{\bar{\gamma}_q}$.
They determine the
$(d+1)$-dimensional cones $\Delta_{\bar{\tau}_q}$, $q=0,1,\dots,r$.
The cones associated with
compact faces of ${\Gamma(h)}$ (or $\Gamma(f)$) give a partition of
$\br_{>0}^{d+1}$ in the disjoint union
$\cup \Delta_{\bar \tau}.$

If $\barphi\in\cl_0(X)\setminus\cl_0(S)$ and
$\bdk:=\bdk(\barphi)$,
the order of the differential form $\bar \omega$
is $
m(\bdk):=\sigma_{\bar \omega}(\bdk)-\bdk.
$
If  $\bdk\in\bp^{d+1}$ is fixed,
the semialgebraic subset $V_{n,m}^\bdk$ is naturally defined. It
is empty unless $m=m(\bdk)$. Then
for each compact face $\bar \tau$
define
$$
Z_{DL}^{\bar \tau}(h,\bar \omega,T):=\sum_{\bdk\in\Delta_{\bar \tau}\cap\bp^{d+1}}
Z_{DL}^{\bdk}(T).
$$

\begin{paso} Description of the Newton polyhedron $\Gamma^t(f^t)$
of $f^t\in k[[t,\bdx]][z]$.
\end{paso}

The discriminant condition
(\ref{discriminant}) implies that $f^t$ is a $k$-quasi-ordinary
power series in good coordinates in $(d+2)$ variables.
Then its Newton polyhedron
$\Gamma^t(f^t)$ has only $0$ and
$1$-dimensional compact faces,
say $m$ edges $\tilde{\gamma}_1,\ldots,\tilde{\gamma}_m$
and the $m+1$ corresponding vertices.
We can write
\begin{equation}\label{t-newton1}
f^t_{ND(f^t)}=z^\varepsilon \prod_{p=1}^m \prod_{s=1}^{v(p)}(z^{\tilde n_1^p}-\beta_s^p t^{\tilde b_0^p}x_1^{\tilde b_1^p}\ldots x_d^{\tilde b_d^p})^{m_{p,s}},
\end{equation}
where $\gcd(\tilde n_1^p,\tilde b_0^p,\tilde b_1^p,\dots,\tilde b_d^p)=1$
and $\varepsilon=0$ or $1$ depending if  $I=I'$ or not.

The characteristic exponents associated with $\Gamma^t(f^t)$ are
of type $(\frac{\tilde b_0^p}{\tilde n_1^p},\blambda^t_p)\in \bq^{d+1}$
where $\tilde n_1^p\blambda^t_p=(\tilde b_1^p,\ldots,\tilde b_d^p)\in\bn^d$. We need also to define
$n_1^p,b_1^p,\dots, b_d^p$ such that  $\gcd(n_1^p,b_1^p,\dots,b_d^p)=1$ and $n_1^p\blambda^t_p=(b_1^p,\ldots,b_d^p)$.

Given $\blambda_q$ from a compact face $\bar \gamma_q$
of $\Gamma(h)$, $q=1,\dots,r$,
we collect the set $S_q$ of rationals $\frac{b}{n}$
such that there exists a compact face $\tilde\gamma_p$ of $\Gamma^t(h^t)$
whose associated characteristic exponent is
$(\frac{b}{n},\blambda_q)$; all characteristic exponents
of $h^t$ are obtained in this way.

Let $u(q)$ be the cardinality
of $S_q$; we denote by $\frac{\alpha_{q,j}}{n_{q,j}}$, $j=1,\dots,u(q),$
the elements of $S_q$, where $\alpha_{q,j}\in\bn$ and $n_1^{q,j}=\tilde n_1^p$ for the corresponding $p\in\{1,\dots,m\}$. We denote also
$s_1^{q,j}:=\frac{n_1^{q,j}}{n_1^q}\in\bn$.
Then identity (\ref{t-newton1}) can be rewritten as follows
$$
f^t_{ND(f^t)}=z^\varepsilon \prod_{q=1}^r \prod_{j=1}^{u(q)}
\prod_{w=1}^{w(q,j)}((z^{n_1^q})^{s_1^{q,j}}-\beta_w t^{\alpha_{q,j}}(x_1^{b_1^q}
\ldots x_d^{b_d^q})^{s_1^{q,j}})^{m_{q,j,w}},
$$
where $v(q)=\sum_{j=1}^{u(q)} s_1^{q,j} w(q,j)$
since $f^t$ and $f$ have by hypothesis the same characteristic
exponents.

\begin{obs} If $h^t$ is $\tww$-quasi-ordinary and
Newton compatible then $r=m$, it means that for each $q=1,\dots,r$ we have
$\# S_q=1$ and $S_q\in\bn$. In this case
$s_1^{q,j}=1$, $\tilde b_j^q=b_j^q$.
\end{obs}

In the dual space $(\br^{d+2})^*$, with coordinates
$(v_0,v_1,\ldots,v_d,v_{d+1}),$ the dual decomposition induced by $\Gamma^t(h^t)$
is given by dual cones
$\Delta_{\tilde{\gamma}_p}$ associated with the hyperplanes,
$\tilde l^p:\tilde b^p_0 v_0+\tilde b^p_1 v_1+\ldots+\tilde b^p_d v_d-\tilde n^p_1 v_{d+1}=0,$
 for $p\in\{1,\ldots,m\}$,
and
dual cones $\Delta_{\tilde{\tau}}$ associated with a vertex $\tilde{\tau}$.
We identify the affine hyperplane $H=\{v_0=1\}$ with
$(\br^{d+1})^*$ and $H\cap (\br_+^{d+2})^*$
with $(\br_+^{d+1})^*$ and
consider the decomposition induced
by $\Delta_{\tilde{\gamma}_p}$ and $\Delta_{\tilde{\tau}}$
in $(\br_+^{d+1})^*$.
This decomposition has the same properties as the initial one.
It has exactly $m$ ordered $d$-dimensional planes in $H$
and $m+1$ ``pieces'' of dimension $(d+1)$
between them.
Let us describe them better.

For $q=1,\ldots,r$, and
$j\in\{1,\dots,u(q)\}$ define the cone $\Delta_{{q,j}}^t$
as the intersection with $\br_{>0}^{d+1}$ of
the hyperplane $l^{q,j}_t$ with integral equation:
\begin{equation} \label{hyperplane}
s_1^{q,j}\eta_q(v_1,\dots,v_d,v_{d+1})+\alpha_{q,j}=0,\text{ see (\ref{eq-hyper}).}
\end{equation}
Observe that $m=\sum_{q=1}^r u(q)$.

The $(d+1)$-dimensional convex rational polyhedra
$\Delta_s^t$, $s=1,\ldots,m-1$,
are either the region $M^{q,j}_t$ (that we will call
of type $M$) contained between two parallel
hyperplanes $l^{q,j}_t$ and $l^{q,j+1}_t$
or the region $N^q_t$ (of type $N$)
contained between the hyperplanes $l^{q,u(q)}_t$ and $l^{q+1,1}_t$
for some $q=1,\ldots,m-1.$
The first one $\Delta_0^t$ and the last one $\Delta^t_m$
are
\begin{equation}
\begin{split}
\Delta_0^t:=\{v\in \br_{>0}^{d+1}:
s_1^{1,1}\eta_1(v_1,\dots,v_d,v_{d+1})+\alpha_{1,1}<0\},\\
\Delta_m^t:=\{v\in \br_{>0}^{d+1}:
s_1^{r,u(r)}\eta_r(v_1,\dots,v_d,v_{d+1})+\alpha_{r,u(r)}>0\}.
\end{split}
\end{equation}

Thus we fix the partition of
$\br_{>0}^{d+1}$ (identified with $H$)
as disjoint union of three distinct types of convex rational
polyhedra $\Delta_\tau^t$:
\begin{itemize}
\item $d$-dimensional cones  $\Delta_{q,j}^t$
corresponding to the hyperplanes $l_t^{q,j},$

\item convex rational
polyhedra of type $N$, and

\item convex rational
polyhedra of type $M$.
\end{itemize}

Following the conventions of this work the
 convex rational polyhedra of type $M$ or $N$ will be called 
 \emph{vertices} and anyone of the first type  an \emph{edge}.

Let $\barphi\in\cl_0(X)\setminus\cl_0(S)$ and let
$\bdk:=\bdk(\barphi)$.
The order of the differential form $\omega$ is
the same as for $\bar \omega$, i.e., 
$m(\bdk)=
\sigma_{\omega}(\bdk)-\bdk.$
Fix $\bdk\in\bp^{d+1}$ and define
$$
V_{n,m}^{\bdk,t}:=\{\barphi\in \cl_0(X)\setminus \cl_0(S ): \,
\bdk(\barphi)=\bdk,\ord(h^t\circ \barphi)=n+\theta,
\ \ord(\omega\circ\barphi)=m\};
$$
These sets are empty unless $m=m(\bdk)$, then
we define
\begin{equation*}
Z_{DL}^{\bdk,\ww}(T):=\sum_{n\in\bn}
\bl^{-m(\bdk)}\mu_X(V_{n,m(\bdk)}^{\bdk,t}) T^n.
\end{equation*}
For each convex rational polyhedron $\Delta_\tau^t$
in the previous
partition,  define
$$
Z_{DL}^{\Delta_\tau^t,\ww}(h^t,\omega,T):=\sum_{\bdk\in\Delta_\tau^t\cap\bp^{d+1}}
Z_{DL}^{\bdk,\ww}(T).$$

\begin{obs} Each hyperplane $l_q$ is parallel to the hyperplanes
$l_t^{q,j}$. The region associated
with a vertex $\bar\tau_q$ is a cone delimited
by two hyperplanes which are parallel to the
hyperplanes delimiting $N_t^q$ (different from the coordinate hyperplanes).
\end{obs}

\begin{obs} If $h^t$ is $\tww$-quasi-ordinary and Newton compatible then no region of type $M$ exists.
\end{obs}

\begin{paso}Newton polyhedron of $h^t$ as $K$-quasi-ordinary function.
\end{paso}
\label{K-newton}
This Newton polyhedron $\Gamma^K(h^t)$
of $h^t$ is
the projection
over the hyperplane ${t=0}$ of the Newton polyhedron
of $h^t$ and it coincides with $\Gamma(h)$.
For $q\in\{1,\ldots,r\}$,
we define
\begin{equation*}
{\widehat f}^t_{\blambda_q}:=\prod_{j=1,\alpha_{q,j}\in S_q}^{u(q)}
{\widehat f}^t_{q,j},\,\text{ where }
{\widehat f}^t_{q,j}:=\prod_{w=1}^{w(q,j)}((z^{n_1^q})^{s_1^{q,j}}-{\beta}_w t^{\alpha_{q,j}}(x_1^{b_1^q}
\ldots x_d^{b_d^q})^{s_1^{q,j}})^{m_{q,j,w}}.
\end{equation*}

\medskip
\subsection{Vertices of the dual decomposition}
\mbox{}
\medbreak

In the following, the  aim is to compute $
Z_{DL}^{\Delta_\tau^t,\ww}(h^t,\omega,T)$ for $\tau$ a vertex of type $M$
or $N$ and compare
it with the computations we have already done for
 $Z_{DL}^{\bar \tau}(h,\bar \omega,T),$ see (\ref{eqA}).

Take a convex rational polyhedron $\Delta_\tau^t$
of the dual decomposition of $\br^{d+1}_+$
which belongs to
any of the two types $M$ or $N$.
It corresponds to
a monomial $t^\delta \bdx^{\balpha} z^n$ whose support  is the intersection
of two $1$-dimensional compact faces of $\Gamma^K(h^t).$
The cone $\Delta_\tau^t$ is
the positive region delimited by two inequalities, say
\begin{equation*}
\begin{split}
\tilde\eta_1(v_1,\dots,v_d,v_{d+1})+\tilde\alpha_{1}>0,\\
\tilde\eta_{2}(v_1,\dots,v_d,v_{d+1})+\tilde\alpha_{2}<0,
\end{split}
\end{equation*}
where $\tilde\eta_1$ and $\tilde\eta_2$ can define parallel (or not) hyperplanes.
Let $\barphi\in\cl_0(X)\setminus\cl_0(S)$ and let
$\bdk:=\bdk(\barphi)$ and
$n_\tau(\bdk):=
\alpha_1 k_1+\dots+\alpha_d k_d + n k_{d+1}$.
We have
$
\ord(h^t\circ\barphi)=n_\tau(\bdk)+N(\bdk)+\delta,$
with
$N(\bdk):=\sum_{j=1}^{d} N_j k_j$
(we are forgetting the exponent $\theta$).

\begin{lema}\mbox{}
\label{2002}
If $\bdk\in\Delta_{\tau}^t,$ then
$$
\bl^{-m(\bdk)}\mu_X(V_{n_\tau(\bdk)+\delta +N(\bdk),m(\bdk)}^{\bdk,t})=
\bl^{-(d+1+\sigma_\omega(\bdk))}(\bl-1)^{d+1}.
$$
\end{lema}

\begin{proof}
The measure of the cylindrical sets can be computed as in Step 2
in~\S~\ref{secnondeg}. Denote $n_\tau(\bdk)+\delta+N(\bdk)$ by  $\bar n$;
hence $\mu_X(V_{{\bar n},m(\bdk)}^{\bdk,t})=
[\pi_{{\bar n}}
(V_{{\bar n},m(\bdk)}^{\bdk,t})]
\bl^{-({\bar n}+1)(d+1)}.
$
It comes from a monomial since it is a vertex then
$$
[\pi_{{\bar n}}
(V_{{\bar n},m(\bdk)}^{\bdk,t})]
=(\bl-1)^{d+1}
\bl^{(d+1){\bar n}-(k_1+\ldots+k_{d+1})}.
$$
\end{proof}

Hence the contribution in terms
of generating functions of convex rational polyhedra is written as
\begin{equation*}
\begin{split}
Z_{DL}^{\Delta_\tau^t,\ww}=&T^\delta (\bl-1)^{d+1}
\bl^{-(d+1)}
\sum_{\bdk\in\Delta_\tau^t\cap \bp^{d+1}}
\bl^{-\sigma_\omega(\bdk)}
T^{n_\tau(\bdk)+N(\bdk)}=\\
=& T^\delta (\bl-1)^{d+1}
\bl^{-(d+1)} \Phi_{\Delta_\tau^t}(\bdy),
\end{split}
\end{equation*}
where $\bdy:=(\bl^{-\nu_1}T^{\alpha_1},\ldots,\bl^{-\nu_d}T^{\alpha_d},
\bl^{-1}T^{n}).$

 For any
integer linear form $\eta(v_1,\dots,v_d,v_{d+1})$ and any
$\alpha\in \bn$ define
$\Delta_\eta^\alpha:=\bp^{d+1}\cap \{\eta+\alpha=0\}$
and $\Delta_\eta^{c,\alpha}:=\bp^{d+1}\cap \{\eta+\alpha<0\}.$
The indicator functions verify the identity:
$$
[\Delta_{\eta}^{c,\alpha}]=
[\Delta_{\eta}^{c,0}]-\sum_{0<\alpha'\leq \alpha}
[\Delta_{\eta}^{\alpha'}].
$$
In paragraph \ref{dualsection} we have already used $\Delta_{\eta}^{c,0}$
without the superscript $0.$
Since
$
[\Delta_\tau^t\cap \bp^{d+1}]=[\Delta_{\eta_2}^{c,\alpha_2}]-\left(
[\Delta_{\eta_1}^{c,\alpha_1}]+[\Delta_{\eta_1}^{\alpha_1}]\right)
$
then
\begin{equation}\label{indicator}
[\Delta_\tau^t\cap \bp^{d+1}]=
\left([\Delta_{\eta_2}^{c,0}]-[\Delta_{\eta_1}^{c,0}]\right)-
\left(\sum_{0<\alpha\leq \alpha_2}
[\Delta_{\eta_2}^{\alpha}]-\sum_{0<\alpha< \alpha_1}
[\Delta_{\eta_1}^{\alpha}]\right).
\end{equation}

If the corresponding hyperplanes are parallel,
then $[\Delta_\tau^t\cap \bp^{d+1}]=\sum_{\alpha_2<\alpha<\alpha_1}
[\Delta_{\eta}^{\alpha}]$ being $\eta=\tilde\eta_1=\tilde\eta_2$.

Since $\eta_q$ has integer coefficients,
$\Delta_{\eta_q}^{c,0}$
is a convex (simplicial) cone with the origin as
vertex.
But in general the vertex of the convex rational cone
$\Delta_{\eta_q}^{c,\alpha}$ is a rational $(d+1)$-uple
which may have a non-integer coordinate and then
its generating
function cannot
be computed from the one
of $\Delta_{\eta_q}^{c,0}$ with an integer translation.
To avoid this problem we do the following.

Let $H_q^{\alpha}$ be the hyperplane of equation
$\eta_q+\alpha=0$,
$\alpha\in\bn$.
The difference of two elements in a given
$H_q^{\alpha}$ belongs to $l^q\cap\bz^{d+1}$. Recall that
\begin{equation*}
l^q\cap\bz^{d+1}=G^q+\bz\bdw_1^q+\dots+\bz\bdw_d^q\quad \text{and}\quad
l^q\cap\bp^{d+1}=G^q+\bn\bdw_1^q+\dots+\bn\bdw_d^q,
\end{equation*}
where $\Delta_{\eta_q}^{0}=l^q\cap\bp^{d+1}$ and $G^q$
is the fundamental set of $\Delta_{\eta_q}^0:$
$$
G^q:=\bn^d\cap \left\{\sum_{l=1}^d \mu_l\bdw_l^q:\, 0<\mu_l\leq 1\,
\,\,\text{for} \,l=1,\ldots,d\right\}.
$$

We need the following lemma which we will prove later.

\begin{lema}\label{biyec}
For every $\alpha\in \bn,$ there exists a subset $G^{q,\alpha}
\subset H_q^{\alpha}$ with a bijection between $G^q$ and $G^{q,\alpha}$
such that $H_q^{\alpha}\cap\bp^{d+1}
=G^{q,\alpha}+\bn\bdw_1^q+\dots+\bn\bdw_d^q$.

In particular for any $\bdk\in \Omega^{q,\alpha}$ there exists a unique
$\bdg\in G^{q,\alpha}$ and $(l_1,\ldots,l_d)\in \bn^d$ such that
$\bdk=\bdg+l_1\bdw_1^q+\dots+l_d\bdw_d^q.$
\end{lema}

As a consequence of equation (\ref{simplicial}),
for every $\alpha\in \bn,$ we have
\begin{equation}\label{simpli2}
\Phi_{\Delta_{\eta_q}^{\alpha}}(\bdy)=\frac{
\left(\sum _{\bdg\in G^{q,\alpha}} \bdy ^{\bdg}\right)}
{\prod _{l=1}^{d}\left(1-\bdy^{\bdw_l^q}\right)}.
\end{equation}

Assume $\Delta_\tau^t$ is of type $N$ and is
the region delimited by the inequalities
$
\eta_q+\alpha_{q,u(q)}>0,
\eta_{q+1}+\alpha_{q+1,1}<0$.
Assume also that $\Delta_\tau^t$ is different from
$\Delta_m^t$ and $\Delta_0^t$ (these will be studied
later on).
The associated vertex $\bar \tau$ of $\Gamma(h)$
gives a convex cone $\Delta_{\bar \tau}$ delimited by the inequalities
$\eta_q>0,
\eta_{q+1}<0.$
Thus the monomial ${\bar \bdx}^{\balpha} \bar{z}^n$
appears in the sum defining $h.$ Moreover
if $\bdk\in\Delta_{\bar \tau}$ and $\bdk(\barphi)=\bdk$
then
$\ord(h\circ\barphi)=n_\tau(\bdk)+N(\bdk)$.
As in lemma
\ref{2002},
$\bl^{-m(\bdk)}\mu_X(V_{n_\tau(\bdk)+N(\bdk),m(\bdk)}^{\bdk})=
\bl^{-(d+1)}(\bl-1)^{d+1}\bl^{-\sigma_{\bar \omega}(\bdk)}.$

In this case we can be more explicit
in equality (\ref{indicator}) because
$[\Delta_{\eta_{q+1}}^{c,0}]-[\Delta_{\eta_{q}}^{c,0}]=[\Delta_{{\bar \tau}_q}]+[\Delta_{{\bar \gamma}_q}],$
see equation (\ref{desmulti}). Moreover both generating functions are evaluated
at the same value $\bdy:=(\bl^{-\nu_1}T^{\alpha_1},\ldots,\bl^{-\nu_d}T^{\alpha_d},
\bl^{-1}T^{n})$.
Thus
$$
Z_{DL}^{\Delta_\tau^t,\ww}(T)= T^\delta (\bl-1)^{d+1}
\bl^{-(d+1)} \left(\Phi_{\Delta_{{\bar \tau}_q}}(\bdy)+
\sum_{0<\alpha\leq \alpha_2}
\Phi_{\Delta_{\eta_{q+1}}^{\alpha}}(\bdy)-
\sum_{0\leq \alpha< \alpha_1}
\Phi_{\Delta_{\eta_q}^{\alpha}}(\bdy)\right).
$$

\begin{lema} \label{2002b}
For each $\Delta_\tau^t$ of type $N$
its contribution to the topological zeta function is
$\chi_{\text{\rm top}}(Z_{DL}^{\Delta_\tau^t,\ww}(h^t,\omega,\bl^{-s}))=
\chi_{\text{\rm top}}(Z_{DL}^{\bar \tau}(h,\bar \omega,\bl^{-s})).
$
\end{lema}

\begin{proof} The proof follows directly first from
the description in equation (\ref{simpli2}) and the fact that
in the denominator of any
$\Phi_{\Delta_{\eta}^{\alpha}}(\bdy)$ there are at most
$d$ factors.
\end{proof}

\begin{obs} This fact will be used for both theorems \ref{thm-polos} and \ref{thm-top}.
\end{obs}

The indicator function $[\Delta_0^t]$
for the first convex rational polyhedron
$\Delta_0^t=\{\eta_1+\alpha_{1,1}<0\}$
also verifies the identity $[\Delta_0^t]=[\Delta_{\eta_1}^{c,0}]-
\sum_{0<\alpha\leq \alpha_{1,1}} [\Delta_{\eta_1}^{\alpha}].$
In particular one has for this convex polyhedron
a lemma similar to Lemma \ref{2002b}
too.

For the last convex rational polyhedron
$\Delta_m^t=\{\eta_r+\alpha_{r,u(r)}>0\},$
 $[\Delta_m^t]=[\Delta_{\eta_r}^{c,0}]+
\sum_{0\leq\alpha< \alpha_{r,u(r)}} [\Delta_{\eta_r}^{\alpha}].$
Thus
\begin{equation} \label{lastvertex}
Z_{DL}^{\Delta_m^t,\ww}(h^t,\omega,T)= T^\delta (\bl-1)^{d+1}
\bl^{-(d+1)} \left(\Phi_{\Delta_{{\bar \tau}_r}}(\bdy)+
\sum_{0\leq\alpha< \alpha_{r,u(r)}}
\Phi_{\Delta_{\eta_{r}}^{\alpha}}(\bdy)\right).
\end{equation}
And again one gets for this polyhedron
a lemma similar to Lemma \ref{2002b}.

If the convex rational polyhedron
$\Delta_\tau^t$  is of type $M$
limited by two parallel hyperplanes,
say for instance  $l^t_{q,j}$ and
$l^t_{q,j+1}$, then $\Delta_\tau^t$
also comes from a vertex, say
$t^\delta \bdx^{\balpha}z^n$.
Consider ${\bar \gamma}_q$
the corresponding compact face of $\Gamma(h)$
such that it defines $l^q$ in the dual.
One of the vertices
of ${\bar \gamma}_q$ has to be ${\bar \bdx}^{\balpha}{\bar z}^n$
because it is connected with one of the vertices of type
$N.$ We have already mention that
$[\Delta_\tau^t\cap \bp^{d+1}]$
is nothing but $\sum_{\alpha_{q,j+1}<\alpha<\alpha_{q,j}}
[\Delta_{\eta_q}^{\alpha}].$ Then
 \begin{equation} \label{formulat}
Z_{DL}^{\Delta_\tau^t,\ww}(h^t,\omega,T)= T^\delta (\bl-1)^{d+1}
\bl^{-(d+1)} \left( \sum_{\alpha_{q,j+1}<\alpha<\alpha_{q,j}}
\Phi_{\Delta_{\eta_{q}}^{\alpha}}(\bdy)\right),
\end{equation}
where $\bdy:=(\bl^{-\nu_1}T^{\alpha_1},\ldots,\bl^{-\nu_d}T^{\alpha_d},
\bl^{-1}T^{n}).$ We observe that $(\alpha_1,\ldots,\alpha_d,n)$ also belongs
to the closure of ${\bar \gamma}_q.$
Therefore after identity (\ref{simpli2})
we have proved the next proposition.

\begin{prop} \label{parallel}
If
$\Delta_\tau^t$ is of type $M$, then
$
\chi_{\text{\rm top}}(Z_{DL}^{\Delta_\tau^t,\ww}(h^t,\omega,\bl^{-s}))=0.
$
\end{prop}

We summarize the previous results and we collect also the information about the denominators of the partial \emph{motivic} zeta functions, using the above
arguments and the
description of the corresponding cones given in
lemmas \ref{1connuevo}, \ref{01con}, \ref{01con2}
and \ref{01con3}.

\begin{prop} \label{tvertex} For each vertex
${\bar \tau}\in \Gamma(h)$ and each of the vertices $\Delta_\tau^t$
of type $M$ or $N$ of $\Gamma^t(h^t)$,
the functions
$Z_{DL}^{\Delta_\tau^t,\ww}(h^t,\omega,T)$ and $
Z_{DL}^{{\bar \tau}}(h,\bar \omega,T)$
belong to the subring
$\bz[\bl,\bl^{-1},(1-\bl^{-1}T)^{-1},(1-\bl^{-\sigma_\omega(\bdw_j^q)}
T^{m_h(\bdw_j^q)})^{-1},(1-\bl^{-\nu_i}T^{N_i})^{-1}][T],\,
j,i\in\{1,\ldots,d\},$ and $q=1,\ldots,r.$

There is a
one-to-one correspondence between elements of type $N$ and
vertices of $\Gamma(h).$  If $\Delta_\tau^t$ is of type $N$
and under this bijection
corresponds with ${\bar \tau}\in \Gamma(h)$ then
$
\chi_{\text{\rm top}}(Z_{DL}^{\Delta_\tau^t,\ww}(h^t,\omega,\bl^{-s}))=
\chi_{\text{\rm top}}(Z_{DL}^{{\bar \tau}}(h,\bar \omega,\bl^{-s})).
$

Moreover, if
$\Delta_\tau^t$ is of type $M$ then
$
\chi_{\text{\rm top}}(Z_{DL}^{\Delta_\tau^t,\ww}(h^t,\omega,\bl^{-s}))=0.
$
\end{prop}

\begin{proof}[Proof of Lemma \ref{biyec}]
Let $\bdh\in H_q^{\alpha}\cap\bz^{d+1}$, and let $\bdg \in G^q$.
Write
$$
\bdh+\bdg=(h_1,h_2,...,h_d,h_{d+1}).
$$
Using the division algorithm,
for $l=1,...,d+1$, we decompose
$h_l=s_l p_l^q+r_l$, where $s_l \in \bz$ and
$0 < r_l\leq p_l^q.$
Since $\bdh\in H_q^{\alpha}$ and $\bdg \in l^q$, then $\bdh+\bdg \in H_q^{\alpha}$.
This implies that
$$
\sum_{l=1}^d b_l^q r_l+\alpha=
h_{d+1}n_1^q- \sum_{l=1}^d b_l^q s_l p_l^q
=(h_{d+1}- \sum_{l=1}^d \bar{b}_l^q s_l) n_l^q,\text{ i.e. }
h_{d+1}- \sum_{l=1}^d \bar{b}_l^q s_l \in \bp.
$$
Let us denote
$$\bdh^1_\bdg=(r_1,...,r_d,h_{d+1}- \sum_{l=1}^d \bar{b}_l^q s_l)
\in H_q^{\alpha} \cap \bp^{d+1}.
$$
We have $\bdh+\bdg=\bdh^1_\bdg+s_1\bdw_1^q+\dots+s_d\bdw_d^q$.

\smallbreak
Let $\bdh' \in H_q^{\alpha} \cap \bp^{d+1}$.
Then $\bdh'-\bdh \in  l^q \cap \bz^{d+1}$.
There exist $\bdg \in G^q$ and $(u_1,...,u_d)\in\bz^d$ such that
$$
\bdh'=\bdh+\bdg+u_1\bdw_1^q+\dots+u_d\bdw_d^q,\text{ i.e. }
\bdh'=\bdh^1_\bdg+\sum_{l=1}^d (u_l-s_l)\bdw_l^q.
$$
If $\bdh'=(h'_1,\dots,h'_d,h'_{d+1})$ then, for $l=1,...,d,$
we have $h'_l=r_l+ (u_l-s_l)p_l^q$.

Since $h'_l\in\bp$ and
$0< r_l \leq p_l^q$, then $(u_l-s_l) \in \bn$.
We define
\begin{equation}\label{gqalpha}
G^{q,\alpha}:=\{\bdh^1_\bdg\ |\ \bdg\in G^q \}.
\end{equation}
We have proved that
$H_q^{\alpha}\cap\bn^{d+1}_{>0}
=G^{q,\alpha}+\bn\bdw_1^q+\dots+\bn\bdw_d^q$.

It is enough to prove that $G^q$ and $G^{q,\alpha}$ are bijective.
Now let us suppose that $\bdh^1_\bdg=\bdh^1_{\bdg'}$
for some $\bdg,\bdg'\in G^q$.
We can write:
$$
\bdh+\bdg=\bdh^1_\bdg+ \sum_{l=1}^d s_l \bdw_l^q
\text{ and }
\bdh+\bdg'=\bdh^1_{\bdg'}+ \sum_{l=1}^d s'_l \bdw_l^q.
$$
We have $\bdg-\bdg'= \sum_{l=1}^d (s_l-s'_l)\bdw_l^q\in l^q\cap\bz^{d+1}$.
Since the elements of this set can be written in a unique form
as elements of $G^q+\bz\bdw_1^q+\dots+\bz\bdw_d^q$
then $\bdg=\bdg'$.
 \end{proof}

\subsection{Edges of the Newton polytope }
\mbox{}
\medskip

We now deal with edges.
Consider the edges
$\bar{\gamma}_1,\dots,\bar{\gamma}_r$ of $\Gamma(h).$ Fix
$q=1,\dots,r$.
Since $f_{\blambda_q}=\prod_{{\bar j}=1}^{v(q)}(z^{n_1^q}-
\bar{\beta}_{\bar j}^q
\bdx^{n_1^q\blambda_q})^{m_{q,j}}$ is the quasihomogeneous part of $f$
corresponding to this face, 
$[N_{{\bar \gamma}_q}]=[\bge\cap\{f_{\blambda_q}=0\}]
=v(q)(\bl-1)^d$, see Remark \ref{classqo}.

In Step \ref{K-newton}, for each $j\in\{1,\dots,u(q)\}$, the polynomial
${{\widehat f}^t}_{q,j}:=\prod_{w=1}^{w(q,j)}
((z^{n_1^q})^{s_1^{q,j}}-\beta_w t^{\alpha_{q,j}}
(\bdx^{n_1^q\blambda_{\kappa,q}})^{s_1^{q,j}})^{m_{q,j,w}}$ was defined.
This polynomial is associated with a region
$M_t^{q,j}$ which is contained in a hyperplane parallel to $l^q$. Next result is trivial
because under the hypothesis of the proposition there is no
points with positive integer coordinates in $M_t^{q,j}.$ 

\begin{prop} If $s_1^{q,j}>1$ then $M_t^{q,j}\cap\bn^{d+1}=\emptyset$ and
$Z_{DL}^{M_t^{q,j},\ww}(h^t,\omega,T)=0$.
\end{prop}

From now on we will suppose $s_1^{q,j}=1$.
Because of this proposition we cannot state the equality for the topological zeta function in Theorem \ref{thm-polos}. For the other edges, following the notation
in \S\ref{secnondeg} we break the sets of arcs in $A$- and $B$-parts.
We will consider,
$$
{\widehat f}^t_{q,j|t=1}:=\prod_{w=1}^{w(q,j)}
(z^{n_1^q}-\beta_w
\bdx^{n_1^q\blambda_{\kappa,q}})^{m_{q,j,w}}.
$$

\begin{obs} Observe that for $\tww$-quasi-ordinary and Newton compatible series, for each $q$ we have $u(q)=1$.
The extra hypothesis $s_1^{q,j}=1$ is fulfilled
by $\tww$-quasi-ordinary and Newton compatible series
and this is why equality for topological zeta functions may be obtained in Theorem \ref{thm-top}.
\end{obs}

\begin{paso} Edges in the A-part.
\end{paso}

Take a convex rational cone $\Delta_{{q,j}}^t$ corresponding to an edge.
We decompose its contribution into two disjoint parts as in the
non-degenerated case.
Given $\bdk\in\Delta_{{q,j}}^t,$ we set:
\begin{multline*}
V_{n,A,m}^{\bdk,t}:=\{\barphi\in \cl_0(X)\setminus\cl_0(S):\\
\bdk(\barphi)=\bdk,\ord(h^t\circ \barphi)=n+\theta,
\ord(\omega\circ\barphi)=m,{\widehat f}^t_{{q,j}|t=1}(\bda(\barphi))\neq 0\}.
\end{multline*}
These semialgebraic subsets are empty unless $m=m(\bdk)$
and $n=\eta_q(\bdk)+\alpha_{q,j}+N(\bdk)$.
In particular  $[\bge\cap\{{\widehat f}^t_{q,j|t=1}=0\}]=w(q,j)(\bl-1)^d
\in K_0(\text{Var}_k)$, see Remark \ref{classqo}.

Let us (re)define
$$
Z_{q,j,A}^{\ww}(h^t,\omega,T):=\sum_{\bdk\in\Delta_{{q,j}}^t}
\bl^{-m(\bdk)} \mu_X(V_{\hat n_q(\bdk)+\alpha_{q,j}+N(\bdk),A,m(\bdk)}^{\bdk,t})
T^{\eta_q(\bdk)+\alpha_{q,j}+N(\bdk)}.
$$
In the same way, given $\bdk\in\Delta_{{\bar \gamma}_q}$ we set:
\begin{multline*}
V_{n,A,m}^{\bdk}:=\{\barphi\in \cl_0(X)\setminus\cl_0(S):\\
\bdk(\barphi)=\bdk,\ord(h\circ \barphi)=n,
\ord(\omega\circ\barphi)=m,f_{\blambda_{\kappa,q}}(\bda(\barphi))\neq 0\}.
\end{multline*}
Again the semialgebraic subsets are empty unless $m=m(\bdk)$
and $n=\eta_q(\bdk)+N(\bdk).$
We recall that
$$
Z_{{{\bar \gamma}_q},A}(h,\bar \omega,T)=\sum_{\bdk\in\Delta_{{\bar \gamma}_q}}
\bl^{-m(\bdk)} \mu_X(V_{\hat n_q(\bdk)+N(\bdk),A,m(\bdk)}^{\bdk})
T^{\eta_q(\bdk)+N(\bdk)}.
$$

\begin{lema} \label{tedgesA}
For each $q=1,\ldots,r$, $j=1,\ldots,u(q),$
the functions
$Z_{q,j,A}^{\ww}(h^t,\omega,T)$ and $Z_{{\bar \gamma}_q,A}(h,\bar \omega,T)$
belong to the subring
$\bz[\bl,\bl^{-1},(1-\bl^{-\sigma_\omega(\bdw_l^q)}
T^{m_h(\bdw_l^q)})^{-1}][T],\,
l=1,\ldots,d.$ Moreover, the following equality
holds:
$$
\sum_{j=1}^{u(q)} \chi_{\text{\rm top}}(Z_{q,j,A}^{\ww}(h^t,\omega,\bl^{-s}))=
\chi_{\text{\rm top}}(Z_{{\bar \gamma}_q,A}(h,\bar \omega,\bl^{-s})).
$$
\end{lema}

\begin{proof}
We can follow the proofs of Theorem \ref{nondeg}, Lemma \ref{2002}
and  Lemma \ref{2002b}. Since we are in the $A$-part of the
decomposition, the measure of the sets
$V_{\eta_q(\bdk)+N(\bdk),A,m(\bdk)}^{\bdk}$
and $V_{\eta_q(\bdk)+\alpha_{q,j}+N(\bdk),A,m(\bdk)}^{\bdk,t}$
can be explicitly computed.
Since $h^t$ and $h$ has the same characteristic
exponents
then the Newton polyhedra  $\Gamma^K(h^t)$ and $\Gamma(h)$ coincide and
we get a sum
as (\ref{formulat})
but with only one $\alpha=\alpha_{q,j}.$
Then
\begin{equation}\label{contface}
Z_{q,j,A}^\ww(h^t,\omega,T)=
\bl^{-(d+1)}\left((\bl-1)^{d+1}-w(q,j)(\bl-1)^d\right)
T^{\alpha_{q,j}}\Phi_{\Delta_{\eta_{q}}^{\alpha_{q,j}}}(\bdy),
\end{equation}
where $\bdy:=(\bl^{-\nu_1}T^{\alpha_1},\ldots,\bl^{-\nu_d}T^{\alpha_d},
\bl^{-1}T^{n})$ and  $(\alpha_1,\ldots,\alpha_d,n)$ is
an element in the closure of ${\bar \gamma}.$ In the same way
\begin{equation*}
Z_{{{\bar \gamma}_q},A}(h,\bar \omega,T)=
\bl^{-(d+1)}\left((\bl-1)^{d+1}-v(q)(\bl-1)^d\right)
\Phi_{\Delta_{\eta_{q}}^{\alpha_{q,j}}}(\bdy).
\end{equation*}
Since there is a bijection between  the sets $G^{q}$ and $G^{q,\alpha_{q,j}}$
and $v(q)=\sum w(q,j)$ for $j=1,\ldots,u(q)$ then
taking the usual Euler characteristic the lemma is proved.
\end{proof}

To sum up we collect the $A$-part in the arc decomposition as follows:
\begin{equation}\label{primera}
\begin{split}
Z_{DL}^A(h,\bar \omega,T)=\sum_{\bar \tau\,\text{ vertex}}
Z_{DL}^{\bar \tau}(h,\bar w,T)+
\sum_{q=1}^r Z_{\bar \gamma_q,A}(h,\bar \omega,T), \\
Z_A^\ww(h^t,\omega,T):=\sum_{\Delta_\tau^t \,\text{ vertex}}
Z_{DL}^{\Delta_\tau^t,\ww}
(h^t,\omega,T)
+\sum_{q=1}^r\sum_{j=1}^{u(q)} Z_{q,j,A}^\ww(h^t,\omega,T).
\end{split}
\end{equation}

\begin{obs}
Up to now we have proved two facts:
\begin{itemize}
\item For a $\ww$-quasi-ordinary power series
 the denominators of the $A$-part of the motivic zeta function are controlled.

\item For a $\tww$-quasi-ordinary power series
 which is Newton compatible, identity (\ref{eqAtop}) also holds.
\end{itemize}
We need to replace $k$ by $k[[t]]$ because of computations of the $B$-part.
Observe that in the $B$-part there are no vertices.
\end{obs}

\begin{paso} Edges in the $B$-part.
\end{paso}

Given $\bdk\in\Delta_{q,j}^t,$ consider the semialgebraic subsets
\begin{multline*}
V_{n,B,m}^{\bdk,t}:=\{\barphi\in \cl_0(X)\setminus\cl_0(S):\\
\bdk(\barphi)=\bdk,\ord(h^t\circ \barphi)=n+\theta,
\ord(\omega\circ\barphi)=m,
{\widehat f}^t_{{q,j}|t=1}(\bda(\barphi))= 0\}.
\end{multline*}
If $\bdk\in\Delta_{{\bar \gamma}_q}$ we define $V_{n,B,m}^{\bdk}$
in the same way.
We have to compute $Z_{{q,j},B}^{\ww}(h^t,\omega,T)$ and
$Z_{{{\bar \gamma}_q},B}(h,\bar \omega,T)$ which have
obvious definitions.

The sets $\Delta_{{\bar \gamma}_q}\cap \bp^{d+1}$, resp.
$\Delta^t_{q,j}\cap \bp^{d+1}$, are the disjoint unions
of the sets $\bdg+\bdw_1^q\bn+\ldots+\bdw_d^q\bn$ where
$\bdg\in G^q$, resp. $\bdg\in G^{q,\alpha_{q,j}}$.
Accordingly, we define the sets
$$
V_{n,B,m}^{{\bar \gamma}_q,\bdg}:=\bigcup_{\bdk\in\bdg+\bdw_1^q\bn+\ldots+
\bdw_d^q\bn} V_{n,B,m}^{\bdk},\qquad \bdg\in G^{q},
$$
and
$$
V_{n,B,m}^{q,j,\bdg,t}:=\bigcup_{\bdk\in\bdg+\bdw_1^q\bn+\ldots+\bdw_d^q
\bn} V_{n,B,m}^{\bdk,t}, \qquad \bdg\in G^{q,\alpha_{q,j}}.
$$

Let $\barphi\in V_{n,B,m}^{{\bar \gamma}_q,\bdg},$ with $\bdg\in G^q$;
this arc is related to exactly one
Newton component of $f$, see subsection \ref{newtonmaps}.
There is a unique $\bar{j}=1,\ldots,v(q)$
such that $\bda(\barphi)\in \bge$
belongs to the zero locus $V_{\bar{j}}^q\subset \ba_k^{d+1}$ of
the quasi-ordinary polynomial
$f_{q,j}=z^{n_1^q}-\bar{\beta}^q_{\bar j}\bdx^{n_1^q\blambda_q}$.

\begin{obs}\label{notodos} Let $J_q$ be the subset of $\{1,\dots,d\}$ of the non-zero coordinates of $\blambda_q$. It may happen that $J_q\subsetneqq \{1,\dots,d\}$. For the sake of simplicity we suppose that there is equality and we will point out where the non-equality may affect.
\end{obs}

We decompose the set $V_{n,B,m}^{{\bar \gamma}_q,\bdg}$
in $v(q)$ disjoint sets according to this property:
$$
V_{n,B,m}^{{\bar \gamma}_q,\bdg}=\bigcup_{\bar{j}=1}^{v(q)} V_{n,B,m}^{{\bar \gamma}_q,\bdg,
\bar{j}}.
$$

In the same way the sets $V_{n,B,m}^{t,q,j,\bdg,w}$,
$w=1,\dots,w(q,j)$, are defined. It is also possible to decompose
$$
V_{n,B,m}^{{\bar \gamma}_q,\bdg,\bar{j}}=
\bigcup_{\bdr:=(r_1,\ldots,r_d)\in\bn^d}
V_{n,B,m}^{{\bar \gamma}_q,\bdg,j,\bdr},\quad
$$
where
$
V_{n,B,m}^{{\bar \gamma}_q,\bdg,\bar{j},\bdr}:=
\{\barphi \in V_{n,B,m}^{{\bar \gamma}_q,\bdg,\bar{j}}:
\bdk(\barphi)=\bdg+r_1\bdw_1^q+\ldots+r_d\bdw_d^q\,\},
$
and consider in the same way $V_{n,B,m}^{t,q,j,\bdg,w,\bdr}$.
Define 
$$
Z^{{\bar \gamma}_q,\bdg,\bar{j},\bdr}(T):=
\sum_{n\geq 1}\left(\sum_{m\geq 1} \bl^{-m}
\mu_X(V_{n,B,m}^{{\bar \gamma}_q,\bdg,\bar{j},\bdr})\right)
T^n,
$$
and also $Z^{\ww,\bdg,w,\bdr}_{q,j}(T)$.
We define
$
Z^{{\bar \gamma}_q,\bdg,j}(T):=\sum_{\bdr\in \bn^d} Z^{{\bar \gamma}_q,\bdg,\bar{j},\bdr}(T)
$
and the corresponding $Z^{\ww,\bdg,w}_{q,j}(T)$;
then we have the decompositions:
\begin{equation} \label{induction}
\begin{split}
Z_{DL}^\ww(h^t,\omega,T)=Z_A^\ww(h^t,\omega,T)+\sum_{q=1}^r
\sum_{j=1}^{u(q)}\sum_{\bdg \in G^{q,\alpha_{q,j}}}
\sum_{w=1}^{w(q,j)} Z_{q,j}^{\ww,\bdg,w}(T)\\
Z_{DL}(h,{\bar \omega},T)=Z_{DL}^A(h,{\bar \omega},T)+\sum_{q=1}^r
\sum_{\bdg \in G^{q}}\sum_{{\bar j}=1}^{v(q)}
Z^{{\bar \gamma}_q,\bdg,{\bar j}} (T),
\end{split}
\end{equation}
where $v(q)=\sum_{j=1}^{u(q)} w(q,j)$.
Using Proposition \ref{tvertex}
and Lemma \ref{tedgesA} we have proved that $Z_A^\ww(h^t,\omega,T)$
verifies  Theorems \ref{thm-polos} and \ref{thm-top}
in their respective cases.
Then we have to prove them for arcs in the $B$-part
of the decomposition. In particular we will describe in the following
steps the computations in the identities of (\ref{induction}).

\begin{paso}\label{unos} Computations for $h$ and
$\bdg=\bdg^1:=\bdw_1^q+\ldots+\bdw_d^q\,\in G^q$, and ${\bar j}=1,\ldots,v(q)$.
\end{paso}

We use the notation of subsection \ref{newtonmaps}.
To compute $Z^{{\bar \gamma}_q,\bdg
^1,\bar{j},\bdr}(T)$, for $\bdr\in\bn^d$,
then
$\bdk=(r_1+1)\bdw_1^q+\ldots+(r_d+1)\bdw_d^q$
has coordinates
$
k_l=(r_l+1)p_l^q$ with $l=1,\ldots,d,$ and
$k_{d+1}=\sum_{l=1}^d (r_l+1){\bar b}_l^q.$
If $\barphi$ verifies that $\bdk(\barphi)=\bdk$ then
$m_\bdk:=
\ord (\bar{\omega}\circ\barphi)=\sum_{l=1}^d (r_l+1)
p_l^q (\nu_l-1).$
Therefore $m_\bdk$ is the unique value of $m$ for
which $V_{n,B,m}^{{\bar \gamma}_q,\bdg
^1,\bar{j},\bdr}$ may be
non-empty. Hence
$$
Z^{{\bar \gamma}_q,\bdg^1,\bar{j},\bdr}(T)=
\sum_{n\geq 1}
\bl^{-m_\bdk}
\mu_X(V_{n,B,m_\bdk}^{{\bar \gamma}_q,\bdg^1,\bar{j},\bdr}) T^n.
$$
Fix the
parametrization $\pi_{\bar{j}}^q:\ba_k^d\to V_{\bar j}^q$
as in subsection \ref{newtonmaps}.
Take $\barphi\in V_{n,B,m_\bdk}^{{\bar \gamma}_q,\bdg^1,\bar{j},\bdr}$
and the unique $\bds^0\in \bgd$ such that $\pi_{\bar j}^q(\bds^0)=\bda(\barphi)$.
The idea is to lift
$\barphi$
to an affine space $Y:=\ba_k^{d+1}$ with coordinates $(\bdy,z_1)$
using the \emph{Newton
map} $\pi_{q,{\bar j}}$ associated with
$f_{\bar j}^{q}$.

\smallbreak

This Newton map defines a
$k[t]$-morphism $\pi_{q,\bar j}:\mathcal{L}(Y)\to
\mathcal{L}(X)$
(in fact nothing is done on the variable $t$).
We can apply to $\pi_{q,{\bar j}}$ the change of variables formula, see theorem
\ref{changevariables}.

For $l=1,\dots,d,$ the $l$-th component of $\barphi$ is  $\varphi_l(t)=t^{(r_l+1)p_l^q}v_l(t),$
where $v_l(t)\in k[[t]]$ such that $v_l(0)=(s_l^0)^{p_l^q}$.
There exists a unique  $w_l(t)\in k[[t]]$ such that
$w_l(0)=s_l^0$ and $w_l(t)^{p_l^q}=v_l(t)$.
If $\psi_l(t):=t^{(r_l+1)}w_l(t)$ then
$$
\varphi_{d+1}(t)=t^{(r_1+1)
{\bar b}_1^q+\ldots+(r_d+1){\bar b}_d^q} v_{d+1}(t)
\text{ and }
v_{d+1}(0)=\prod_{l=1}^d (s_l^0)^{{\bar b}_l^q}.
$$
The \emph{equation}
\begin{equation*}
\begin{split}
t^{(r_1+1){\bar b}_1^q+\ldots+(r_d+1){\bar b}_d^q}
v_{d+1}(t)=&(\psi_{d+1}(t)+{\bar \beta}_{\bar j}^q)\psi_1(t)^{{\bar b}_1^q}
\dots \psi_d(t)^{{\bar b}_d^q}=\\
=&(\psi_{d+1}(t)+{\bar \beta}_{\bar j}^q)
t^{(r_1+1){\bar b}_1^q+\ldots+(r_d+1){\bar b}_d^q}
w_1(t)^{{\bar b}_1^q}
\dots w_d(t)^{{\bar b}_d^q}
\end{split}
\end{equation*}
completely determines an arc $\psi_{d+1}(t)$ such that
$\psi_{d+1}(0)=0$. So there is a unique
$\bpsi:=(\psi_1,\dots,\psi_d,\psi_{d+1})\in\cl_0(Y)$
such that $\pi_{q,{\bar j}}(\bpsi)=\barphi$.

Moreover, if ${\bar h}_{q,{\bar j}}:=h\circ \pi_{q,{\bar j}}$ and
$\omega_{q,{\bar j}}:={\bar \omega}\circ \pi_{q,{
\bar j}}$, we denote by
$B_{n,m_\bdk}^{\bdr}$ the set of arcs $\bpsi\in\cl_{0}(Y)$
such that

\begin{itemize}

\item $\ord({\bar h}_{q,{\bar j}} \circ\bpsi)=n$,

\item $\ord(\omega_{q,{\bar j}}\circ\bpsi)=m_\bdk$,

\item if $\bdk(\bpsi)=(\tilde k_1,\dots,\tilde k_d,\tilde k_{d+1})$ then
$\tilde k_l=r_l+1$, for $l=1,\dots,d$.
\end{itemize}
Then $\pi_{q,{\bar j}}$ defines a bijection between $B_{n,m_\bdk}^{\bdr}$
and $V_{n,B,m_\bdk}^{{\bar \gamma}_q,\bdg^1,{\bar j},\bdr}$.
The order of the Jacobian of the $k[t]$-morphism
$\pi_{q,{\bar j}}$ is constant on $B_{n,m_\bdk}^{\bdr}$ and equals to:
$$
\ord_t(\mathcal{J}_{\pi_{q,{\bar j}}})=\sum_{l=1}^d (r_l+1) (p_l^q+{\bar b}_l^q-1).
$$
The set $B_{n,m_\bdk}^{\bdr}$ is strongly measurable
in $\cl(Y)$. Hence $V_{n,B,m_\bdk}^{{\bar \gamma}_q,\bdg^1,{\bar j},
\bdr}$,
its bijective image by $\pi_{q,{\bar j}}$,
is also strongly measurable in $\cl(X)$,
\cite[Theorem A.8]{dl:02}.
 To measure
$V_{n,B,m_\bdk}^{{\bar \gamma}_q,\bdg^1,{\bar j},\bdr}$
the change variables formula gives
$$
\mu_X(V_{n,B,m_\bdk}^{{\bar \gamma}_q,\bdg^1,{\bar j},\bdr})=
\int_{V_{n,B,m_\bdk}^{{\bar \gamma}_q,\bdg^1,{\bar j},\bdr}}d\mu_X
=\int_{B_{n,m_\bdk}^{\bdr}}
\bl^{-\ord_t \mathcal{J}_{\pi_{q,{\bar j}}}(y)} d\mu_Y.
$$
If ${\widehat \omega}_1$ is the pullback by
$\pi_{q,{\bar j}}$ of the differential form
$d{\bar x}_1\wedge \dots\wedge d{\bar x}_d\wedge d{\bar z}$
then
$$
\mu_X(V_{n,B,m_\bdk}^{{\bar \gamma}_q,\bdg^1,{\bar j},\bdr})\!=\!\sum_{k\in\bp} \bl^{-k}
\mu_Y(B_{n,m_\bdk}^{\bdr}\cap\,\{\ord_t({\widehat \omega}_1)=k\})\!=\!
\bl^{-\sum_{l=1}^d (r_l+1) (p_l^q+{\bar b}_l^q-1)}
\mu_Y(B_{n,m_\bdk}^{\bdr}).
$$
Putting all these terms together
$$
Z^{{\bar \gamma}_q,\bdg^1,{\bar j},\bdr}(T)=
\sum_{n\geq 1} \bl^{-\sum_{l=1}^d (r_l+1) (p_l^q \nu_l+{\bar b}_l^q-1)}
\mu_Y(B_{n,m_\bdk}^{\bdr}) T^n.
$$

Since $\ord_t(\omega_{q,{\bar j}}\circ\bpsi)=\sum_{l=1}^d (r_l+1)
(p_l^q\nu_l+{\bar b}_l^q-1)$ for any arc $\bpsi\in B_{n,m}^r$,
$$
Z^{{\bar \gamma}_q,\bdg^1,{\bar j},\bdr}(T)=\sum_{n\geq 1}
\bl^{-m_\bdk}\mu_Y(B_{n,m_\bdk}^{\bdr})
T^n.
$$
The set $B_{n,m}^{\bdr}$ is empty if $m\neq m_\bdk$ so
\begin{equation}\label{caso11}
Z^{{\bar \gamma}_q,\bdg^1,{\bar j}}(T)=\sum_{\bdr\in\bn^d}\sum_{n\geq 1}
\bl^{-m_\bdk}\mu_Y(B_{n,m_\bdk}^{\bdr})T^n =Z_{DL}({\bar h}_{q,{\bar j}},
\omega_{q,{\bar j}},T),\, \forall {\bar j}=1,\ldots,v(q).
\end{equation}

\begin{obs} If 
$h$ has $\lgt(h)=1$,
no more steps to get a formula for
$Z_{DL}(h,\omega,T)$ are needed,
see section \ref{curvascase}. 
This happen because each fundamental set $G^q$ has only one element.

\end{obs}

\begin{paso} \label{doses} Computations for $h^t$, $\bdg\in G^{q,\alpha_{q,j}}$
 and $w=1,\ldots,w(q,j).$
\end{paso}
By the proof of Lemma \ref{biyec}, in particular its definition in (\ref{gqalpha})
one can write $\bdg=(\mathbf 0,\dfrac{\alpha_{q,j}}
{n_1^q})+\sum_{l=1}^d \mu_l^\bdg \bdw_l^q$. Thus
$$
\bdg=(\mu_1^\bdg p_1^q,\ldots,\mu_d^\bdg p_d^q,
\frac{\alpha_q}{n_1^q}+\sum_{l=1}^d \mu_l^\bdg {\bar b}_l^q)\in \bp^{d+1}.
$$
By definition of $G^{q,\alpha_{q,j}}$,
 $0<\mu_l^\bdg\leq 1$ is rational.
Choose $v_1,\dots,v_d\in\bp$
such that $\mu_l^\bdg=\dfrac{v_l}{p_l^q}$, $l=1,\dots,d$ and let
$c^\bdg\in\bp$ the $(d+1)^\text{th}$-coordinate of $\bdg$.

Fix $\bdr\in\bn^d$. Let $\barphi\in V_{n,B,m}^{t,q,j,\bdg,w,\bdr}=\{\barphi\in
V_{n,B,m}^{t,q,j,\bdg}:\bdk(\barphi)=\bdg+r_1\bdw_1+\ldots+r_d\bdw_d \}$
be an arc and denote:
\begin{itemize}
\item $\bdk:=\bdk(\barphi)=\left((\mu_1^\bdg+r_1)p_1^q,\ldots,
(\mu_d^\bdg+r_d)p_d^q,\dfrac{\alpha_{q,j}}{n_1^q}+\displaystyle\sum_{l=1}^d
(\mu_l^\bdg+r_l){\bar b}_l^q\right)$.

\item $\bda:=\bda(\barphi)\in \bge\cap
\{{\widehat f}^{t,w}_{q,j}=0\}$,  where
${\widehat f}^{t,w}_{q,j}=z^{n_1^q}-
\beta_w \bdx^{n^q_1 \blambda_{\kappa_q}}.$
Consider $\bda=\pi_j^{q,w}(s^0)$,
$s^0\in \bgd$ with the fixed parametrization
of $\{{\widehat f}^{t,w}_{q,j}=0\}$ induced by $\pi_{j}^{q,w}:
\ba_k^d\to \bge\cap \{{\widehat f}^{t,w}_{q,j}=0\}$: $(s_1,\ldots,s_d)\mapsto
(s_1^{p_1^q},\ldots,s_d^{p_d^q}, \alpha_w s_1^{{\bar b}_1^q}\dots
s_d^{{\bar b}_d^q})$ where $(\alpha_w)^{n_1^q}=\beta_w.$
\end{itemize}

For each $\barphi$ with $\bdk(\barphi)=\bdk$ then
$
m_\bdk:=\ord ({\omega}\circ\barphi)=\sum_{l=1}^d
(r_l+\mu_l^\bdg) p_l^q (\nu_l-1).
$
Since $m_\bdk$ is the unique value of $m$ for
which $V_{n,B,m}^{t,q,j,\bdg,w,\bdr}$ may be
non-empty,
$$Z^{\ww,\bdg,w}_{q,j}(T)=\sum_{r\in \bn^d}Z^{\ww,\bdg,w,\bdr}_{q,j}(T)=
\sum_{r\in\bn^d} \bl^{-m_\bdk}
\mu_X(V_{n,B,m_\bdk}^{t,q,j,\bdg,w,\bdr})
T^n.
$$

Let $Y=\ba_k^{d+1}$ and let $\pi_\bdg:\cl(Y)\to \cl(X)$ be the $k[t]$-morphism
defined by
$x_l=t^{v_l} y_l,$ for every $l=1,\ldots,d,$
and $z=t^{c^\bdg}z_1$.

\begin{lema} Under the above conditions:
\begin{enumerate}[\rm(1)]
\item The map $\pi_\bdg$ defines a $k[t]$-morphism.
\item The Jacobian of $\pi_\bdg$ has constant order $
\nu_\bdg:=\ord_t (J_{\pi_\bdg})=\sum_{l=1}^d v_l
c^\bdg$ on the
arcs of $\cl(Y)$.

\item The $K$-quasi-ordinary function
$h^{\{t\}}:=h^t\circ \pi_\bdg$ has the same characteristic exponents (and monomials)
as $h^t$, and consequently the same as $h$ (this follows
directly from Remark \ref{Achange}).
\end{enumerate}
\end{lema}

Fix $\bdr\in\bn^d$, and take an arc
$\barphi\in V_{n,B,m_\bdk}^{t,q,j,\bdg,w,\bdr}.$
There exists a unique arc $\bpsi \in \cl(Y)$ such that
$\barphi=\pi_\bdg\circ\bpsi$. Properties
of $\bpsi$ will depend 
on $\bdr$:
\begin{enumerate}[({2}.a)]

\item \label{caso1}$\bdr\in\bp^d$; in this case $\bpsi\in\cl_0(Y)$,
$\bda(\bpsi)=\bda$ and
$$
\bdk(\bpsi)=\left(r_1 p_1^q,\ldots,
 r_d p_d^q,\displaystyle\sum_{l=1}^d r_l {\bar b}_l^q\right).
$$

\item \label{caso2}$\bdr\in\bn^d\setminus\bp^d$; in this case
$\bpsi\in\cl_{(\bdy,0)}(Y)$.
More precisely; let $J_\bdr \subsetneqq \{1,\dots,d\}$ be the set of indices
such that $r_l=0\Leftrightarrow l\in J_\bdr$;
we observe that $\emptyset\neq J_\bdr\subsetneqq \{1,\dots,d\}$.
The point $\bdy$ corresponds to $(y_1,\dots,y_d)$ where
$y_l=(s^0_l)^{p_l^q}$ if
$l\in J_\bdr$ and $y_l=0$ otherwise.
We include in this case the case where
$\bdr={\mathbf 0}\in\bn^d$ and $J=\{1,\dots,d\}$.
\end{enumerate}

\begin{obs}\label{notodos1} Note that in fact $J\subset J_q$,
see Remark \ref{notodos}.
\end{obs}
\smallbreak
Decompose
$\displaystyle Z^{\ww,\bdg,w}_{q,j}(T)=Z^{\ww,\bdg,w}_{1,q,j}(T)+
\sum_{\emptyset\subsetneqq J\subset\{1,\dots,d\}}
Z^{\ww,\bdg,w}_{J,q,j}(T)$,
where
$$
Z^{\ww,\bdg,w}_{1,q,j}(T):=
\sum_{\bdr\in\bp^d}\sum_{n\geq 1}
\bl^{-m_\bdk} \mu_X(V_{n,B,m_\bdk}^{t,q,j,\bdg,w,\bdr})
T^n
$$
and
$$
Z^{\ww,\bdg,w}_{J,q,j}(T):=\sum_{\bdr\in\bn^d: J=J_\bdr}\sum_{n\geq 1}
 \bl^{-m_\bdk}
\mu_X(V_{n,B,m_\bdk}^{t,q,j,\bdg,w,\bdr}) T^n.
$$

Define
$\omega^{\{t\}}_\bdg:=\pi_\bdg^*\omega$ and for $\bdr\in\bp^d$
denote by $W_{n,B,m_\bdk}^{t,q,j,\bdg,w,\bdr}$ the semialgebraic subset of arcs
$\bpsi\in\cl_0(Y)$ such that:
\begin{itemize}

\item $\ord(h^{\{t\}}\circ\bpsi)=n$.

\item $\bdk(\bpsi)=(r_1 p_1^q,\ldots,
 r_d p_d^q,\displaystyle\sum_{l=1}^d r_l {\bar b}_l^q)$.

\item $\ord(\omega^{\{t\}}_\bdg\circ\bpsi)=m_\bdk+\nu_\bdg$.
The differential form $\omega^{\{t\}}_\bdg$ is defined as $t^{A_\bdg} \omega_\bdg$
where
$\omega_\bdg:=\prod_{l=1}^d y_l^{\nu_l-1} d\bdy\wedge dz_1$
and $A_\bdg:=\sum_{l=1}^d v_l(\nu_l-1).$
Then ${\tilde m}_\bdk:=\ord(\omega_\bdg\circ\bpsi)=m_\bdk-A_\bdg$
and
$$\nu_\bdg+A_\bdg=
\frac{\alpha_q}{n_1^q}+\sum_{l=1}^d \mu_l^\bdg(p_l^q\nu_l+{\bar b}_l^q).$$

\end{itemize}

The map $\pi_\bdg$ induces a bijection between the sets
$V_{n,B,m_\bdk}^{t,q,j,\bdg,w,\bdr}$
and $W_{n,B,m_\bdk}^{t,q,j,\bdg,w,\bdr}$. Since
$\mu_X(V_{n,B,m_\bdk}^{t,q,j,\bdg,w,\bdr})=
\bl^{-\nu_\bdg}\mu_Y(W_{n,B,m_\bdk}^{t,q,j,\bdg,w,\bdr}),$ applying
the change of variables formula to $\pi_\bdg$
we have
\begin{equation}
\begin{split}
Z^{\ww,\bdg,w}_{1,q,j}(T)=&\sum_{\bdr\in\bp^d}\sum_{n\geq 1}
\bl^{-m_\bdk-\nu_\bdg} \mu_Y(W_{n,B,m_\bdk}^{t,q,j,\bdg,w,\bdr})
T^n \\
=& \sum_{\bdr\in\bp^d}\sum_{n\geq 1}
\bl^{-\nu_\bdg-A_\bdg-\ord_t(\omega_\bdg)}
\mu_Y(W_{n,B,m_\bdk}^{t,q,j,\bdg,w,\bdr})
T^n \\
=&\bl^{-(\nu_\bdg+A_\bdg)}
\sum_{\bdr\in\bp^d}\sum_{n\geq 1}
\bl^{-{\tilde m}_\bdk} \mu_Y(W_{n,B,{\tilde m}_\bdk}^{t,q,j,\bdg^1,w,\bdr})T^n.
\end{split}
\end{equation}

Consider the following Newton map associated
with the Newton component of $h^t$ (as a $K$-quasi-ordinary power series)
we are dealing with. Its
Newton polyhedron is defined by
$z^{n_1^q}-\bbeta_{j,w}^{q}(t)x_1^{b_1^q}
\ldots x_d^{b_d^q}$,
where ${\bbeta}_{j,w}^{q}(t)=\beta_w t^{\alpha_{q,j}}+\text{higher order terms}\in k[[t]]$
and ${\beta}_w\in \bg.$ We can write
${\bbeta}_{j,w}^{q}(t)=t^{\alpha_{q,j}}\alpha(t)$
being $\alpha(t)=(\alpha_w+\ldots)^{n_1^q}\in k[[t]]$
because $(\alpha_w)^{n_1^q}=\beta_w.$

Let $\pi_{j,t}^{q,w}$ be the $k[[t]]$-morphism associated
to the transformation  $y_l={\tilde y}_l^{p^q_l}$ for all
$i=1,\ldots,d$ and $z_1=(z_2+\alpha(t))\prod_{l=1}^d {\tilde y}_l
^{{\bar b}_l^q}.$ Recall  the change of
variables formula is also applied here, see Remark \ref{kt}.

Taking into account the pull-back of the discriminants
and applying the same proof as in Step \ref{unos},
the following proposition follows.

\begin{prop}
Let $h^t_{q,j,\bdg,w}:=
h^{\{t\}}\circ \pi_{j,t}^{q,w}=
h^t\circ \pi_{\bdg}\circ \pi_{j,t}^{q,w}$
and
${\tilde\omega}_{q,j,\bdg,w}=(\pi^{q,w}_{j,t})^*\omega_\bdg$.
Then:
\begin{enumerate}[\rm(1)] 

\item The series $h^t_{q,j,\bdg,w}$ is $\ww$-quasi-ordinary
and its  depth is strictly less than the one
of $h^t$. We have
\begin{equation} \label{caso2a}
Z^{\ww,\bdg,w}_{1,q,j}(T)=\bl^{-(\nu_\bdg+A_\bdg)}
Z^\ww_{DL}(h^t_{q,j,\bdg,w},{\tilde\omega}_{q,j,\bdg,w},T).
\end{equation}

\item If $h^t$ is $\tww$-quasi-ordinary and Newton compatible then there exists a weight $\tww{\,}'$ such that $h^t_{q,j,\bdg,w}$
is $\tww{\,}'$-quasi-ordinary and Newton compatible.

\item The series $h^t_{q,j,\bdg,w}$ is $J_q$-bounded, see Remark \ref{notodos}.
\end{enumerate}
\end{prop}

\begin{proof} The only non trivial facts are the two last ones. 
It is easy to prove from the Newton map that $h^t_{q,j,\bdg,w}$
is $\tww{\,}'$-quasi-ordinary, where
$$
\tww{\,}':=(p n_1^q;(p v_1+g_1) c_1^q,\dots,(p v_d+g_d) c_d^q).
$$
This is due to the following equalities:
$$
\alpha_{q,j}=\frac{1}{p}\sum_{l=1}^d b_l^q g_l,\quad
c^\bdg=\frac{\alpha_{q,j}}{n_1^q}+\sum_{l=1}^d \frac{v_l}{p_l^q} {\bar b}_l^q=
\frac{\alpha_{q,j}}{n_1^q}+\sum_{l=1}^d \frac{v_l}{n_1^q} b_l^q=
\frac{1}{n_1^q}\sum_{l=1}^d (v_l+\frac{g_l}{p}) b_l^q.
$$
The $\tww{\,}'$ depends only on $\tww$ and on the Newton map. It implies that we can choose a \emph{generic}  $\tww$-quasi-ordinary series $\tilde h^t$ with the same characteristic monomials and then Newton compatible.
Then, $\tilde h^t_{q,j,\bdg,w}$ is also $\tww{\,}'$-quasi-ordinary.
If we pass to good coordinates the weight $\tww{\,}'$ remains and
both power series have the same Newton polyhedron. For $\tilde h^t_{q,j,\bdg,w}$ all the coefficients in its restriction
to the Newton polyhedron will be non zero and it is easily seen that it implies that $\tilde h^t_{q,j,\bdg,w}$ is Newton compatible and it is also the case for $h^t_{q,j,\bdg,w}$.

The fact that $h^t_{q,j,\bdg,w}$ is $J_q$-bounded follows from the behaviour of the $t$-morphism $\pi_\bdg$; note that each variable of $J_q$ is replaced by this variable times a power of $t$; then for a fixed power of $t$ only a finite number of monomials in $t$ and the variables of $J_q$ can contribute.
\end{proof}

Consider now arcs in Case (2.\ref{caso2}).
Let $\emptyset \subsetneq J \subset \{1,\ldots,d\}$;
in fact, from Remark \ref{notodos1}, we are only interested
in the case $\emptyset \subsetneq J \subset J_q$.
Let us fix a point
$(\bdy,0)$ corresponding to the subset $J$.
Since $h^t_{q,j,\bdg,w}$ is $\ww$-quasi-ordinary $J$-bounded
power series  the germ of $h^t_{q,j,\bdg,w}$ at $(\bdy,0)$ is also
a $\ww$-quasi-ordinary power series, see Remark \ref{translat}. We can indeed evaluate at $(\bdy,0)$
because the coefficients of its monomials in $t$
belong to $k[\widetilde{\bdy},z]$.
Since $h^t_{q,j,\bdg,w}$ is $\ww$-quasi-ordinary
its discriminant is computed from the pull-back of
$D_{z}(h^t)=t^\beta x_1^{\alpha_1}\ldots x_d^{\alpha_d}
\varepsilon (\bdx,t)$
where $\varepsilon(0,0)\ne 0.$ Since $\bdg\in \bp^{d+1}$
then $\tilde \varepsilon:=\varepsilon \circ \pi_{\bdg}\circ
\pi_{t,j}^{q,w}(\tilde{\bdy},t)$
also verifies $\tilde \varepsilon (0,0)\ne 0.$
The characteristic
exponents do not depend on the particular
$\bdy$ but only on $J$ and the depth has decreased.
Namely, by induction $Z_{DL}^\ww((h^t_{q,j,\bdg,w})_{(\bdy,0)},
({\tilde\omega}_{q,j,\bdg,w})_{(\bdy,0)},T)$
satisfies the theorem.
Observe also that the variety of such $\bdy$'s is
${\mathbb G}_{m,k}^{\#J}$.

\begin{prop} The following equality holds:
\begin{equation} \label{caso2b}
Z_{J,q,j}^{\ww,\bdg,w}(T)=
(\bl-1)^{\#J} \bl^{-(\nu_\bdg+A_\bdg)}
Z_{DL}^\ww((h^t_{q,j,\bdg,w})_{(\bdy,0)},({\tilde\omega}_{q,j,\bdg,w})_{(\bdy,0)},T),
\end{equation}
where ${}_{(\bdy,0)}$ means the germ at this point.
\end{prop}
\begin{proof}
This result is a consequence of Proposition \ref{fam-qo}.
\end{proof}

\begin{obs}
We note that
$(h^t_{q,j,\bdg,w})_{(\bdy,0)}$ is a $\ww$-quasi-ordinary series
in good coordinates  but its  depth is strictly less than
$\dpt(h^t_{q,j,\bdg,w})$.
Moreover its Newton polyhedron is obtained projecting
the Newton polyhedron of $h^t_{q,j,\bdg,w}$ on the corresponding coordinate planes.
\end{obs}

\begin{paso}
Computations for $h$ and $\bdg=\mu_1^\bdg\bdw_1^q+\ldots+\mu_d^\bdg\bdw_d^q\,\in G^q$
such that $0<\mu_l^\bdg<1$ for some $l\in\{1,\dots,d\}$.
\end{paso}

This case can be computed as in step \ref{doses} since
we may consider $h$ as $\tww$-quasi-ordinary where
$\tww=(1;0,\dots,0)$.

We are going to collect all the results in this section 
to compare motivic  zeta functions for the
functions $h^t$ and $h$ in order to obtain recurrence formul{\ae}.

By the identities in (\ref{induction}),
\begin{equation} \label{despues}
\begin{split}
Z_{DL}^\ww(h^t,\omega,T)&=Z_A^\ww(h^t,\omega,T)+ \\
&+\sum_{q=1}^r\underset{s_1^{j,q}=1}{\sum_{j=1}^{u(q)}}
\left( \sum_{\bdg \in G^{q,\alpha_{q,j}}}
\bl^{-(\nu_\bdg+A_\bdg)} (\,\sum_{w=1}^{w(q,j)}
(Z_{DL}^\ww (h^t_{q,j,\bdg,w},
{\tilde \omega}_{q,j,\bdg,w}^t,T)+\right. \\
&\left.+\sum_{\emptyset \subsetneq J \subset
\{1,\ldots,d\}}
(\bl-1)^{\#J}
Z_{DL}^\ww((h^t_{q,j,\bdg,w})_{(\bdy,0)},
({\tilde\omega}_{q,j,\bdg,w})_{(\bdy,0)},T))\right),
\end{split}
\end{equation}
where $\nu_\bdg+A_\bdg=\frac{\alpha_{q,j}}{n_1^q}+\sum_{l=1}^d
\mu_l^\bdg(p_l^q\nu_l+{\bar b}_l^q).$
In the same way,
\begin{equation}
\label{despues2}
\begin{split}
Z_{DL}(h,{\bar \omega},T) & =Z_{DL}^A(h,{\bar \omega},T)+
\sum_{q=1}^r \sum_{{\bar j}=1}^{v(q)} (
Z_{DL}({\bar h}_{q,{\bar j}},{\bar \omega}_{q,{\bar j}},T) \\
&+\sum_{\bdg \in G^{q},\bdg\ne \bdg^1 }
(\bl^{-(\nu_\bdg+A_\bdg)} (
 Z_{DL}^\ww(\tilde h^t_{q,\bdg,{\bar j}},
\tilde\omega_{q,\bdg,{\bar j}},T)+ \\
&+\sum_{\emptyset \subsetneq J \subset \{1,\ldots,d\}}
(\bl-1)^{\#J}
Z_{DL}^\ww((\tilde h^t_{q,\bdg,{\bar j}})_{(\bdy,0)},
(\tilde\omega_{q,\bdg,{\bar j}})_{(\bdy,0)},T))).
\end{split}
\end{equation}
where $v(q)=\sum_{j=1}^{u(q)} w(q,j)$ and
$\nu_\bdg+A_\bdg=\sum_{l=1}^d
\mu_l^\bdg(p_l^q\nu_l+{\bar b}_l^q).$

\begin{proof}[Proof of Theorem \ref{thm-polos}]
The starting point of the induction has been done
in Step \ref{start-ind}.

Let us attack the general case.
Since $h^t$ and $h$ (as $K$-quasi-ordinary power series) 
have the same characteristic exponents then
by Lemma \ref{charexp} there exist one-to-one correspondences
between the sets of faces of their Newton polyhedra and their
 sets of roots on such faces for $h^t$ and $h$. 
 Assume that under these bijections 
 the root $\bbeta_{j,w}^q(t)$, 
 where $s_1^{q,j}=1$, corresponds to ${\bar \beta}^q_{\bar j}$,
 namely the polynomial $z^{n_1^q}-\bbeta_{j,w}^q(t) x_1^{b_1^q}\dots x_d^{b_d^q}$
corresponds to ${\bar z}^{n_1^q}-
{\bar \beta}_{\bar j}^q {\bar x}_1^{b_1^q}\dots
{\bar x}_d^{b_d^q}$.
Let $h^t_{\bbeta}$ and ${\bar h}_{q,{\bar j}}$ denote the pull-back
of $h^t$ and $h$ under the corresponding Newton maps
associated with the roots $\bbeta_{j,w}^q(t)$ and ${\bar \beta}^q_{\bar j}$.
\begin{enumerate}

\item By Lemma \ref{charexp},
the characteristic exponents of $h^t_{\bbeta}$
and ${\bar h}_{q,{\bar j}}$ are equal.
\item For each $\bdg\in G^q$ (resp. $\bdg \in G^{q,\alpha_{q,j}}$)
the $K$-quasi-ordinary power series $h\circ \pi_\bdg$ (resp.
$h^t\circ \pi_\bdg$) has the same characteristic exponents as
$h$ (resp. $h^t$). This is Remark \ref{Achange}.
\item Applying (1) and (2)
for each $\bdg\in G^q\setminus \{\bdg^1\}$, then
${\tilde h}_{q,\bdg,{\bar j}}$ has the same
characteristic exponents as  ${\bar h}_{q,{\bar j}}$
(recall that ${\bar h}_{q,{\bar j}}$ is obtained
in Step \ref{unos}).

\item In the same way, for each $\bdg\in G^{q,\alpha_{q,j}}$,
${h}^t_{q,j,\bdg,w}$ has the same
characteristic exponents as $h^t_{\bbeta}$ which has the same as
${\bar h}_{q,{\bar j}}$. In all these cases, the depth of all
these quasi-ordinary power series
are smaller than  $\dpt(h)$. In particular we can follow
by induction.

\item For the statement of the Theorem the terms $\bl^{-\nu}$
are not essential then we have that
$$Z_{DL}^\ww(\tilde h^t_{q,\bdg,{\bar j}},
\tilde\omega_{q,\bdg,{\bar j}},T)\in \bz [\bl,\bl^{-1}, (1-\bl^{-\nu}T^N)^{-1}][T],
$$
where $(N,\nu)\in CP ({\bar h}_{q,{\bar j}},{\bar \omega}_{q,{\bar j}} ).$
\end{enumerate}
\end{proof}

\begin{proof}[Proof of Theorem \ref{thm-top}]
The starting point of the induction was done
in Step \ref{start-ind}.

For the general case we consider the formul{\ae} (\ref{despues})
and (\ref{despues2}). In this case the restriction
$s_1^{j,q}=1$ is empty and the series appearing in the second line of (\ref{despues}) are ${\tww\,}'$-quasi-ordinaries for some
weights and they have  depth less than $h$. The series in the third line are only $\ww$-quasi-ordinaries.

Because of Theorem \ref{thm-polos} we can compute 
$\chi_{\text{\rm top}}(\bullet,(\bl^{-s}))$ for each term.
For the second line terms we deduce by induction that
$$
\chi_{\text{\rm top}}(\bl^{-(\nu_\bdg+A_\bdg)}
Z^\ww_{DL}(h^t_{q,j,\bdg,w},{\tilde\omega}_{q,j,\bdg,w},
(\bl^{-s})))=
\zlo({\bar h}_{q,{\bar j}},\omega_{q,{\bar j}},s).
$$
By Theorem \ref{thm-polos} and (\ref{caso2b}),
 the terms in the third line  vanish
when applying $\chi_{\text{\rm top}}(\bullet(\bl^{-s}))$.
\end{proof}

\subsection{Zeta functions along strata}
\mbox{}

To finish the proof of  Theorems \ref{thm-polos}
and \ref{thm-top}
we must prove Proposition \ref{caso2b}.
We will proceed by induction in \emph{families}
of $\ww$-quasi-ordinary power series.

\begin{defini} Let $h^t\in k[[t]][[\bdx,\bdy]][z]$ be
a $\ww$-quasi-ordinary power series, $\bdx=(x_1,\dots,x_d)$,
$\bdy=(y_1,\dots,y_e)$ which is bounded for the $\bdy$-variables. We say that
$h^t$ is \emph{a family of $\ww$-quasi-ordinary power series}
if for any $\bdy^0\in{\mathbf G}_{k,m}^e$, the well-defined germ
$h^t_{\bdy^0}:=h^t(\bdx,\bdy+\bdy^0)$
is also $\ww$-quasi-ordinary power series with the
same characteristic monomials.
\end{defini}

\begin{defini} The Denef-Loeser zeta function of a
family of $\ww$-quasi-ordinary power series
and a form $\omega$
satisfying the support condition (\ref{spc}) is defined
by
\begin{equation*}
\begin{split}
Z_{DL}^{\ww,\bdy}(h^t,\omega,T):=T^\theta \sum_{n\in\bn} \left(\sum_{m\in\bn}
\bl^{-m}\mu_X(V_{n,m}^t)\right)T^n\in \widehat{\mm}_k[[T]] ,\\
Z^{\ww}_{\text{top},0}(h^t,\omega,s):=
\chi_{\text{top}}(Z^{\ww}_{DL}(h^t,\omega,\bl^{-s})).
\end{split}
\end{equation*}
where we are considering 
arcs having their origin along $\{0\}\times {\mathbf G}_{k,m}^e$.
\end{defini}

\begin{obs} The powers of the variables in $\bdy$ for $\omega$ are
negligible.
\end{obs}

\begin{lema} $Z_{DL}^{\ww,\bdy}(h^t,\omega,T)=Z_{DL}^{\ww,\bdy}(\tilde h^t,\omega,T)$
where $\tilde h^t(\bdx,\bdy)=h^t(\bdx,\bdy^\bdn)$.
\end{lema}

\begin{prop}\label{fam-qo} Let $h^t$ be a
family of $\ww$-quasi-ordinary power series and $\omega$ a form
verifying the support condition (\ref{spc}). Then for all
$\bdy^0\in {\mathbf G}_{k,m}^e$,

$$
Z_{DL}^{\ww,\bdy}(h^t,\omega,T)=
(\bl-1)^e  Z_{DL}^{\ww}(h^t_{\bdy^0},\omega,T).
$$
\end{prop}

\begin{proof}
The proof is by induction on the depth of $h^t$;
 if $\dpt h^t=0$, the result is trivial. We compute $Z_{DL}^{\ww,\bdy}(h^t,\omega,T)$
following the Newton polyhedron and the Newton maps.

We begin with the $A$-part. If we consider a vertex, it is related
to a monomial $t^\alpha\bdx^\bdn\bdy^\bdm z^n$.
It is clear that both terms in the equality differ by
an $(\bl-1)^e$ factor.
Let us fix an edge of the Newton polytope and a factor
of the principal part which will give the key for
the Newton map, i.e., we fix some factor
$(z^n-\beta t^\alpha \bdx^\bdn\bdy^\bdm)^m$.

The mapping
$$
(\bdx,\bdy,z)\mapsto(\bdx,y_1^n,\dots,y_e^n,y_1^{m_1}\dots y_e^{m_e} z)
$$
does not change $Z_{DL}^{\ww,\bdy}(h^t,\omega,T)$. Thus we may suppose
that the factor is $(z^n-\beta t^\alpha \bdx^\bdn)^m$.
Since the computations for the $A$-part depend only on this
Newton principal part,
it is clear that both terms in the equality differ by
an $(\bl-1)^e$ factor.

Let us consider now the $B$-part. For the sake of simplicity
we can suppose that all coordinates of $\bdn$ are non zero.
Following the argument of Step \ref{doses}, the zeta function
corresponding to the $B$-part is the zeta function
associated with a power series $\tilde h^t$ and a form $\tilde\omega$
where we must consider arcs based at
$\ba_k^d\times {\mathbf G}_{k,m}^e$. We can perform a stratification
of this space and decompose the result in a sum of
zeta functions of families of $\ww$-quasi-ordinary power series,
where variables corresponding to $J\subset\{1,\dots,d\}$,
pass to $\bdy$. Since these power series have strictly less depth
we can also obtain the $(\bl-1)^e$ factor.
\end{proof}
The Theorems \ref{thm-polos} and \ref{thm-top} are proved.

\section{Consequences of the main theorems}

\subsection{Essential variables}
\mbox{}

Let $h\in k[[\bdx]][z]$ be a quasi-ordinary
power series. For any subset
$\emptyset\subsetneq I\subsetneq \{1,\ldots,d\},$ let $k_I$ be an
algebraic closure of the quotient field of the domain
$k[[x_i]]_{i\in I}.$ Thus
$h_I\in k_I[[{\widehat {\bdx}}_I]][z]$
is a $k_I$-quasi-ordinary power series,
where ${\widehat\bdx}_I$ are the variables not in $I.$

Assume  $h\in k[[\bdx]][z]$ is a quasi-ordinary
function with $\lgt(h)=e<d$. Relabeling the variables we assume that $x_1,\ldots,x_e$
are the essential variables. Set $I=\{e+1,\ldots,d\}.$
Let $\omega=\prod_{l=1}^d x_j^{\nu_j -1} d\bdx \wedge dz$
be a differential form such that $(h,\omega)$
verifies the support condition \ref{spc}.
Whenever the
variable $x_j$ is not essential then the support condition
implies $\nu_j=1.$ In particular if
$\omega_I=\prod_{l=1}^e x_j^{\nu_j -1}
d{\widehat\bdx}_I \wedge dz$ then the pair $(h_I,w_I)$
also verifies the support condition \ref{spc} over the field $k_I$.

\begin{cor} \label{essen-var}
In the previous conditions, the following formal
identity holds:
$$
Z_{DL}(h,\omega,T)=\bl^{-(d+1-e)} Z_{DL} (h_I,w_I,T).
$$
\end{cor}

Here formal means that in the RHS, resp. LHS,
the varieties are over $k$, resp. $k_I.$
The proof of the lemma
follows essentially
from the proof of theorems \ref{thm-polos} and
\ref{thm-top}. In fact
Newton maps do affect the non-essential variables.

\subsection{Curves case}
\mbox{}

Suppose a $\ww$-quasi-ordinary power series
$h^t$ has $\lgt(h^t)=1$.
After relabeling the variables $x_i,$
 assume 
$D_z(h^t)=t^\alpha x_1^{\alpha_1} \varepsilon(t,\bdx),$
with $\varepsilon(0,0)\ne 0.$
We assume the same for the quasi-ordinary power series $h.$
Then the compact faces of their Newton polyhedra are
contained in $2$-dim plane $x_2=\ldots=x_d=0.$
In particular, in the dual
space one has the same property. The dual hyperplanes
$l_{q,j}^t$ are defined by
$\alpha_{q,j}+s_1^{q,j}(b_1^qv_1-n^q_1 v_{d+1})=0$,
where $\gcd(b_1^q,n_1^q)=1.$
The dual hyperplanes are
determined by the linearly independent
vectors $\bdw^q_1$ and $\bdw^q_l=e_l$ for $l=2,\ldots,d.$
Therefore the sets $G^q$ and $G^{q,\alpha_{q,j}}$ have
one element. Then we can conclude that
\begin{equation}\label{fcurvascase}
\begin{split}
Z_{DL}^\ww(h^t,\omega,T)=&Z_A^\ww(h^t,\omega,T)
+\sum_{q=1}^r
\underset{s_1^{q,j}=1}{\sum_{j=1}^{u(q)}}
\bl^{-(\nu_\bdg+A_\bdg)}\sum_{w=1}^{w(q,j)}
Z_{DL}^\ww (h^t_{q,j,\bdg,w},
{\tilde \omega}_{q,j,\bdg,w}^t,T),  \\
Z_{DL}(h,{\bar \omega},T)=&Z_{DL}^A(h,{\bar \omega},T)+\sum_{q=1}^r
\sum_{{\bar j}=1}^{v(q)}
Z_{DL}({\bar h}_{q,{\bar j}},\omega_{q,{\bar j}}, T),
\end{split}
\end{equation}
where $\nu_\bdg+A_\bdg=\frac{\alpha_{q,j}}{n_1^q}+
\mu_1^\bdg(p_1^q\nu_l+{\bar b}_1^q)$ and $v(q)=\sum_{j=1}^{u(q)} w(q,j).$

\begin{cor} \label{curvascase} For a quasi-ordinary power series 
$h$ with $\lgt(h)=1$ and a form $\omega$ such that $(h,\omega)$
verifies condition (\ref{spc}) the second line
in equations (\ref{fcurvascase}) can be used as recursive
formula to compute $Z_{DL}(h,{\bar \omega},T).$
\end{cor}
A general formula for the topological zeta function will
be given in Theorem \ref{mainformula}.

\subsection{The topological zeta function}
\label{formula}
\mbox{}

Assume the quasi-ordinary power series
$h(\bdx,z):=\prod_{l=1}^d x_l^{N_l}f(\bdx,z)u(\bdx,z)$,
$N_l\in\bn$,
is given in a good system of coordinates
and that the differential form
$\omega=\prod_{j=1}^d x_j^{\nu_j-1}d x_1\wedge\ldots\wedge d x_d\wedge dz, \,\nu_j\geq 1,$
verifies condition~(\ref{spc}).
Assume that $\gamma_1,\ldots,\gamma_r$ (resp. $\tau_0,\ldots,\tau_r$)
are the $1$-dimensional (resp. $0$-dimensional) compact
faces of its Newton polyhedron. Assume  for each  $q=1,\ldots,r,$
on  $\gamma_q$ there are exactly $v(q)$ distinct roots.
Let $\pi_{q,j}$ be the Newton map associated with
the Newton component $f_j^q$ of $f.$

\smallbreak

\begin{thm} \label{mainformula}
Under the above conditions the following equality holds
\begin{equation*}
\begin{split}
Z_{\text{\rm top},0}(h,w,s)=&\sum_{l=0}^r J_{\tau_l}(h,\omega,s)
-\sum_{q=1}^r
{v(q)}
J_{\gamma_q}(h,\omega,s)+\\
&+\sum_{q=1}^r
\mult(\Delta_{\gamma_q})\sum_{j=1}^{v(q)}
Z_{\text{\rm top},0}(h\circ\pi_{q,j},w\circ\pi_{q,j},s).
\end{split}
\end{equation*}

Moreover if $h$ has non-degenerated principal part then
for each $q=1,\ldots,r$ the integer number $v(q)$ is nothing
but $V(\gamma_q)$ and for any $j=1,\ldots,v(q),$ the
following identity is verified
$$
Z_{\text{\rm top},0}(h\circ\pi_{q,j},w\circ\pi_{q,j},s)=
\frac{J_{\gamma_q}(h,\omega,s)}{(s+1)\mult(\Delta_{\gamma_q})}.
$$
\end{thm}

\begin{proof} The same ideas as in the construction of
the above section
give $Z_{DL}(h,\omega,T)$ as
$$
Z_{DL}(h,\omega,T)=Z_{DL}^A(h,\omega,T)+
\sum_{q=1}^r  \sum_{j=1}^{v(q)}
\sum_{\bdg\in G^q} Z^{\gamma_q,\bdg,j}(T),
$$
 see equation (\ref{induction}).
The contribution of the A-part was given in equation (\ref{eqAtop}).

\begin{ctrb}
Fix a $1$-dim face $\gamma_q$, a root
$\beta_j^q$ of $\gamma_q$ and $\bdg$ in $G^q$.

\begin{enumerate}[(1)]
\item If $\bdg=\bdg^1$ then $
Z^{\gamma_q,\bdg^1,j}(T)=Z_{DL}(\tilde h,\tilde\omega,T),
$
where $ \tilde h=h\circ \pi_{q,j}$ and
$\tilde \omega=\omega\circ \pi_{q,j},$
see equation (\ref{caso11}) in Step \ref{unos}.

\item If $\bdg\ne \bdg^1$ then
$\displaystyle Z^{\gamma_q,\bdg,j}(T)=Z_1^{\gamma_q,\bdg,j}(T)+
\sum_{\emptyset\subsetneqq J\subset\{1,\dots,d\}}
Z_J^{\gamma_q,\bdg,j}(T)$, see Step \ref{doses}.
\end{enumerate}
\end{ctrb}

Recall that the following  properties are also verified.
\begin{enumerate}
\item  $\chi_{\text{\rm top}}(Z_1^{\gamma_q,\bdg,j}(\bl^{-s}))=
\zlo(h\circ \pi_{q,j},\omega\circ \pi_{q,j},s),$ see equation (\ref{caso2a}).
\item $\chi_{\text{\rm top}}(Z_J^{\gamma_q,\bdg,j}(\bl^{-s}))=0$, see
equation (\ref{caso2b}).
\end{enumerate}

Hence
the contribution to $Z_{\text{\rm top},0}(h,w,s)$ of the $B$-part
is
$$ \sum_{q=1}^r
\sum_{j=1}^{v(q)} \sum_{\bdg\in G^{q}}
\zlo(h\circ \pi_{q,j},\omega\circ \pi_{q,j},s).$$
Finally the cardinality
of $G^q$ is
mult$(\Delta_{\gamma_q})$.

The second part of the theorem follows from the results of
section~\S~\ref{secnondeg}, in particular from Lemma \ref{lemab}.
\end{proof}

\begin{algthm} Our method gives
an effective algorithm
to compute $Z_{\text{\rm top},0}(h,w,s)$.
After a Newton map $\pi_{q,j}$, 
$h\circ \pi_{q,j}(\bdy,z_1)=y_1^{N_1}\cdots y_d^{N_d}
\tilde{f}(\bdy,z_1)\tilde{u}(\bdy,z_1)$
where $\tilde{f}\in k[[\bdy]][z_1]$
is a $z_1$-quasi-ordinary power series.
Hence
$\tilde{f}(\mathbf 0)=0$ and
$\displaystyle
N_l=\inf_{\bdx\in \Gamma(h)} \{\bdw_l\cdot \bdx\}
$
for any $l=1,\ldots,d$. Apply again
Theorem \ref{mainformula} until depth $0$.
This algorithm can be implemented in Maple based on a
Maple program made by  K. Hoornaert and D. Loots
for the non-degenerated case, \cite{kt:00}.
\end{algthm}

\medskip
\subsection{A special candidate pole}
\mbox{}

In this paragraph we will show that some candidate poles
disappear when we add their local contributions.
For a candidate pole $s_0=(N,\nu)$ coming
from a compact $1$-dimensional face $\gamma,$
its \emph{local contribution}
is defined as
the sum of the $A$-part of the local Denef-Loeser motivic
zeta function
which corresponds to $\gamma$ and its two vertices
plus the sum of the contributions of the $B$-part.
The latter is nothing but the contribution of
the highest vertex in the new Newton polyhedron
after all Newton maps of Newton components associated with
$\gamma$ and at all
possible new quasi-ordinary power series.

Consider the following formal power series
\begin{equation*}
\begin{split}
h^t(\bdx,z):=&x_1^{N_1}\dots x_{d-1}^{N_{d-1}}
z^{n-n_1 m}(z^{n_1}-\beta t^\alpha x_1^{b_1}\dots x_{d-1}^{b_{d-1}} x_d)^m+\dots,\\
h(\bdx,z):=&x_1^{N_1}\dots x_{d-1}^{N_{d-1}}
z^{n-n_1 m}(z^{n_1}-\beta
x_1^{b_1}\dots x_{d-1}^{b_{d-1}} x_d)^m+\dots,
\end{split}
\end{equation*}
where the other terms appears \emph{behind}
their Newton polyhedra and $\alpha\in \bp,
\beta\in {\mathbf G}_{m,k}.$
Assume  $h^t\in k[[t]][[\bdx]][z]$ is
$\ww$-quasi-ordinary and $h\in k[[\bdx]][z]$ is
a quasi-ordinary power series.

Consider the differential form
$\omega:=x_1^{\nu_1-1}\dots x_{d-1}^{\nu_{d-1}-1}
d x_1\wedge\dots\wedge d x_{d}\wedge d z,$ in such a way $(h^t,\omega)$
and $(h,\omega)$ verify
the support condition \ref{spc}.
In other words the pair
$(N_d,b_d)=(0,1).$
In this section 
the local contribution of the candidate pole
$s_0:=1-\bl^{-(n_1+1)} T^n$ which corresponds with
the $d$-th coordinate of the first characteristic exponent
of $h^t$ and $h$ is computed. Under the above conditions we have the following
result.

\begin{prop}\label{tspecial} If there is no pair
$(N_i,b_i)$ equals to $(0,1)$ other than $(N_d,b_d)$
then $s_0=1-\bl^{-(n_1+1)} T^n$ does not appear in the denominator
of the local contribution of $s_0$ to
$Z_{DL}^{\ww}(h^t,\omega,T)$ and $Z_{DL}(h,\omega,T).$
\end{prop}

\begin{proof}
The compact faces of $\Gamma(h)$ and
$\Gamma^t (h^t)$ are the same and
consist of a face
$\gamma$ with
vertices $\tau_0:=\!(N_1+m b_1,\dots,
N_{d-1}+ m b_{d-1}, m,n-n_1 m)$ and
$\tau_1:=\! (N_1,\dots,N_{d-1},0,n)$.

Let $\eta:=
b_1 v_1+\dots b_{d-1} v_{d-1}+ v_d- n_1 v_{d+1}.$
In the dual space, the dual decomposition induced by
$\Gamma(h)$ has only three
convex rational cones
$\Delta_{\tau_1}=\{\eta>0\},$
$\Delta_{\tau_0}=\{\eta<0\}$ and
$\Delta_{\gamma}=\{\eta=0\}.$
In the same way, the decomposition
induced by $\Gamma^t(h)$ is given
by the convex rational polyhedron
$\Delta^t_{\tau_1}=\{\eta+\alpha>0\},$
and the cones
$\Delta^t_{\tau_0}=\{\eta+\alpha<0\}$ and
$\Delta^t_{\gamma}=\{\eta+\alpha=0\}.$
Since there is only one exponent $\alpha$ of $t$, convex polyhedra
delimited by parallel hyperplanes do not exist.

To compute the contributions of the candidate pole
$s_0:=1-\bl^{-(n_1+1)} T^n$
we multiply $Z_{DL}^{\ww}(h^t,\omega,T)$ and $Z_{DL}(h,\omega,T)$
by $s_0$
and simplify the result under the condition
$\bl^{-(n_1+1)} T^n=1$.

The contribution coming from
the Newton polyhedra of $h$ and $h^t$
written in terms of generating functions
are:
$$
Z_{DL}^{\Delta_{\tau_0}^t,\ww}(T)= T^{m\alpha} (\bl-1)^{d+1}
\bl^{-(d+1)} \Phi_{\Delta_{\tau_{0}}^{t}}({\bdy}),
$$
$$
Z_{DL}^{\Delta_{\tau_1}^t,\ww}(T)= (\bl-1)^{d+1}
\bl^{-(d+1)} \Phi_{\Delta_{\tau_{0}}^{t}}(\bda),
$$
$$
Z_{DL}^{\Delta_{\gamma}^t,\ww}(T)=
\bl^{-(d+1)} \left((\bl-1)^{d+1}-(\bl-1)^{d}\right)
\Phi_{\Delta_{\gamma}^{t}}(\bda),
$$
$$
Z_{DL}^{\Delta_{\tau_0}}(T)=  (\bl-1)^{d+1}
\bl^{-(d+1)} \Phi_{\Delta_{\tau_{0}}}({\bdy}),
$$
$$
Z_{DL}^{\Delta_{\tau_1}}(T)= (\bl-1)^{d+1}
\bl^{-(d+1)} \Phi_{\Delta_{\tau_{0}}}(\bda),
$$
$$
Z_{DL}^{\Delta_{\gamma}}(T)= 
\bl^{-(d+1)} \left((\bl-1)^{d+1}-(\bl-1)^{d}\right)
\Phi_{\Delta_{\gamma}}(\bda),
$$
where
$\bda:=
(a_1,\dots,a_{d-1},a_d,a_{d+1})=(\bl^{-\nu_1} T^{N_1},\dots,
\bl^{-\nu_{d-1}} T^{N_{d-1}},\bl^{-1},\bl^{-1} T^n)
$ and
$$
\bdy:=
(y_1,\dots,y_d,y_{d+1})=(\bl^{-\nu_1} T^{N_1+m b_1},\dots,\bl^{-\nu_{d-1}} T^{N_{d-1}+m b_{d-1}},\bl^{-1} T^m,\bl^{-1} T^{n-n_1 m}).
$$
If $\Delta_{\mathbf 1}$ denotes the positive cone $\bp^{d+1}$ then the
following identities among indicator functions hold:
$[\Delta_{\tau_1}^t]=[\Delta_{\mathbf 1}]-
[\Delta_{\tau_0}^t]-[\Delta_{\gamma}^t]$
and $[\Delta_{\tau_1}]=[\Delta_{\mathbf 1}]-
[\Delta_{\tau_0}]-[\Delta_{\gamma}]$.
Summing and looking only at the terms
where $1-\bl^{-(n_1+1)} T^n$ appears, the contribution of $s_0$
is given by
\begin{equation} \label{contribuA}
\bl^{-(d+1)}(\bl-1)^{d+1}\left(
-A^t(\bda)+T^{m\alpha}A^t(\bdy)-
\frac{\bl^{-1}\Phi_{\Delta_{\gamma}^t}(\bda)}{(1-\bl^{-1})}\right),
\,\text{ for $h^t$ and}
\end{equation}
\begin{equation*}
\bl^{-(d+1)}(\bl-1)^{d+1}\left(
-A(\bda)+A(\bdy)-
\frac{\bl^{-1}\Phi_{\Delta_{\gamma}}(\bda)}
{(1-\bl^{-1})}\right),\,\text{ for $h$},
\end{equation*}
where
$$
A^t(\bdx):=\frac{\sum_{\bbeta\in G_0^t} {\bdx}^{\bbeta}}{(1-x_{d+1})
\prod_{l=1}^d \left(1-{\bdx}^{\bdw_l}\right)},\qquad
A(\bdx):=\frac{\sum_{\bbeta\in G_0} {\bdx}^{\bbeta}}{(1-x_{d+1})
\prod_{l=1}^d \left(1-{\bdx}^{\bdw_l}\right)},
$$
and $G_0^t$ (resp. $G_0$) is the fundamental set of
$\Delta_{\tau_0}^t$ (resp. $\Delta_{\tau_0}).$

We will compute  the terms
$$A_1^t:=(1-\bl^{-(n_1+1)} T^n)\left(
-A^t(\bda)+T^{m\alpha}A^t(\bdy)-\frac{\bl^{-1}
\Phi_{\Delta_{\gamma}^t}(\bda)}{(1-\bl^{-1})}\right)
$$ 
and
$A_1:=(1-\bl^{-(n_1+1)} T^n)\left(
-A(\bda)+A(\bdy)-\frac{\bl^{-1}\Phi_{\Delta_{\gamma}}
(\bda)}{(1-\bl^{-1})}\right)$
under the given condition $\bl^{-(n_1+1)} T^n=1$.

Consider the vertex $\tau_1$.
The simplification $\bl^{-(n_1+1)} T^n=1$
gives $a_d^{n_1} a_{d+1}=1$.
Define $w_j:=a_j^{p_j} a_{d+1}^{\bar b_j}$, $1\leq j\leq d-1,$ and $c_i:=\gcd(n_1,b_i)$
and $p_i:=\frac{n_1}{c_i}$, $\bar b_i:=\frac{b_i}{c_i}$,
$i=1,\dots,d-1$.
Thus
$
w_j=\bl^{-(\nu_{j} p_j+\bar b_j)} T^{p_j N_{j}+n \bar b_j}=
\bl^{( b_j-\nu_{j})p_j} T^{p_j N_{j}},
$
and $w_j=u_j^{p_j}$ with $u_j:=\bl^{ b_j-\nu_{j}} T^{N_{j}}.$

The elements of $\Delta_{\tau_0}$
are written as
$$
\sum_{j=1}^{d-1}\mu_j (p_j e_j+\bar b_j e_{d+1})
+\mu_d(n_1 e_d+e_{d+1})+\mu_{d+1} e_{d+1}\in\bp^{d+1},
\mu_i\in\bq_{>0}.
$$

The elements of $\Delta_{\tau_0}^t$
can be written as
$$
\sum_{j=1}^{d-1}\mu_j (p_j e_j+\bar b_j e_{d+1})
+\mu_d(n_1 e_d+e_{d+1})+\mu_{d+1} e_{d+1}\in\bp^{d+1},
\mu_i\in\bq_{>0}, \mu_{d+1}>\frac{\alpha}{n_1}.
$$
To parametrize the sets  $G_0^t$ and $G_0$
we may replace $\mu_d$ by $\mu_d+s$,
$s\in\bz$, because of the
substitution $a_d^{-n_1}=a_{d+1}$.
Thus the elements in $G_0^t$ (and in $G_0$ when  $\alpha=0$)
are obtained with
$$
(\mu_1,\dots,\mu_{d-1},\mu_d,\mu_{d+1})=
\left(\frac{i_1}{p_1},\dots,\frac{i_{d-1}}{p_{d-1}},-\frac{\alpha+ k+\sum_{j=1}^{d-1} i_j b_j}{n_1},\frac{\alpha+k}{n_1}\right),
$$
$i_j=1,\dots,p_j$, $j=1,\dots,d-1$, and $k=1,\dots,n_1$.
Therefore $A^t(\bda)$ (resp. $A(\bda)$)
multiplied by $(1-a_d^{n_1}a_{d+1})=1-\bl^{-(n_1+1)} T^n$
is nothing but:
$$
\frac{\displaystyle\sum_{i_1=1}^{p_1}\dots\sum_{i_{d-1}=1}^{p_{d-1}}
\sum_{k=1}^{n_1}  a_d^{-(\alpha+k+\sum_{j=1}^{d-1} i_j b_j)}\prod_{j=1}^{d-1} a_j^{i_j}}{\displaystyle(1-a_{d+1})\prod_{j=1}^{d-1}(1-a_j^{p_j} a_{d+1}^{\bar b_j})}=
\frac{\displaystyle a_d^{-\alpha}(\sum_{k=1}^{n_1} a_d^{-k})\prod_{j=1}^{d-1}\sum_{i_j=1}^{p_j}(a_j a_d^{-b_j})^{i_j}
}{\displaystyle(1-a_{d+1})\prod_{j=1}^{d-1}(1-a_j^{p_j} a_{d+1}^{\bar b_j})}=
$$
\begin{equation}\label{ar}
=-\frac{\displaystyle a_d^{-\alpha}(1-a_d^{-n_1}) \prod_{j=1}^{d-1}\left(a_j a_d^{-b_j}(1-(a_j a_d^{-b_j})^{p_j})\right)
}{\displaystyle(1-a_d)(1-a_{d+1})\prod_{j=1}^{d-1}\left((1-a_j a_d^{-b_j})(1-a_j^{p_j} a_{d+1}^{\bar b_j})\right)}.
\end{equation}
The first summand $A^t(\bda)(1-a_d^{n_1}a_{d+1})$
(resp. $A(\bda)(1-a_d^{n_1}a_{d+1})$ when $\alpha=0$)
is:
\begin{equation}
\label{primera0}
-\frac{\displaystyle a_d^{-\alpha} \prod_{j=1}^{d-1} u_j
}{\displaystyle(1-a_d)\prod_{j=1}^{d-1}(1-u_j)}.
\end{equation}

Consider the term from $\gamma.$ The elements of
$\Delta_\gamma^t$ (resp. $\Delta_\gamma$ when $\alpha=0$)
are
$$
\sum_{j=1}^{d-1}\mu'_j (p_j e_j+\bar b_j e_{d+1})+\ldots
+\mu'_d(n_1 e_d+e_{d+1})+\frac{\alpha}{n_1} e_{d+1}\in\bp^{d+1},
\mu'_i\in\bq_{>0}.
$$
After simplifying,
the fundamental set $G_\gamma^t$ (resp. $G_\gamma$)
of $\Delta_\gamma^t$ (resp. $\Delta_\gamma$)
is obtained from
\begin{equation}\label{ggamma}
(\mu'_1,\dots,\mu'_{d-1},\mu'_d)=
\left(\frac{i_1}{p_1},\dots,\frac{i_{d-1}}{p_{d-1}},-\frac{\alpha+\sum_{j=1}^{d-1} i_j b_j}{n_1}\right),
\end{equation}
$i_j=1,\dots,p_j$, $j=1,\dots,d-1$.

Since we have to multiply the factor by $\frac{\bl^{-1}}{1-\bl^{-1}}$
and $\bl^{-1}=a_d$, then the product
$(1-a_d^{n_1}a_{d+1}
)\frac{\bl^{-1}\Phi_{\Delta_{\gamma}^t}(\bda)}
{(1-\bl^{-1})}$ equals
$$
\frac{\displaystyle a_d\sum_{i_1=1}^{p_1}\dots\sum_{i_{d-1}=1}^{p_{d-1}}
a_d^{-(\alpha+\sum_{j=1}^{d-1} i_j b_j)}\prod_{j=1}^{d-1} a_j^{i_j}}{\displaystyle(1-a_d)\prod_{j=1}^{d-1}(1-a_j^{p_j} a_{d+1}^{\bar b_j})}=
\frac{\displaystyle a_d^{1-\alpha}
\prod_{j=1}^{d-1}\sum_{i_j=1}^{p_j}(a_j a_d^{-b_j})^{i_j}
}{\displaystyle(1-a_d)
\prod_{j=1}^{d-1}(1-a_j^{p_j} a_{d+1}^{\bar b_j})}=
$$
\begin{equation}
\label{segunda}
=\frac{\displaystyle a_d^{1-\alpha} \prod_{j=1}^{d-1}\left(a_j a_d^{-b_j}(1-(a_j a_d^{-b_j})^{p_j})\right)
}{\displaystyle(1-a_d)\prod_{j=1}^{d-1}\left((1-a_j a_d^{-b_j})(1-a_j^{p_j}
a_{d+1}^{\bar b_j})\right)}=
\frac{\displaystyle a_d^{1-\alpha} \prod_{j=1}^{d-1}u_j
}{\displaystyle(1-a_d)\prod_{j=1}^{d-1}(1-u_j)}.
\end{equation}
The term $(1-a_d^{n_1}a_{d+1}
)\frac{\bl^{-1}\Phi_{\Delta_{\gamma}}(\bda)}
{(1-\bl^{-1})}$ is obtained from equation (\ref{segunda}) doing
$\alpha=0.$

For the vertex $\tau_0,$ after the
simplification, 
$y_d^{n_1} y_{d+1}=\bl^{-(n_1+1)} T^n=1$, $y_j y_d^{-b_j}=u_j$
and $y_j^{p_j} y_{d+1}^{\bar b_1}=w_j$, for $j=1,\ldots,d-1.$
The term $A^t(\bdy)$ multiplied by $1-\bl^{-(n_1+1)} T^n=
(1-y_d^{n_1}y_{d+1})$
is the same as (\ref{ar}) but $a_i$ replaced by $y_i$. Therefore
the term we are interested in is
\begin{equation}\label{tercera}
-\frac{\displaystyle T^{m\alpha} y_d^{-\alpha} \prod_{j=1}^{d-1} u_j
}{\displaystyle(1-y_d)\prod_{j=1}^{d-1}(1-u_j)}=
-\frac{\displaystyle a_d^{-\alpha} \prod_{j=1}^{d-1} u_j
}{\displaystyle(1-y_d)\prod_{j=1}^{d-1}(1-u_j)},
\end{equation}
because $T^{-1}y_d=a_d$. Furthermore $A(\bdy)(1-\bl^{-(n_1+1)} T^n)$
is obtained from equation (\ref{tercera})
when $\alpha=0.$ Thus
$$
A_1^t=(\ref{tercera})-(\ref{primera0})-(\ref{segunda})=
-\frac{\displaystyle a_d^{-\alpha} y_d \prod_{j=1}^{d-1} u_j
}{\displaystyle (1-y_d)\prod_{j=1}^{d-1}(1-u_j)}.
$$
And substituting $\alpha=0$ in the last identity
we get $A_1:$
\begin{equation}
\label{tercera1}
A_1=
-\frac{\displaystyle  y_d \prod_{j=1}^{d-1} u_j
}{\displaystyle (1-y_d)\prod_{j=1}^{d-1}(1-u_j)}.
\end{equation}

Finally we compute the contribution of the $B$ part in the
arc decomposition for both power series.
Let $\pi_\gamma$ be the map defined by
$x_l={\bar x}_l^{p_l}$ and $z=(z_1+\bar \beta)
\prod_{l=1}^d {\bar x}_l^{{\bar b}_l}$
where $\bar \beta^{n_1}=\beta.$
We compute the contribution after the Newton map
$\pi_\gamma$
associated with the
unique Newton component of $h$ and the corresponding
$k[t]$-morphisms associated with each $\bdg\in G_\gamma$.
For a pair $(h,\omega)$ in equation
(\ref{despues2})
we show that the contribution is
\begin{equation}\label{despB}
\begin{split}
Z_{DL}({\bar h},{\bar \omega},T)+&
\sum_{\bdg \in G_\gamma,\bdg\ne \bdg^1 }
\bl^{-(\nu_\bdg+A_\bdg)}
 Z_{DL}^\ww(\tilde h^t_{\bdg},
\tilde\omega_{\bdg},T) +\\
+\sum_{\bdg \in G_\gamma,\bdg\ne \bdg^1 }  \bl^{-(\nu_\bdg+A_\bdg)}
&\sum_{\emptyset \subsetneq J \subset \{1,\ldots,d-1\}}
(\bl-1)^{\#J}
Z_{DL}^\ww((\tilde h^t_{\bdg})_{(\bdy,0)},
(\tilde\omega_{\bdg})_{(\bdy,0)},T),
\end{split}
\end{equation}
where the pair
$(\bar h,\bar \omega)$ is the pull-back by $\pi_\gamma$
of $(h,\omega),$ that is
$$\bar h={\bar x}_1^{\tilde N_1}\dots
{\bar x}_d^{\tilde N_d}(z_1^{m_r}+\ldots)
u({\bar \bdx},z_1), \quad u({\mathbf 0},0)\neq 0,
$$
with $\tilde N_l=p_lN_l+{\bar b}_ln,\, l=1,\dots,d,$
and $\bar\omega=
{\bar x}_1^{{\tilde \nu}_1-1}\dots {\bar x}_d^{{\tilde \nu}_d-1}
d{\bar \bdx}\wedge dz_1.$

For $(\tilde h^t_{\bdg},\tilde\omega_{\bdg})$
the form $\tilde\omega_{\bdg}$ is nothing but $\bar\omega.$
Fix $\bdg\in  G_\gamma\setminus \{\bdg^1\}.$
The map $\pi_\gamma\circ \pi_{\bdg}$ is
$x_l=t^{\mu_l^\bdg p_l} {\bar x}_l^{p_l}$ and $z=t^{c^\bdg}(z_1+\bar \beta)
\prod_{l=1}^d {\bar x}_l^{{\bar b}_l}.$
Then $\tilde h^t_{\bdg}=h\circ \pi_\gamma\circ \pi_{\bdg}$,
that is
$$\tilde h^t_{\bdg}:=t^{nc^{\bdg} +\sum_{l=1}^{d-1}\mu_l^\bdg p_l N_l }{\bar x}_1^{\tilde N_1}\dots
{\bar x}_d^{\tilde N_d}(z_1^{m_r}+\ldots)
u({\bar \bdx},z_1), \quad u({\mathbf 0},0)\neq 0.$$

For $h^t$ 
the maps
associated with the
  corresponding
$k[t]$-morphisms associated with each $\bdg\in G_\gamma^t$ appear.
Looking at equation
(\ref{despues}) we have $r=1,u(r)=1$ and $s_1^{1,1}=1.$
It shows that the contribution is
\begin{equation}\label{despBt}
\begin{split}
\sum_{\bdg \in G_\gamma^t}
\bl^{-(\nu_\bdg+A_\bdg)}& \left(
 Z_{DL}^\ww(\tilde h^t_{\gamma,\bdg},
\tilde\omega_{\gamma,\bdg},T) +
\vphantom{\sum_{\emptyset \subsetneq J \subset \{1,\ldots,d-1\}}}\right.\\
+\sum_{\emptyset \subsetneq J \subset \{1,\ldots,d-1\}}
(\bl-1)^{\#J}&\left.\vphantom{\sum_{\emptyset \subsetneq J \subset \{1,\ldots,d-1\}}}
Z_{DL}^\ww((\tilde h^t_{\gamma,\bdg})_{(\bdy,0)},
(\tilde\omega_{\gamma,\bdg})_{(\bdy,0)},T)\right).
\end{split}
\end{equation}

In this case, once $\bdg\in  G_\gamma^t$ is fixed,
the
mapping  considered
$\pi_\gamma\circ \pi_{\bdg}$ is defined by
$x_l=t^{\mu_l^\bdg p_l} {\bar x}_l^{p_l}$ and $z=t^{c^\bdg}(z_1+\bar \beta)
\prod_{l=1}^d {\bar x}_l^{{\bar b}_l}.$
The pull-back
$\tilde h^t_{\gamma,\bdg}$ is
$$t^{nc^{\bdg} +\sum_{l=1}^{d-1}\mu_l^\bdg p_l N_l }{\bar x}_1^{\tilde N_1}\dots
{\bar x}_d^{\tilde N_d}(z_1^{m_r}+\ldots)
u({\bar \bdx},z_1), \quad u({\mathbf 0},0)\neq 0,
$$
and
$\tilde\omega_{\gamma,\bdg}=\bar \omega.$

Computations for (\ref{despBt})
and for all terms but the first one in (\ref{despB})
are similar.
In both cases, if either $\bdg\in G_\gamma^t$ or
$\bdg\in G_\gamma\setminus\{\bdg^1\}$  are fixed, then the
local
contribution from $Z^{\ww}_{DL}
(\tilde h^t_{\gamma,\bdg},\tilde\omega_{\gamma,\bdg},T)$
or from $Z^{\ww}_{DL}(\tilde h^t_{\bdg},\tilde\omega_{\bdg},T)$
we are interested in
is given by the $z_1$-highest vertex
${\tilde \tau}$ which corresponds in both cases
to the monomial
$$t^{\sum_{l=1}^{d-1}\mu_l^\bdg p_l N_l+ nc^{\bdg}}{\bar x}_1^{\tilde N_1}\dots
{\bar x}_d^{\tilde N_d}z_1^{m_r}.$$

We show in
Step 6 of the proof of theorems \ref{thm-polos} and
\ref{thm-top} that the factor coming from ${\tilde \tau}$
is nothing but
$T^{\sum_{l=1}^{d-1}\mu_l^\bdg p_l N_l+nc^{\bdg}}
\bl^{-(d+1)}(\bl-1)^{d+1}
\Phi_{\Delta_{{\tilde \tau}}}({\bdz}),$
where
$$
{\bdz}:=(
\bl^{-(\nu_1 p_1+\bar b_1)} T^{N_1 p_1+n\bar b_1},
\dots,\bl^{-(\nu_{d-1} p_{d-1}+\bar b_{d-1})} T^{N_{d-1}
p_{d-1}+n\bar b_{d-1}},\bl^{-(n_1+1)}
T^n,\bl^{-1} T^m).
$$
Since the corresponding cone
$\Delta_{{\tilde \tau}}^c$ is the positive cone $\bp^{d+1}$ (whose fundamental set
$\tilde G$ has only one element $\tilde G=\{\mathbf 1\}$),
its contribution is
$$
T^{\sum_{l=1}^{d-1}\mu_l^\bdg p_l N_l+nc^{\bdg}}
\bl^{-(d+1)}(\bl-1)^{d+1} S^1({\bdz}),
$$ 
where
\begin{equation}\label{contribu1}
S^1(\bdx):=\frac{\bdx^{\mathbf 1}}{\prod_{j=1}^{d+1}\left(1-x_j\right)}
=\frac{x_1x_2\dots x_{d+1}}{\prod_{j=1}^{d+1}\left(1-x_j\right)}.
\end{equation}

\medskip

Take any $\emptyset \subsetneq J \subsetneq \{1,\ldots,d-1\}.$
The number $e$ of essential variables of
$(\tilde h^t_{\gamma,\bdg})_{(\bdy,0)}$ and of
$(\tilde h^t_{\bdg})_{(\bdy,0)}$
is $e=d+1-\#J.$
By Corollary \ref{essen-var}, 
$
Z_{DL}^\ww((\tilde h^t_{\gamma,\bdg})_{(\bdy,0)},
(\tilde\omega_{\gamma,\bdg})_{(\bdy,0)},T)$ is
$\bl^{-\#J} (\bl-1)^e \bl^{-e} T^{\sum_{l=1}^{d-1}\mu_l^\bdg p_l N_l
+nc^{\bdg}}
S^1_J({\bdz})$
where
$
S^1_J(\bdx):=\frac{1}{(1-x_{d+1})}\prod_{j\not\in J}
\frac{x_j}{\left(1-x_j\right)}.
$
And the same identity is true for the zeta function $
Z_{DL}^\ww((\tilde h^t_{\bdg})_{(\bdy,0)},
(\tilde\omega_{\bdg})_{(\bdy,0)},T)$

Thus we get all  contributions we are interested in
for (\ref{despBt})
adding all previous results:
$$
\bl^{-(d+1)}(\bl-1)^{d+1}
\sum_{\bdg \in G_{\gamma}^t}
\bl^{-(\nu_\bdg+A_\bdg)}T^{\sum_{l=1}^{d-1}\mu_l^\bdg p_l N_l
+nc^{\bdg}}
 \left(S^1({\bdz})+
\sum_{\emptyset \subsetneq J \subset \{1,\ldots,d-1\}}
S^1_J({\bdz})\right).$$
It is proved by induction the formula
$
S^1(\bdx)+\sum_{\emptyset \subsetneq J \subset \{1,\ldots,d-1\}}
S^1_J(\bdx)
=S(\bdx)
$ where
$
S(\bdx):=\frac{x_d x_{d+1}}{\prod_{l=1}^{d+1} (1-x_l)}.$
Thus we must simplify the contribution
$$
B_1^t:=(1-\bl^{-(n_1+1)} T^n)
\left(\sum_{\bdg \in G^t_{\gamma}}
\bl^{-(\nu_\bdg+A_\bdg)}T^{\sum_{l=1}^{d-1}\mu_l^\bdg p_l N_l
+nc^{\bdg}} S({\bdz})\right)$$
under our hypothesis $\bl^{-(n_1+1)} T^n=1.$

In the other case, 
equation (\ref{despB}),
the contribution
of
$Z_{DL}({\bar h},{\bar \omega},T)$ 
comes  from the $z_1$-highest vertex
${\tilde \tau}$ which corresponds to the monomial
${\bar x}_1^{\tilde N_1}\dots
{\bar x}_d^{\tilde N_d}z_1^{m_r}$.
Thus it comes 
from $\bl^{-(d+1)}(\bl-1)^{d+1}
\Phi_{\Delta_{{\tilde \tau}}}({\bdz}),$
where
$(\bdz)$ was defined before.
In fact since the corresponding cone
$\Delta_{\tilde \tau}^c$ is the positive cone $\bp^{d+1}$
then the contribution is
$\bl^{-(d+1)}(\bl-1)^{d+1} S^1(\bdz)$ where
$
S^1(\bdx)$ (cf. equation (\ref{contribu1})).
Thus  all  contributions after the Newton maps
of (\ref{despB})  are
$$
\bl^{-(d+1)}(\bl-1)^{d+1} \left(S^1({\bdz}) +
\sum_{\bdg \in G_{\gamma},\bdg\ne\bdg^1}
\bl^{-(\nu_\bdg+A_\bdg)}T^{\sum_{l=1}^{d-1}\mu_l^\bdg p_l N_l
+nc^{\bdg}}
S({\bdz})\right).$$

We should simplify under our hypothesis $\bl^{-(n_1+1)} T^n=1$
the formula
$$
B_1:=(1-\bl^{-(n_1+1)} T^n) \left(S^1({\bdz})+
\sum_{\bdg \in G_{\gamma},\bdg\ne\bdg^1}
\bl^{-(\nu_\bdg+A_\bdg)}T^{\sum_{l=1}^{d-1}\mu_l^\bdg p_l N_l
+nc^{\bdg}}
S({\bdz})\right).
$$
Recall that the $(d+1)$-tuple $\bdz$
is nothing but $(w_1,\dots,w_{d-1},1,y_d)$.

Let us study the case $B_1^t$.
For every $\bdg$ parametrized by $(i_1,\dots,i_{d-1})$ from
equation (\ref{ggamma})
we have the following identity,
which is only valid after the simplification we are doing:
$$
\nu_\bdg+A_\bdg=
\frac{\alpha}{n_1}+\sum_{j=1}^{d-1}\frac{i_j}{p_j} (p_j\nu_j+{\bar b}_j)
-\frac{\displaystyle\alpha+\sum_{j=1}^{d-1} i_j b_j}{n_1}(n_1+1).
$$
It is turn out
$
\nu_\bdg+A_\bdg=
-\alpha-\sum_{j=1}^{d-1}i_j(b_j-\nu_j).
$
Since $c^{\bdg}=\frac{\alpha}{n_1}+\sum_{l=1}^{d}\mu_l \bar b_{l}$
then, under our
simplification, $c^{\bdg}=0$ and
the exponent of $T$ is nothing but $\sum_{j=1}^{d-1}i_j N_j$.

Therefore $B_1^t$ equals:
$$
\frac{\displaystyle\bl^{\alpha}z_{d+1}
\prod_{j=1}^{d-1}\sum_{k=1}^{p_j}(\bl^{(b_j-\nu_j)} T^{N_j})^k}
{\displaystyle(1-z_{d+1})\prod_{j=1}^{d-1}(1-z_j)}=
\frac{\displaystyle\bl^{\alpha} z_{d+1} \prod_{j=1}^{d-1}\left(u_j(1-u_j^{p_j})\right)}{\displaystyle(1-z_{d+1})
\prod_{j=1}^{d-1}\left((1-z_j)(1-u_j)\right)}.
$$

Written as before we have
\begin{equation}\label{cuarta}
B_1^t=\frac{\displaystyle a_d^{-\alpha} y_d \prod_{j=1}^{d-1} u_j}{\displaystyle(1-y_d)
\prod_{j=1}^{d-1}(1-u_j)}.
\end{equation}
Thus
$$
A_1^t+B_1^t=(\ref{tercera})-(\ref{primera0})-(\ref{segunda})+
(\ref{cuarta})=0.
$$

Let us study the case $B_1$. For every $\bdg\ne \bdg^1$ parametrized by $(i_1,\dots,i_{d-1})$ from
equation (\ref{ggamma}) but now
$i_j\in\{1,\ldots,p_j-1\},$
we have
$
\nu_\bdg+A_\bdg=-\sum_{j=1}^{d-1}i_j(b_j-\nu_j).
$
In the same way the exponent $T$ is nothing but
$\sum_{j=1}^{d-1}i_j N_j$.

The first term of $B_1$ is
\begin{equation}\label{cuarta1}
(1-z_d)S^1(\bdz)=\frac{y_d\displaystyle\prod_{j=1}^{d-1} w_j}
{\displaystyle(1-y_d)\prod_{j=1}^{d-1} (1-w_j)}=
\frac{\displaystyle y_d \prod_{j=1}^{d-1}u_j^{p_j}}
{\displaystyle(1-y_d)\prod_{j=1}^{d-1} (1-u_j^{p_j})}.
\end{equation}
The other terms comes from
$$
\frac{z_{d+1}}{(1-z_{d+1})
\displaystyle\prod_{j=1}^{d-1}(1-z_j)}\prod_{j=1}^{d-1}
\sum_{k=1}^{p_j-1}\left(\bl^{(b_j-\nu_j)} T^{N_j}\right)^k.
$$
This is nothing but
\begin{equation}\label{quinta1}
\frac{y_d}{\displaystyle(1-y_d)\prod_{j=1}^{d-1} (1-u_j^{p_j})}
\prod_{j=1}^{d-1}
\sum_{k=1}^{p_j-1} (u_j)^k=
\frac{y_d}
{\displaystyle(1-y_d)\prod_{j=1}^{d-1} (1-u_j^{p_j})}
\prod_{j=1}^{d-1} u_j \left(\frac{1-u_j^{p_j}}{1-u}-u_j^{p_j-1}\right).
\end{equation}
Thus
$
A_1+B_1=(\ref{tercera1})+
(\ref{cuarta1})+(\ref{quinta1})=0.
$

\end{proof}

\medskip

Let $h\in k[[\bdx]][z]$ be a quasi-ordinary power series
in good coordinates
with
$h(\bdx,z)=x_1^{N_1}x_2^{N_2}\dots x_d^{N_d}g(\bdx,z)$,
where no $x_i$ divides $g(\bdx,z)$
and $N_l\geq 0$ for any $l=1,\ldots,d.$
Write
$
h(\bdx,z)=x_1^{N_1}x_2^{N_2}\dots x_d^{N_d}f(\bdx,z)u(\bdx,z),
$
where $f(\bdx,z)$ is a quasi-ordinary $z$-polynomial of degree
$n$
in $k[[\bdx]][z].$ Let $\omega=\prod_{j=1}^d x_j^{\nu_j-1}
d x_1\wedge\ldots\wedge d x_d\wedge dz, \,\nu_j\geq 1$ define a form
such that $(h,\omega)$
verifies the support condition ~(\ref{spc}).
 The
set $CP(h,\omega)$ of \emph{candidate poles} was defined in  (\ref{candidate}).
Proposition  \ref{tspecial} will allow
to consider a smaller set of candidate poles than the set $CP(h,\omega)$.

\begin{defini} A
compact $1$-dimensional face $\gamma_q$
in $ND(h)$
will be called \emph{special} in the $i$-th coordinate if
the pair
$(h_{\gamma_q},\omega)$
verifies the conditions of Proposition \ref{tspecial}
in the $i$-th coordinate (instead of $d$-th coordinate).
 The face $\gamma_q$
will be called \emph{special} if it is special in one of the
coordinates.
If  $\gamma_q$
 is special in the $i$-th and $j$-th, ($i\ne j$), coordinates then the corresponding  candidate poles coincide.
Define
$$CP(h,\omega)_{q}:=
\left\{\left(\frac{M_l^q}{c_l^q},
\nu_lp^q_l+\bar{b}_l^q\right)
\right\}_{l=1,\ldots,\widehat i,\ldots,d}
\text{ if $\gamma_q$ is special,}
$$
where we assume that $\gamma_q$ is special in the $i$-th coordinate
and $\widehat i$ means that we omit $i.$
In particular,
if $\gamma_q$ is special in more than one variable then the corresponding pair 
$(N,\nu)$
appears in $CP(h,\omega)_{q}.$
Otherwise, define
$$
CP(h,\omega)_{q}:=
\left\{\left(\frac{M_l^q}{c_l^q},
\nu_lp^q_l+\bar{b}_l^q\right)\right\}_{l=1}^d
\text{ if $\gamma_q$ is not special. }
$$
Let us define by induction the sets
$$
{\widetilde {CP}}(h,\omega):=\bigcup_{q=1}^r CP(h,\omega)_{q} \cup
\bigcup\, {\widetilde {CP}}({\bar{h}_{q,j}},{\bar{\omega}_{q,j}})
$$
where the union is over all pull-back $(\bar{h}_{q,j},{\bar{\omega}_{q,j}})$ under the
Newton maps associated with Newton components $f_{j}^q$ of $f.$
The set of \emph{strong candidate poles}
of a pair $(h,\omega)$ is defined by
$$
SCP(h,\omega):=\left\{(N_i,\nu_i)\right\}_{i=1}^d
\cup {\widetilde {CP}}(h,\omega).$$
\end{defini}

After Proposition \ref{tspecial} and Theorem \ref{mainformula}
the following result is clear.

\begin{prop} If $N s+\nu$ is a pole of
$Z_{top,0}(h,\omega,s)$ then $(N,\nu)\in SCP(h,\omega)$.
\end{prop}

In fact Proposition \ref{tspecial} and equations (\ref{despues})
and (\ref{despues2}) show that the same result is true for the local
Denef-Loeser motivic zeta function.

\begin{prop}
$$
Z_{DL}(h,\omega,T)\in\bz[\bl,\bl^{-1},
(1-\bl^{-\nu}T^{N})^{-1}][T]_{(N,\nu)\in SCP(h,\omega)}.
$$
\end{prop}

\begin{obs}\label{specialproj} Let
$h=\prod_{l=1} x_l^{N_l}f(\bdx,z) u(\bdx,z),$
with $u(\mathbf 0,0)\ne 0.$ Take any compact $1$-dimensional
face $\gamma$ of $\Gamma(h)$.
For each $j\in\{1,\ldots,d\}$ such that
$\gamma$ is not contained in the $(x_j,z)$-plane, consider the
$j$-transversal section $h^{0}_{j}$
with root $\alpha=0$.  Let $\tilde \gamma$ be the compact face
of $\Gamma(h_j^0)$ on which $\gamma$ is projected.
The face $\tilde \gamma$ can be the projection
of several distinct compact faces of $\Gamma(h).$
Thus if $\gamma$ is not special
then its projection $\tilde \gamma$ cannot be special in any coordinate.
\end{obs}

\begin{stev1}
Assume  $h=x_1^{N_1}f(x_1,z)u(x_1,z)\in k[[x_1]][z]$
is a quasi-ordinary power series in good coordinates
with $\lgt(h)=1.$ Let $\gamma$ be a compact
$1$-dim face of $\Gamma(h).$
Then $\gamma$ is special if
and only if $\gamma$ is the highest (with respect to $z$)
compact face of $\Gamma(h)$ and
$h_\gamma$ is of type
$z^{n-n_1m}(z^{n_1}-\beta x_1)^{n-n_1},\, \beta\in {\mathbf G}_{m,k}$.
In particular $N_1=0$. We observe that $\gamma$ is special
if and only we can apply inversion formula to $h$
permuting coordinates. If
$\gamma$ is  special
then
$CP(h,\omega)_\gamma=\emptyset,$ for each pair
$(h,\omega)$ verifying the support condition (\ref{spc}).
After a Newton map associated with a Newton component of $f,$
the pull-back of $h$ is of type $y^a g(y,z_1)$ with $a>0.$
In particular Proposition \ref{tspecial} cannot be applied anymore.
\end{stev1}

If $\lgt(h)>1,$ more general facts occur.
Take for instance any of the following examples:  $f(x,y,z)=(z^2-x^3)^2+x^{11}y$
and  $g(x,y,u,z)=((z^2-x^3y)^2+x^7y^2)((z^ 1)2-x^3yu)^2+x^7y^2u^3).$
For $f$, after the unique possible Newton map $\pi,$
a special face for $f\circ \pi$ appears.
In the other case, $\Gamma(g)$ has two compact $1$-dimensional faces
and both are special, in different coordinates.

\begin{stev2}
Let $h\in k[[x_1,x_2]][z]$ be a quasi-ordinary
power series in good coordinates
with $\lgt(h)=2$.
Assume that we can decompose $h=x_1^{N_1}x_2^{N_2} f(x_1,x_2,z) u(x_1,x_2,z)$,
with $u(\mathbf 0)\ne 0$ and $f$ is a Weierstrass polynomial of degree
$n$. Let $\gamma$ denote the
highest (with respect to $z$) compact face of $\Gamma(h)$ and let
$h_\gamma=x_1^{N_1}x_2^{N_2}z^{n-n_1m}\prod_{j=1}^{u}
(z^{n_1}-\beta_j x_1^{b_1}x_2^{b_2})^{m_{j}}$ with
$\beta_j \in {\mathbf G}_{m,k}$ and $m=\sum m_j.$
Let $\omega=x_1^{\nu_1-1}x_2^{\nu_2-1}d x_1\wedge d x_2\wedge
dz$ be a differential form such that
$(h,\omega)$ verifies (\ref{spc}).
\end{stev2}

\begin{lema} \label{special2}
$(n,n_1+1)\in CP(h,w)_\gamma$ if and only if either
\begin{enumerate}
\item there exists $i\in\{1,2\}$ such that $b_i=1$, $N_i=0$ (this implies $\nu_i=1$)
and $u>1,$ or 
\item $u=1$ and $b_1=b_2=1$, $N_1=N_2=0$ (this implies $\nu_1=\nu_2=1$)
and therefore $\gamma$ is special in both coordinates.
\end{enumerate}
\end{lema}

\begin{proof} By definitions and Lemma \ref{1connuevo},
$\gamma$ is defined by
$n_1x_1+b_1z=n_1N_1+b_1n$ and $n_1x_2+b_2z=n_2N_2+b_2n$.
The candidate poles from $\gamma$ are $\left(\frac{n_1N_i+b_in}{\gcd(n_1,b_i)},
\frac{n_1\nu_i+b_i}{\gcd(n_1,b_i)}\right)$
with $i=1,2.$ In particular
$\left(\frac{n_1N_i+b_in}
{\gcd(n_1,b_i)},\frac{n_1\nu_i+b_i}{\gcd(n_1,b_i)}\right)=(n,n_1+1)$
if and only if either $b_i=0$ and then $\gcd(n_1,b_i)=n_1,$ which implies $(N_i,\nu_i)=(n,n_1+1)\in  CP(h,w)_\gamma$ which is absurd.
Or one of the two conditions of the lemma is verified.
\end{proof}

In case $(1)$, since $u>1$, $(n,n_1+1)$ is a strong candidate
pole for $h^{\alpha}_j$, $j\ne i).$ In case $(2)$
the $i$-transversal section $h^{\alpha}_i$ of $h$
has only $\alpha=0$ as a root because $(\frac{1}{n},\frac{1}{n})$ is the smallest
characteristic exponent of $f,$ it is the $z$-highest face of $NP(h).$
In particular $h_i^0$ has at least one Newton component
associated with the projection $\tilde \gamma$ of $\gamma$ on the plane $x_iz.$
Thus
$(n,n_1+1)$  is a candidate pole for
$(h_i^0,\omega_i).$ In fact, $\tilde \gamma$ is the highest face with respect to
$z$. Applying the case $\lgt(h)=1$ then $(n,n_1+1)$ is not a strong candidate pole
for both $h_i^0.$

\begin{prop} \label{long2poles} Under the above conditions one has
$$
SCP(h,\omega)\subset
\bigcup_{i=1}^2 \bigcup_{m=1}^{v_i}
SCP(h^{\alpha_m}_{i},\omega_i)\cup \{(n,n_1+1)\},
$$
where $h^{\alpha_m}_{i}$ are the corresponding $i$-transversal
section at the root $\alpha_m,$ cf. equation (\ref{cpp}) in Proposition \ref{transsection}.
\end{prop}

\begin{proof}
Let us prove by induction on $\dpt(h),$ that if we are in good coordinates,
with the above definitions of $n,$ $n_1$ and $\gamma$ then
$$
SCP(h,\omega)\subset
\bigcup_{i=1,2}\, \bigcup_{m=1}^{v_i}
SCP(h^{\alpha_m}_{i},\omega_i)\cup \{(n(h),n_1(h)+1)\}.
$$
We write $(n(h),n_1(h)+1)$
to remark that at each step of the induction they depend on $h$.

If $\dpt(h)=0$ then
$h=x_1^{N_1}x_2^{N_2} z\,u(\bdx,z),$ with $u(\mathbf 0)\ne 0$
The result follows
easily.
Assume that we have proved the result
for $\dpt(h)<m.$ Let us decompose
$h=x_1^{N_1}x_2^{N_2} f(x_1,x_2,z) u(x_1,x_2,z)$,
with $u(\mathbf 0)\ne 0$ and  $\dpt(h)=m.$
Let $(N,\nu)\in SCP(h,\omega) \subset CP(h,\omega).$
If $(N,\nu)=(n,n_1+1)$ we are done.
If $(N,\nu)=(N_i,\nu_i)\ne (0,1),$ $i=1,2$ then any of the
factors $h^{\alpha}_{j}$ of the
$j$-transversal section with $i\ne j$ has
$x_i^{N_i}$ as a factor which implies $(N,\nu)=(N_i,\nu_i)\in SCP(h^{\alpha}_{j},\omega_j).$

If $(N,\nu)\in\bigcup_{q=1}^r CP(h,\omega)_{q}$ where
$\Gamma(h)$ has $r$ compact $1$-dim faces. Assume $(N,\nu)$ appears
in the compact face $\gamma_q,$ say in the $i$-th coordinate.
Take the transversal section
$h_j^0,\, j\ne i$. Let $\tilde \gamma$ be the compact face
of $\Gamma(h_j^0)$ on which $\gamma_q$ is projected. 

If $\gamma_q$ is contained
in the $(x_i,z)$-plane,
then $\tilde \gamma$ contains $\gamma_q$. If we can apply proposition
\ref{tspecial} to $h_j^0$ and $\tilde \gamma$ then
$(h_j^0)|_{\tilde \gamma}=z^{n-mn_q}(z^{n_q}-\beta x_i)^m, \beta\in {\mathbf G}_{m,k_j}$. Then $\gamma_q$ must corresponds to the smallest 
characteristic exponent of $h$. Thus $\gamma_q$ is special in this coordinate. 
This contradicts 
$(N,\nu)\in SCP(h,\omega)$.

Otherwise, 
Remark \ref{specialproj}  contradicts
$(N,\nu)\in CP(h,\omega)_{q}$ if $\gamma_q$ is not special
in the $i$-th coordinate.
If $\gamma_q$ is special in the $i$-th coordinate
but not in the $j$-coordinate then
$(N,\nu)\not\in CP(h,\omega)_{q}$ which is absurd. Therefore if $\gamma_q$ is
special in both  coordinates  then $\gamma_q$ is the highest
compact face and by Lemma \ref{special2} $(2)$, $(N,\nu)=(n,n_1+1).$
Finally we apply induction.

To finish the proof of the proposition is enough
to remark that the conditions in Lemma \ref{special2} $(2)$
cannot be obtained after any Newton map. The reason for that
is that
after any Newton map
the pull-back of $h$ has at least a factor of type $y_i^{N_i}$
with $N_i>0.$ Therefore the proof is finished.
\end{proof}

\begin{stevd}
Let $h\in k[[x_1,\ldots,x_d]][z]$ be a quasi-ordinary
power series in good coordinates
with $\lgt(h)=d>2$.
Assume that $h=\prod_{l=1}^d x_l^{N_l} f(\bdx,z) u(\bdx,z),$
with $u(\mathbf 0)\ne 0$ and $f$ is a Weierstrass polynomial of degree
$n.$
Let $\omega$ be a differential form such that $(h,\omega)$ verifies
the support condition \ref{spc}.
\end{stevd}

\begin{prop} \label{longdpoles} Under the above conditions one has
$$
SCP(h,\omega)\subset
\bigcup_{i=1,\ldots,d,}\bigcup_{m=1}^{v_i}
SCP(h^{\alpha_m}_{i},\omega_i),
$$
where $h^{\alpha_m}_{i}$ are the corresponding $i$-transversal
section at the root $\alpha_m,$ see lemma~\ref{cpp}.
\end{prop}

\begin{proof} The proof is by induction on the
$\dpt(h)$ and follows the same ideas as proof of proposition
\ref{long2poles} but in this case we do not need to take care
of the special compact $1$-dimensional faces.
\end{proof}

\begin{anlt}
If $k=\bc$ and we work with
convergent  complex quasi-ordinary
power series $h\in \bc\{\bdx\}[z]$ all results presented in the last three
sections are valid too. We leave the details to the reader.
The pull-back under the Newton maps of convergent
quasi-ordinary power series are
again convergent.
The transversal sections can be seen now as follows.
We write
$h(\bdx,z)=\prod_{l=1}^d x_l^{N_l}f(\bdx,z)u(\bdx,z)$,
where no $x_i$ divides $f(\bdx,z)$, $u(0,0)\neq 0,$
$N_l\geq 0$ for any $l=1,\ldots,d,$
and $f$ is a Weierstrass polynomial in $z$.
After Lemma \ref{essen-var} we may assume that $\lgt(h)=d$.
We may assume  $f$ is convergent in a polydisk
$\Delta^d_\varepsilon\times\Delta^1_\delta$, where
$0\ll\delta\ll\varepsilon\ll 1$ and if $\bdx^0\in\Delta^d_\varepsilon$ all roots of
$f(\bdx^0,z)$ lie in $\Delta^1_\delta$.
\end{anlt}

Let us consider the germ $(\text{Sing}(V),0)$ of singular points of $V=h^{-1}(0)$
in a neighborhood of the origin.
The condition of the discriminant implies that $(\text{Sing}(V),0)$ is contained in the
intersection of $h^{-1}(0)$ with the hyperplanes
$x_l=0$, $1\leq l\leq d$. Fix $i\in\{1,\ldots,d\}$.
Consider the polynomial, over $\bc\{x_i\}$, of degree $n$
$$
f({\mathbf 0}_i(x_i),z)=z^n+\text{ lower degree terms}.
$$
where the $d$-tuple ${\mathbf 0}_i$
has all coordinates $0$ but the $i$-th coordinate which is $x_i$.
We have two possibilities:

\begin{enumerate}

\item If $f({\mathbf 0}_i(x_i),z)=z^n$ then
the points $({\mathbf 0}_i(x^0),0)$ belong
to $h^{-1}(0)$, for all $x^0\in\Delta^1_\varepsilon$.
The power series $h_{i}^{0,x^0}(\bdx,z):=h(\bdx+{\mathbf 0}_i(x^0),z)$
are quasi-ordinaries. They are the analytic
equivalent to
the formal $i$-transversal section at the root $\alpha=0.$

\item Otherwise, i.e. if $f({\mathbf 0}_i(x_i),z)\not=z^n,$
given $x^0\in\Delta^1_\varepsilon$ there
exist $z_1^0,\dots,z_l^0\in\Delta^1_\delta$,
pairwise distinct, and
$r_1,\dots,r_l\in\bp$, $\sum r_j=n$, such that
$$
f({\mathbf 0}_i(x^0),z)=
\prod_{j=1}^l(z-z_j^0)^{r_j}.
$$
We can suppose that $\varepsilon$ is small enough
in order to have that $r_1,\dots,r_l$
independent of $x^0$.
Again $h_{i}^{\alpha, x^0}(\bdx,z):=h(\bdx+{\mathbf 0}_i(x^0),z+z_j^0)$
are quasi-ordinary  power series.
In this case they are the equivalent in the analytic case to
the $i$-transversal at some root $\alpha\ne 0.$
\end{enumerate}

In particular we have the analytic analogue of proposition
\ref{transsection}.

\section{Monodromy conjecture for
quasi-ordinary power series}

Let $h:(\bc^{d+1},0)\to(\bc,0)$ be a germ of
complex analytic function such that $h(0)=0$.
Fix $U$ a sufficiently small neighbourhood
of $0$ where $h$ is defined.
Let $F$ be the Milnor fibre of the Milnor fibration
at the origin associated with $h$. Let $m_F:F\to F$
be the monodromy transformation. The zeta-function
of the monodromy of $h$ is
$$
\zeta(h)(t):=\prod_{q\geq 0}\det(I-t{m_F}_q)^{(-1)^q}
$$
where ${m_F}_q:
H_q(F,\bc)\to H_q(F,\bc),$ $q\geq 0,$ are the
homological monodromy transformations.

In this section we solved in the quasi-ordinary case
the 
monodromy conjectures stated in the Introduction.
We will focus firstly on the
the topological monodromy conjecture of \cite{dl:92}
and on the \emph{motivic monodromy conjecture} of \cite[section 2.4]{dl:98}
as stated in the introduction but for convergent quasi-ordinary series.


The monodromy conjectures deal with eigenvalues
of the complex algebraic monodromy at some points
of the zero locus of $h$.  First we will prove the monodromy
conjecture for $\lgt(h)=1.$ As we mention in the introduction the result is
known (several proofs by Loeser, Veys) but we present here an independent proof
following our ideas which by the way will be useful for the proof in the general case.
In the general case, i.e $\lgt(h)>1,$
we need to consider the monodromy at some different points
of the singular locus of $h.$
These points
will be at some transversal sections of $h$.
Thus we follow by induction on $\lgt(h).$

\subsection{Monodromy conjecture for curves}
\mbox{}

Recall that any germ $h\in\bc \{x,z\}$ of curve
is quasi-ordinary. 
We will prove that 
each Newton process produces eigenvalues of
the zeta function of the monodromy of the curve defined by $h$.
This is true in all cases except for
the \emph{special candidate pole} but the special pole
candidate (if it appears) is not a pole for the motivic zeta function (hence for the topological zeta function)
because of Proposition~\ref{tspecial}.

\begin{thm} \label{DLcurvas} Let $h(x,z)\in\bc\{x\}[z]$ be a germ of curve
and $\omega$ a regular differential form such that $(h,\omega)$
satisfies the support condition (\ref{spc}).
If $(N,\nu)\in SCP(h,\omega)$ let $q:=-\frac{\nu}{N}$, then 
 $\exp(2 i\pi q)$ is an eigenvalue
of the zeta function of the monodromy at some point of
$h^{-1}(0)$. 
\end{thm}


\begin{proof}
The case $\dpt(h)=0$ is trivial, since $h=x^N z$. In this case $SCP(h)=\{(N,\nu),(1,1)\}$,
where $\nu$ comes from the differential form; at some points
of $h^{-1}(0)$, the zeta function of the monodromy 
is either $1-t^N$ or $1-t$
and the result follows.

We  prove  by induction on the depth $\geq 1$ that if $(N,\nu)\in SCP(h,\omega)$ and $q:=-\frac{\nu}{N}$, then 
$\exp(2 i\pi q)$ is an eigenvalue
of the zeta function of the monodromy.  We will not consider the
case $q\in\bz$ since in this case it gives always the eigenvalue $1$ of the monodromy.

Let $h$ be a curve singularity with depth $m+1$ in good coordinates,
$m\geq 0$. Assume $h=x^Nz^\varepsilon g(x,z)$,
$\varepsilon=0,1$ and $\omega=x^{\nu-1} d x\wedge d z$, where neither $x$ nor $z$ divide $g$ and $(h,\omega)$ verifies the support condition \ref{spc}.
Let $\gamma_1,\dots,\gamma_r$ be the compact edges of the Newton diagram of $h$; for each $q=1,\dots,r,$ we have denoted by $v(q)$ the number of non-zero distinct roots of $h_{\gamma_q}$.
Each one of these roots define a Newton map and let $h_{q,j}$
$j=1,\dots,v(q)$, be the pull-backs of $h$ by these Newtons maps.
Recall that the vertex $\tau_r$ of $\gamma_r$ is the $z$-highest vertex.
We will define also
$\alpha,\beta$ such that $\alpha=0$ (resp. $\beta=0$) if the slope of $\gamma_r$ (resp. the inverse of the slope of $\gamma_1$) is an integer and $x$ (resp. $z$) does not divide $h$.
Otherwise we set $\alpha=1$ (resp. $\beta=1$).
Besides the candidate poles $(N,\nu)$ and $(1,1)$, denote
by $(N_q,\nu_q)$, $q=1,\dots,r$, the candidate poles corresponding to the edges. 
It is clear that with the poles $(N,\nu)$ and $(1,1)$ we can argue as above to conclude.
Note that:
\begin{itemize}
\item $\dpt(h)=1$ if and only if $\dpt(h_{q,j})=0$; the corresponding
monodromy zeta functions for these pull-backs are equal to $1$.

\item $h$ is in good coordinates if and only if $v(1)+\beta>1$.

\item $(N_r,\nu_r)$ is a strong candidate pole if and only if $v(r)+\alpha>1$; recall that $(N_q,\nu_q)$, $q=1,\dots,r-1$, are strong candidate poles.
\end{itemize}   
One can prove, see e.g. \cite{PiP}, that at the corresponding points of the strict transform of
$h$ after Newton mappings and after the partial resolution induced by the
Newton polygon of $h$ the total transforms of $h$ are isomorphic curves. 
Using A'Campo's formula
we deduce that the monodromy zeta function of $h$ at the origin is
\begin{equation}\label{monocur}
\zeta (h)(t)=
\frac{(1-t^{N'_r})^a}{(1-t^{N_r})^{\alpha-1}}
\frac{(1-t^{N''_1})^b}{(1-t^{N_1})^{\beta-1}}
\prod_{q=1}^{r}
\left(
(1-t^{N_q})^{-v(q)}
\prod_{j=1}^{v(q)}\zeta(h_{q,j})(t)
\right),
\end{equation}
where 
\begin{itemize}
\item $N'_1$ divides $N_1$, $N''_r$ divides $N_r$.

\item If $x$ (resp. $z$) divides $h$ then $a=0$ (resp. $b=0$)
and the corresponding numerator does not appear, otherwise $a=1$
(resp. $b=1$).

\item If $x$ does not divide $h$, it implies $(N,\nu)=(0,1)$,
 we can express $h_{\gamma_r}=z^{u_r}\tilde h_r$ where $\tilde h_r$ is a product
of $m_r=\sum_{j=1}^{v(r)} m_{r,j}$ factors of the form $(z^{n_1^r}-v x^{b^r})$, counting multiplicities. Recall that in the case of curves we have always $\gcd({n_1^r},b^r)=1$. Thus
\begin{equation}
\label{xnodiv}
N_r= (u_r+m_r n_1^r)b^r,\quad\nu_r=n_1^r+b^r,\quad N'_r=u_r+m_r n_1^r.
\end{equation}

\item If $z$ does not divide $h$, that is $\varepsilon=0$, we can express $h_{\gamma_1}=x^{u_1}\tilde h_1$  where $\tilde h_1$ is a product
of $m_1$ factors of the form $(z^{n_1^1}-\tilde v x^{b^1})$, again $\gcd(n_1^1,b^1)=1.$
In this case
\begin{equation}
\label{znodiv}
N_1= (u_1+m_1 b^1) n_1^1,\quad\nu_1=\nu n_1^1+b^1,\quad N''_1=u_1+m_1 b^1.
\end{equation}
\end{itemize} 
Recall that the inverse of the zeta function of the monodromy is, up to a factor $t$ and $(1-t)$ 
a polynomial. In order to prove the theorem we must take special care of numerators in the formula (\ref{monocur}).
They can arise in several situations:

\begin{enumerate}[(1)]

\item $r>1$, $x$ does not divide $h$, $v(r)=1$. 
The equations of (\ref{xnodiv}) imply that $\exp(-2 i\pi\frac{\nu_r}{N_r})$ is a root of the polynomial which is the inverse of the quotient $\frac{1-t^{N'_r}}{1-t^{N_r}}$.

\item $r>1$, $z$ does not divide $h$, $v(1)=1$. 
We proceed using equations (\ref{znodiv}).

\item $r=1$, either $x$ or $z$ divide $h$, $v(1)=1$. 
We proceed using equations (\ref{xnodiv}) or (\ref{znodiv}).

\item $r=1$, neither $x$ nor $z$ divide $h$, $v(1)=1$. 
Using equations (\ref{xnodiv}) it is easily seen that $(N_1,\nu_1)$
provides a root of the inverse of the monodromy zeta function. In such a case the zeta function is written as
\begin{equation}
\label{ultimocaso}
\frac{(1-t^{N'_1})(1-t^{N''_1})}{(1-t^{N_1})}\zeta(h_{1,1})(t).
\end{equation} 
\end{enumerate}

With this arguments, if $\dpt(h)=1$, the proof of the result is finished. 
Let us assume now that $\dpt(h)>1$. 
From the above computations we must prove that
$SCP(h_{q,j})$ produce also roots of $(\zeta(h)(t))^{-1}$.
From the induction hypothesis this is true whenever
$(\zeta(h_{q,j})(t))^{-1}$ divides $(\zeta(h)(t))^{-1}$. The unique case
where this fact does not happen is in (\ref{ultimocaso}), since
the first factor is not a polynomial. Let us rewrite the above formula
using (\ref{xnodiv}), (\ref{znodiv}) and $u_1=u_r=0$:
\begin{equation}
\label{ultimocaso1}
\frac{(1-t^{m_1 n_1^1})(1-t^{m_1 b^1})}{(1-t^{m_1 n_1^1 b^1})} \zeta(h_{1,1})(t).
\end{equation} 
Since the inverse zeta of the monodromy is, up to a factor $t(1-t),$ a polynomial,
it must exist a factor $(1-t^{m_1})$ which divides $(\zeta(h_{1,1})(t))^{-1}$. 
Let us note that $h_{1,1}=x^{m_1 n_1^1 b^1}(z^{m_1}+\dots)$ and the new differential form is
$x^{n_1^1+b^1-1}d x\wedge d z$. The highest vertex of the new Newton polygon is
$(m_1 n_1^1 b^1,m_1)$; let us assume that we have 
$s$~edges and denote $(\tilde N,\tilde\nu)$ 
the strong candidate pole associated to highest edge. 
Let us denote $(\tilde n,\tilde b)$ the coprime integers which provide the slope of the edge.

We have $\tilde N=m_1(\tilde n  n_1^1 b^1+\tilde b)$, $\tilde\nu=\nu_1\tilde n+\tilde b=(n_1^1+b^1)\tilde n+\tilde b$.
The function $\zeta(h_{1,1})(t)$ is 
\begin{equation}\label{monocur2}
\zeta (h_{1,1})(t)=(1-t^{\tilde N})^{-v(s)}
\zeta(\bar h_{s,j})(t)
\prod_{q=1}^{s-1}
\Delta_q,
\end{equation}
where $\Delta_q$ is the product of the factors  associated to the $s-1$ first edges.
 The factor $(1-t^{\tilde N})$ is a multiple of $(1-t^{m_1})$.
This means that $(\tilde N,\tilde\nu)$ is the only candidate pole of which may not give an eigenvalue of the monodromy of $h.$
Since 
$$
(\tilde N,\tilde\nu)=(m_1(\tilde n n_1^1 b^1+\tilde b),(n_1^1+b^1)\tilde n+\tilde b),
$$
it follows that
$\exp(-2 i\pi\frac{\tilde \nu}{\tilde N})^{m_1}=
\exp(-2 i\pi\frac{(n_1^1+b^1)\tilde n+\tilde b}{\tilde n n_1^1 b^1+\tilde b}).$
It easy to check that the rational number $\frac{(n_1^1+b^1)\tilde n+\tilde b}{\tilde n n_1^1 b^1+\tilde b}$
is not a positive integer. In particular $\exp(-2 i\pi\frac{\tilde \nu}{\tilde N})^{m_1}\neq 1$.
\end{proof}

\subsection{Monodromy conjecture: general case}
\mbox{}

We will need
some facts about the zeta-function of the monodromy
of quasi-ordinary analytic power series.
A formula to compute
the zeta function of the monodromy has been obtained
by L.J.~McEwan, A~N\'emethi, and P.D.~Gonz\'alez P\'erez,
see \cite{g:02,g:03}.

Assume $h$ is a quasi-ordinary power series
defined by  $f(\bdx)u(\bdx,z)$ with $u(0)\ne 0$ and
$f=z^n+a_1(\bdx)z^{n-1}+\ldots$ is in good coordinates.
If $\Lambda_{CE}=\emptyset$ then $f$ is irreducible and smooth.
Otherwise one
reorders the variables $\bdx$ in such a way
that the first entry of min $\Lambda_{CE}\in\bq^d$ is non-zero.  Since
 the set $\Lambda_{ND}(f)\subset \Lambda$ is totally ordered, 
min $\Lambda_{CE}$ exists.

\begin{thm} \cite{g:02} \cite{g:03} \label{nemethi}
With the above hypothesis,
after the reordering of the coordinates
described before
one has $\zeta(h)(t)=\zeta(f)(t)=\zeta\left(f|_{\{x_2=\ldots=x_d=0\}}\right)(t).$

In fact, if min $\Lambda_{CE}$ has at least two non-zero entries,
then $\zeta(f)(t)=(1-t^n)$.
\end{thm}

\begin{cor} If $h\in \bc\{\bdx\}[z]$ is quasi-ordinary 
satisfying $\lgt(h)=1$ 
and $\omega$ a regular differential form such that $(h,\omega)$
satisfies the support condition (\ref{spc}), then 
 $Z_{DL}(h,\omega,T)$
verifies the monodromy conjecture.
\end{cor}

We can apply Theorem \ref{nemethi} and Corollary \ref{essen-var}
to conclude the proof of the corollary
from Theorem \ref{DLcurvas}.

\begin{ejem}
If $f=z^n+x_1\cdots x_r$. Then $f$ has a
non-degenerated Newton polyhedron.
One can compute the topological zeta function
(using Denef and Loeser formula or our algorithm)
and the monodromy zeta function $\zeta_f(t)$
(using the above result or Varchenko formula, \cite{va:76}).
In fact $\zeta(f)=(1-t^n)$ and
$Z_{top,0}(f,s)$ has only two poles $s=-1$ and $s=-\frac{n+1}{n}.$
In particular the monodromy conjecture is verified.
\end{ejem}

\begin{thm} \label{maincoro} Let $h(x,z)\in\bc\{\bdx\}[z]$ be a quasi-ordinary power series
and $\omega$ a regular differential form such that $(h,\omega)$
satisfies the support condition \ref{spc}.
Then 
 $Z_{DL}(h,\omega,T)$ and $\zlo(h,\omega,s)$
verify the monodromy conjecture.
\end{thm}

\begin{proof}
Let $h\in\bc\{x_1,\dots,x_d\}[z]$ be a quasi-ordinary
analytic power series in good coordinates.
We proceed by induction on $\lgt(h)=d>1$.

Assume $h\in\bc\{x_1,x_2\}[z]$ has $\lgt(h)=2$.
By Proposition \ref{long2poles}, every strong candidate pole
is either $(n,n_1+1)$ or it is a strong candidate pole
of one of the transversal sections.
In the former case,
the highest compact $1$-dimensional face of $\Gamma(h)$
is special in both coordinates. Thus
 min $\Lambda_{CE}(h)$ has at least two non-zero entries,
which implies $\zeta(h)(t)=(1-t^n)$, after Theorem \ref{nemethi}.
In particular $\exp(-2i\pi (n_1+1)/n)$ is an eigenvalue
of the monodromy of $h$ at the origin.

Otherwise, since we have proved
the result for $\lgt(h)=1$, $\exp(-2i\pi (n_1+1)/n)$
is an eigenvalue of the monodromy
of one of the transversal sections
at some point which  is, in fact,
contained in $h^{-1}(0).$

Moreover for any quasi-ordinary analytic
power series $h$ such that $\lgt(h)=2$ the monodromy conjectures
for motivic and topological zeta functions
hold, after Corollary
\ref{essen-var}.

By Proposition \ref{longdpoles} and the induction
we can conclude that
the monodromy conjecture
for motivic and topological zeta functions
hold for quasi-ordinary power series $h$ such that $\lgt(h)=d$.
\end{proof}

\subsection{Monodromy conjecture for the Igusa zeta-function}
\mbox{}

Let $p$
be a prime number and let $K$ be a $p$-adic field, i.e
$[K:\bq_p]<\infty$.
Let $R$ be the valuation ring of $K$, $P$ the maximal ideal of $R$, and
$\bar K = R / P$ the residue field of $K$.
Let $q$ denote the cardinality of
$\bar K$, so $\bar K \simeq \ff_q$. For $z$ in $K$,  let
${\rm ord} \, z$ denote the valuation of $z$, and  set
$\vert z \vert = q^{- {\rm ord} \, z}$.
Let $h$ be a non constant element in $K [x_1, \ldots, x_d]$.
The $p$-adic \emph{Igusa local zeta function} $I_0(h,K,s)$ associated with $h$
(relative to
the trivial multiplicative character) is defined as the
$p$-adic integral
\begin{equation}
I_0(h,K,s) := \int_{PR^d} \vert h (x) \vert^s \vert d x \vert,
\end{equation}
for $s \in \bc$, ${\Re} (s) > 0$,
where $\vert d x \vert$ denotes the Haar measure on
$K^d$ normalized in such of way that $R^d$ is of volume 1.
Igusa proved that $I_0(h,K,s)$ is a rational function on $q^{-s},$
see \cite{de:91}.

\begin{obs} The method described here can be used to computed
$I_0(h,K,s)$ in case $h\in K [x_1, \ldots, x_d,z]$ be a quasi-ordinary
polynomial.
\end{obs}

Assume now that $h$ is a non-constant polynomial in $F [x_1, \ldots, x_d],$
for some number field $F\subset \bc.$ Igusa's monodromy conjecture
states that for almost all $p$-adic completion $K$ of $F,$ if $s_0$
is a pole of $I_0(h,K,s)$, then $\exp(2i\pi\Re(s_0))$
is an eigenvalue of the local monodromy of $h$ at some point of
$h^{-1}(0).$

It is known and deduced from \cite{dl:98}, see also \cite{dl:01},
that for almost all finite places
of the number field $F,$ the real parts $Ns+\nu$ of poles
of $I_0(h,K,s)$ come from factors $(1-\bl^{-\nu}T^N)$
in the denominator of $Z_{\text{naive}}(h,T)$
(which it is essentially nothing but $Z_{\text{DL}}(h,T)).$
In particular our proof can be applied to
polynomials $h\in F[\bdx,z]$ with coefficients in a number field
$F$ which are quasi-ordinary polynomials.
Then  the Igusa monodromy conjecture
is also true in this case.

\providecommand{\bysame}{\leavevmode\hbox to3em{\hrulefill}\thinspace}

\providecommand{\MR}{\relax\ifhmode\unskip\space\fi MR }
\providecommand{\MRhref}[2]{%
  \href{http://www.ams.org/mathscinet-getitem?mr=#1}{#2}
}
\providecommand{\href}[2]{#2}

\end{document}